\documentclass[11pt]{amsart}
\usepackage{amsmath,amsfonts,amscd,amssymb,epsf}
\newtheorem{theorem}{Theorem}[section]
\newtheorem{proposition}[theorem]{Proposition}
\newtheorem{lemma}[theorem]{Lemma}
\newtheorem{corollary}[theorem]{Corollary}
\theoremstyle{definition}
\newtheorem{definition}[theorem]{Definition}

\theoremstyle{remark} \newtheorem{remark}[theorem]{Remark}
\numberwithin{equation}{section}
\newcommand{\field}[1]{\ensuremath{\mathbb{#1}}}
\newcommand{\CC}{\field{C}}

\newcommand{\RR}{\field{R}}

\newcommand{\ZZ}{\field{Z}}
\newcommand{\complex}[1]{\mathsf{#1}} 

\newcommand{\SSS}{\complex{S}}
\newcommand{\BBB}{\complex{B}}
\newcommand{\KKK}{\complex{K}}
\newcommand{\AAA}{\complex{A}}
\newcommand{\CCC}{\complex{C}}

\DeclareMathOperator{\id}{id} 
\DeclareMathOperator{\Tot}{Tot} \DeclareMathOperator{\Hom}{Hom}
 
 \DeclareMathOperator{\Diff}{Diff}
\DeclareMathOperator{\Mob}{M\ddot{o}b}
\DeclareMathOperator{\im}{Im}
 \DeclareMathOperator{\PSL}{PSL}
\DeclareMathOperator{\re}{Re}
\newcommand{\del}{\partial}
\newcommand{\delb}{\bar\partial}

\newcommand{\ptheta}{\check{\theta}}
\newcommand{\pu}{\check{u}}

\newcommand{\T}{\mathfrak{T}}
\newcommand{\Pro}{\mathfrak{P}}
\newcommand{\D}{\mathfrak{D}}
\newcommand{\ihalf}{\frac{i}{2}}

\newcommand{\up}{{\mathbb{U}}}
\newcommand{\lo}{{\mathbb{L}}}
\newcommand{\Z}{{\mathbb{Z}}}
\newcommand{\U}{{\mathbb{U}}}
\newcommand{\C}{{\mathbb{C}}}
\newcommand{\al}{\alpha}
\newcommand{\be}{\beta}
\newcommand{\de}{\delta}
\newcommand{\g}{\gamma}
\newcommand{\s}{\gamma}

\newcommand{\N}{\mathbb{N}}
\newcommand{\cl}{\overline}
\newcommand{\bk}{\backslash}
\newcommand{\Ga}{\Gamma}

\newcommand{\pa}{\partial}
\newcommand{\la}{\langle}
\newcommand{\ra}{\rangle}

\newcommand{\Ka}{K\"{a}hler\:}
\newcommand{\Te}{Teichm\"{u}ller\:}

\newcommand{\bu}{\bullet}
\newcommand{\ov}{\overline}
\newcommand{\ep}{\epsilon}
\newcommand{\vep}{\varepsilon}

\newcommand{\z}{\bar{z}}

\newcommand{\ma}[4]{(\begin{smallmatrix}
              #1 & #2 \\ #3 & #4
             \end{smallmatrix})}

\begin{document}
\title[Liouville action and Weil-Petersson metric]{Liouville action and
Weil-Petersson metric on deformation spaces, global Kleinian reciprocity and
holography}
\author{Leon A. Takhtajan} \address{Department of Mathematics \\
SUNY at stony Brook \\ Stony Brook, NY 11794-3651 \\ USA}
\email{leontak@math.sunysb.edu}
\author{Lee-Peng Teo} \address{Department of Mathematics \\
SUNY at stony Brook \\ Stony Brook, NY 11794-3651 \\ USA}
\email{lpteo@math.sunysb.edu}
\begin{abstract}
We rigorously define the Liouville action functional for finitely generated, purely
loxodromic quasi-Fuchsian group using homology and cohomology double complexes
naturally associated with the group action. We prove that the classical action --- the
critical  value of the Liouville action functional, considered as a function on the
quasi-Fuchsian deformation space, is an antiderivative of a 1-form given by the
difference of Fuchsian and quasi-Fuchsian projective connections. This result can be
considered as global quasi-Fuchsian reciprocity which implies McMullen's
quasi-Fuchsian reciprocity. We prove that the classical action is a \Ka potential of
the Weil-Petersson metric. We also prove that Liouville action functional satisfies
holography principle, i.e., it is a regularized limit of the hyperbolic volume of a
3-manifold associated with a quasi-Fuchsian group. We generalize these results to a
large class of Kleinian groups including finitely generated, purely loxodromic
Schottky and quasi-Fuchsian groups, and their free combinations.
\end{abstract}
\maketitle \tableofcontents
\section{Introduction}

Fuchsian uniformization of Riemann surfaces plays an important role in the \Te theory.
In particular, it is built into the definition of the Weil-Petersson metric on \Te
spaces. This role became even more prominent with the advent of the string theory,
started with Polyakov's approach to non-critical bosonic strings~\cite{Pol}. It is
natural to consider the hyperbolic metric on a given Riemann surface as a critical
point of a certain functional defined on the space of all smooth conformal metrics on
it. In string theory this functional is called \emph{Liouville action functional} and
its critical value --- \emph{classical action}. This functional defines
\emph{two-dimensional theory of gravity with cosmological term} on a Riemann surface,
also known as \emph{Liouville theory}.

From a mathematical point of view, relation between Liouville theory and complex
geometry of moduli spaces of Riemann surfaces was established by P.~Zograf and the
first author in~\cite{ZT85, ZT87a, ZT87b}. It was proved that the classical action is
a K\"{a}hler potential of the Weil-Petersson metric on moduli spaces of pointed
rational curves~\cite{ZT87a}, and on Schottky spaces~\cite{ZT87b}. In the rational
case the classical action is a generating function of accessory parameters of Klein
and Poincar\'{e}. In the case of compact Riemann surfaces, the classical action is an
antiderivative of a 1-form on the Schottky space given  by the difference of Fuchsian
and Schottky projective connections. In turn, this 1-form is an antiderivative of the
Weil-Petersson symplectic form on the Schottky space.

C.~McMullen~\cite{McM} has considered another 1-form on \Te space given  by the
difference of Fuchsian and quasi-Fuchsian projective connections, the latter
corresponds to Bers' simultaneous uniformization of a pair of Riemann surfaces. By
establishing quasi-Fuchsian reciprocity, McMullen proved that this 1-form is also an
antiderivative of the Weil-Petersson symplectic form, which is bounded on the \Te
space due to the Kraus-Nehari inequality. The latter result is important in proving
that moduli space of Riemann surfaces is K\"{a}hler hyperbolic~\cite{McM}.

In this paper we extend McMullen's results along the lines of~\cite{ZT87a, ZT87b} by
using homological algebra machinery developed by E.~Aldrovandi and the first author
in~\cite{AT}. We explicitly construct a smooth function on quasi-Fuchsian deformation
space and prove that it is an antiderivative of the 1-form given by the difference of
Fuchsian and quasi-Fuchsian projective connections. This function is defined as a
classical action for Liouville theory for a quasi-Fuchsian group. The symmetry
property of this function is the global quasi-Fuchsian reciprocity, and McMullen's
quasi-Fuchsian reciprocity~\cite{McM} is its immediate corollary. We also prove that
this function is a K\"{a}hler potential of the Weil-Petersson metric on the
quasi-Fuchsian deformation space. As it will be explained below, construction of the
Liouville action functional is not a trivial issue and it requires homological algebra
methods developed in~\cite{AT}. Furthermore, we show that the Liouville action
functional satisfies \emph{holography principle} in string theory (also called
\emph{AdS/CFT correspondence}). Specifically, we prove that the Liouville action
functional is a regularized limit of the hyperbolic volume of a 3-manifold associated
with a quasi-Fuchsian group. Finally, we generalize these results to a large class of
Kleinian groups including finitely generated, purely loxodromic Schottky and
quasi-Fuchsian groups, and their free combinations. Namely, we define the Liouville
action functional, establish the holography principle, and prove that the classical
action is an antiderivative of a 1-form on the deformation space given by the
difference of Fuchsian and Kleinian projective connections, thus establishing global
Kleinian reciprocity. We also prove that the classical action is a \Ka potential of
the Weil-Petersson metric.

Here is a more detailed description of the paper. Let $X$ be a Riemann surface of
genus $g>1$,  and let $\{U_{\alpha}\}_{\alpha\in A}$ be its open cover with charts
$U_{\al}$, local coordinates $z_\al:U_\al\rightarrow \CC$, and transition functions
$f_{\al\be}: U_\al\cap U_\be\rightarrow\CC$.  A (holomorphic) projective connection on
$X$ is a collection $P=\{p_\al\}_{\al\in A}$, where $p_\al$ are holomorphic functions
on $U_\al$ which on every $U_\al\cap U_\be$ satisfy
\begin{equation*}
p_\be = p_\al\circ f_{\al\be} \left(f_{\al\be}'\right)^2 + \mathcal{S}(f_{\al\be}),
\end{equation*}
where prime indicates derivative. Here $\mathcal{S}(f)$ is the Schwarzian derivative,
\begin{equation*}
\mathcal{S}(f)=\frac{f'''}{f'} - \frac{3}{2}\left(\frac{f''}{f'}\right)^2.
\end{equation*}
The space $\mathcal{P}(X)$ of projective connections on $X$ is an affine space modeled
on the vector space of holomorphic quadratic differentials on $X$.

The Schwarzian derivative satisfies the following properties.
\begin{itemize}
\item[\textbf{SD1}] $\mathcal{S}(f \circ g) = \mathcal{S}(f) \circ g\;(g')^2 +
\mathcal{S}(g)$.
\item[\textbf{SD2}]
$\mathcal{S}(\s) = 0~\text{for all}~\s\in\PSL(2,\CC)$.
\end{itemize}
Let $\pi: \Omega\rightarrow X$ be a holomorphic covering of a compact Riemann
surface $X$ by a domain $\Omega\subset\hat\CC$ with a group of deck transformations 
being a subgroup of $\PSL(2,\CC)$. It follows from \textbf{SD1-SD2} that 
every such covering defines a projective connection on $X$ by 
$P_\pi =\{\mathcal{S}_{z_\al}(\pi^{-1}) \}_{\al\in
A}$. The Fuchsian uniformization $X\simeq \Ga\bk\up$ is the covering $\pi_F: \up
\rightarrow X$ by the upper half-plane $\up$ where the group of deck transformations
is a Fuchsian group $\Ga$, and it defines Fuchsian projective connection $P_F$. The
Schottky uniformization $X\simeq\Ga\bk\Omega$ is the covering $\pi_S:\Omega
\rightarrow X$ by a connected domain $\Omega\subset\hat{\CC}$ where the group of deck
transformations $\Ga$ is a Schottky group --- finitely-generated, strictly loxodromic,
free Kleinian group. It defines Schottky projective connection $P_S$.

Let $\mathfrak{T}_g$ be the \Te space of marked Riemann surfaces of genus $g>1$ (with
a given marked Riemann surface as the origin), defined as the space of marked
normalized Fuchsian groups, and let $\mathfrak{S}_g$ be the Schottky space, defined as
the space of marked normalized Schottky groups with $g$ free generators. These spaces
are complex manifolds of dimension $3g-3$ carrying Weil-Petersson K\"{a}hler metrics,
and the natural projection map $\T_g\rightarrow \mathfrak{S}_g$ is a complex-analytic
covering. Denote by $\omega_{WP}$ the symplectic form of the Weil-Petersson metric on
spaces $\T_g$ and $\mathfrak{S}_g$, and by $d=\pa + \bar{\pa}$ --- the de Rham
differential and its decomposition. The affine spaces $\mathcal{P}(X)$ for varying
Riemann surfaces $X$ glue together to an affine bundle $\mathfrak{P}_g\rightarrow
\T_g$, modeled over holomorphic cotangent bundle of $\T_g$. The Fuchsian projective
connection $P_F$ is a canonical section of the affine bundle
$\mathfrak{P}_g\rightarrow \T_g$, the Schottky projective connection is a canonical
section of the affine bundle $\mathfrak{P}_g\rightarrow \mathfrak{S}_g$, and their
difference $P_F - P_S$ is a $(1,0)$-form on $\mathfrak{S}_g$. This 1-form has the
following properties~\cite{ZT87b}. First, it is $\pa$-exact --- there exists a smooth
function $S: \mathfrak{S}_g\rightarrow \RR$ such that
\begin{equation} \label{I}
P_F - P_S=\frac{1}{2}\,\del S.
\end{equation}
Second, it is a $\bar\pa$-antiderivative, and hence a $d$-antiderivative by \eqref{I},
of the Weil-Petersson symplectic form  on $\mathfrak{S}_g$
\begin{equation} \label{II}
\delb(P_F - P_S)= - i\,\omega_{WP}.
\end{equation}
It immediately follows from \eqref{I} and \eqref{II} that the function $-S$ is a
K\"{a}hler potential for the Weil-Petersson metric on $\frak{S}_g$, and hence on
$\frak{T}_g$,
\begin{equation} \label{Kahler}
\del\delb S = 2i\,\omega_{WP}.
\end{equation}

Arguments using quantum Liouville theory (see, e.g., \cite{T} and references therein)
confirm formula \eqref{I} with function $S$ given by the classical Liouville action,
as was already proved in~\cite{ZT87b}. However, general mathematical definition of the
Liouville action functional on a Riemann surface $X$ is a non-trivial problem
interesting in its own right (and for rigorous applications to quantum Liouville
theory). Let $\mathcal{CM}(X)$ be a space (actually a cone) of smooth conformal
metrics on a Riemann surface $X$. Every $ds^2\in \mathcal{CM}(X)$ is a collection
$\left\{e^{\phi_\alpha}|dz_\alpha|^2\right\}_{\al\in A}$, where functions
$\phi_\alpha\in C^\infty(U_\al,\RR)$ satisfy
\begin{equation} \label{glue}
\phi_\alpha\circ f_{\alpha\beta} + \log |f'_{\alpha\beta}|^2
=\phi_\beta\quad\text{on}\quad U_\al\cap U_\be.
\end{equation}
According to the uniformization theorem, $X$ has a unique conformal metric of constant
negative curvature $-1$, called hyperbolic, or Poincar\'{e} metric. Gaussian curvature
$-1$ condition is equivalent to the following nonlinear PDE for functions $\phi_\al$
on $U_\al$,
\begin{equation} \label{Liouville}
\frac{\pa^2\phi_\al}{\pa z_\al \pa\bar{z}_\al}=\frac{1}{2}\,e^{\phi_\al}.
\end{equation}
In the string theory this PDE is called the Liouville equation. The problem is to
define Liouville action functional on Riemann surface $X$ --- a smooth functional
$S: \mathcal{CM}(X)\rightarrow\RR$ such that its Euler-Lagrange equation is the
Liouville equation. At first glance it looks like an easy task. Set $U=U_\al,\,
z=z_\al$ and $\phi=\phi_\al$, so that $ds^2=e^\phi|dz|^2$ in $U$. Elementary calculus
of variations shows that the Euler-Lagrange equation for the functional
\begin{equation*}
\frac{i}{2}\iint\limits_U \left(|\phi_z|^2 + e^\phi\right)dz\wedge d\z,
\end{equation*}
where $\phi_z = \pa\phi/\pa z$, is indeed the Liouville equation on $U$. Therefore, it
seems that the functional $\frac{i}{2}\iint_X \omega$, where $\omega$ is a 2-form on
$X$ such that
\begin{equation} \label{form}
\left.\omega\right|_{U_\al} =
\omega_\alpha=\left(\left|\frac{\partial\phi_\alpha}{\partial z_\alpha}\right|^2 +
e^{\phi_\alpha}\right)dz_\alpha\wedge d\bar{z}_\alpha,
\end{equation}
does the job. However, due to the transformation law \eqref{glue} the first terms in
local 2-forms $\omega_\alpha$ do not glue properly on $U_\al\cap U_\be$ and a 2-form
$\omega$ on $X$ satisfying \eqref{form} does not exist!

Though the Liouville action functional can not be defined in terms of a Riemann
surface $X$ only, it can be defined in terms of planar coverings of $X$. Namely, let
$\Ga$ be a Kleinian group with region of discontinuity $\Omega$ such that
$\Ga\bk\Omega\simeq X_1\sqcup\dots\sqcup X_n$ --- a disjoint union of compact Riemann
surfaces of genera $>1$ including Riemann surface $X$. The covering $\Omega\rightarrow
X_1\sqcup\dots\sqcup X_n$ introduces a global ``\'{e}tale'' coordinate, and for large
variety of Kleinian groups (Class $A$ defined below) it is possible, using
methods~\cite{AT}, to define a Liouville action functional $S:
\mathcal{CM}(X_1\sqcup\dots\sqcup X_n)\rightarrow\RR$ such that its critical value is
a well-defined function on the deformation space $\D(\Ga)$. In the simplest case when
$X$ is a punctured Riemann sphere such global coordinate exists already on $X$, and
Liouville action functional is just $\frac{i}{2}\iint_X \omega$, appropriately
regularized at the punctures ~\cite{ZT87a}. When $X$ is compact, one possibility is to
use the ``minimal'' planar cover of $X$ given by the Schottky uniformization $X\simeq
\Ga\bk\Omega$, as in~\cite{ZT87b}. Namely, identify $\mathcal{CM}(X)$ with the affine
space of smooth real-valued functions $\phi$ on $\Omega$ satisfying
\begin{equation} \label{field-0}
\phi\circ\s + \log|\s'|^2 =\phi\quad\text{for all}\quad\s\in\Ga,
\end{equation}
and consider the 2-form $\omega[\phi]=(|\phi_z|^2 + e^\phi)dz\wedge d\z$ on $\Omega$.
The 2-form $\omega[\phi]$ can not be pushed forward on $X$, so that the integral
$\frac{i}{2}\iint_F\omega$ depends on the choice of a fundamental domain $F$ for a
marked Schottky group $\Ga$. However, one can add boundary terms to this integral to
ensure the independence of the choice of a fundamental domain 
for a marked Schottky group $\Ga$, and to guarantee that
its Euler-Lagrange equation is the Liouville equation on $\Ga\bk\Omega$. The result is
the following functional introduced in~\cite{ZT87b}
\begin{align} \label{Schottky}
S[\phi] = & \frac{i}{2}\iint\limits_F\left(|\phi_z|^2 + e^\phi\right) dz\wedge d\z \\
          & + \frac{i}{2}\sum_{k=1}^g \int_{C_k}
          \left(\phi - \frac{1}{2}\log|\s'_k|^2
          \right)\left(\frac{\s''_k}{\s'_k}dz - \frac{\ov{\s''_k}}{\ov{\s'_k}}d\z\right)
          \nonumber \\
          & + 4\pi\sum_{k=1}^g \log|c(\s_k)|^2. \nonumber
\end{align}
Here $F$ is the fundamental domain of the marked Schottky group $\Ga$ with free
generators $\s_1,\dots,\s_g$, bounded by $2g$ nonintersecting closed Jordan curves
$C_1,\dots, C_g, C'_1,\dots,C'_g$ such that $C'_k = -\s_k(C_k)$, $k=1,\dots,g$, and
$c(\s)=c$ for $\s=\left(\begin{smallmatrix} a & b \\ c & d \end{smallmatrix}\right)$.
Classical action $S:\mathfrak{S}_g\rightarrow\RR$ that enters \eqref{I} is the
critical value of this functional.

In \cite{McM} McMullen considered quasi-Fuchsian projective connection $P_{QF}$ on a
Riemann surface $X$ which is given by Bers' simultaneous uniformization of $X$ and a
fixed Riemann surface $Y$ of the same genus and opposite orientation. Similar to
formula \eqref{II}, he proved
\begin{equation} \label{McM}
d(P_F - P_{QF}) = -i\,\omega_{WP},
\end{equation}
so that the 1-form $P_F - P_{QF}$ on $\mathfrak{T}_g$ is a $d$-antiderivative of the
Weil-Petersson symplectic form, bounded in \Te and Weil-Petersson metrics due to
Kraus-Nehari inequality. Part $\bar\pa(P_F - P_{QF}) = -i\,\omega_{WP}$ of
\eqref{McM} actually follows from \eqref{I} since $P_S - P_{QF}$ is holomorphic
$(1,0)$-form on $\mathfrak{S}_g$. Part $\pa(P_F - P_{QF})=0$ follows from
McMullen's quasi-Fuchsian reciprocity.

Our first result is the analog of the formula \eqref{I} for the quasi-Fuchsian case,
giving the 1-form $P_F - P_{QF}$ the same treatment as to the 1-form $P_F - P_S$.
Namely, let $\Ga$ be a finitely generated, purely loxodromic quasi-Fuchsian group with
region of discontinuity $\Omega$, so that $\Ga\bk\Omega$ is the disjoint union of two
compact Riemann surfaces with the same genus $g>1$ and opposite orientations. Denote
by $\D(\Ga)$ the deformation space of $\Ga$ --- a complex manifold of complex
dimension $6g-6$, and by $\omega_{WP}$ --- the symplectic form of the Weil-Petersson
metric on $\D(\Ga)$.  To every point $\Ga'\in\D(\Ga)$ with the region of discontinuity
$\Omega'$ there corresponds a pair $X, Y$ of compact Riemann surfaces with opposite
orientations simultaneously uniformized by $\Ga'$, that is, $X\sqcup Y\simeq
\Ga'\bk\Omega'$. We will continue to denote by $P_F$ and $P_{QF}$ projective
connections on $X\sqcup Y$ given by Fuchsian uniformizations of $X$ and $Y$ and Bers'
simultaneous uniformization of $X$ and $Y$ respectively. Similarly to \eqref{I}, we
prove in Theorem \ref{variation1} that there exists a smooth function $S:
\D(\Ga)\rightarrow\RR$ such that
\begin{equation} \label{F-QF}
P_F - P_{QF}=\frac{1}{2}\,\pa S.
\end{equation}

The function $S$ is Liouville classical action for the quasi-Fuchsian group $\Ga$
---- the critical value of the Liouville action functional $S$ on
$\mathcal{CM}(X\sqcup Y)$. Its construction uses double homology and cohomology
complexes naturally associated with the $\Ga$-action on $\Omega$. Namely, the homology
double complex $\KKK_{\bu,\bu}$ is defined as a tensor product over the integral group
ring $\Z\Ga$ of the standard singular chain complex of $\Omega$ and the canonical
bar-resolution complex for $\Ga$, and cohomology double complex $\CCC^{\bu,\bu}$ is
bar-de Rham complex on $\Omega$. The cohomology construction starts with the 2-form
$\omega[\phi]\in\CCC^{2,0}$, where $\phi$ satisfies \eqref{field-0}, and introduces
$\theta[\phi]\in\CCC^{1,1}$ and $u\in\CCC^{1,2}$ by
\begin{equation*}
\theta_{\s^{-1}}[\phi] = \left(\phi -  \frac{1}{2}\log|\s'|^2\right)
\left(\frac{\s''}{\s'} dz - \frac{\ov{\s''}}{\ov{\s'}}d\z\right),
\end{equation*}
and
\begin{align*}
u_{\s_1^{-1},\s_2^{-1}}= & -\frac{1}{2}\log|\s_1'|^2
\left(\frac{\s_2''}{\s_2'}\circ\s_1\, \s_1'\, dz -
\frac{\ov{\s_2''}}{\ov{\s_2'}}\circ\s_1\, \ov{\s_1'}\,d\z\right) \\
 & + \frac{1}{2}\log|\s_2'\circ\s_1|^2\left(\frac{\s_1''}{\s_1'} dz -
\frac{\ov{\s_1''}}{\ov{\s_1'}}d\z\right).
\end{align*}
Define $\Theta\in\CCC^{0,2}$ to be a group 2-cocycle satisfying $d\Theta=u$. The
resulting cochain $\Psi[\phi] = \omega[\phi] - \theta[\phi] - \Theta$ is a cocycle of
degree 2 in the total complex $\Tot\CCC$. Corresponding homology construction starts
with fundamental domain $F\in\KKK_{2,0}$ for $\Ga$ in $\Omega$ and introduces chains
$L\in\KKK_{1,1}$ and $V\in\KKK_{0,2}$ such that $\Sigma = F + L - V$ is a cycle of
degree 2 in the total homology complex $\Tot\KKK$. The Liouville action functional is
given by the evaluation map,
\begin{equation} \label{Evaluation}
S[\phi] = \frac{i}{2}\left\la\Psi[\phi],\Sigma\right\ra,
\end{equation}
where $\la~,~\ra$ is the natural pairing between $\CCC^{p,q}$ and $\KKK_{p,q}$.

In case when $\Ga$ is a Fuchsian group, the Liouville action functional on
$X\simeq \Ga\bk\U$, similar to \eqref{Schottky}, can be written explicitly as follows
\begin{align*}
S[\phi] = & \ihalf\iint\limits_F\omega[\phi] + \ihalf\sum_{k=1}^g\left(
\int_{a_k}\theta_{\al_k}[\phi] - \int_{b_k}\theta_{\be_k}[\phi]\right) \\ & +
\ihalf\sum_{k=1}^g\left(\Theta_{\al_k,\be_k}(a_k(0)) - \Theta_{\be_k,\al_k}(b_k(0)) +
\Theta_{\g_k^{-1},\al_k\be_k}(b_k(0))\right) \\ & -
\ihalf\sum_{k=1}^g\Theta_{\g_g^{-1}\ldots\g_{k+1}^{-1},\g_k^{-1}} (b_g(0)),
\end{align*}
where
\begin{equation*}
\Theta_{\s_1,\s_2}(z)=\int_{p}^{z}u_{\s_1,\s_2} + 4\pi i\varepsilon_{\s_1,\s_2}(2\log
2+ \log|c(\s_2)|^2),
\end{equation*}
$p\in\RR\setminus\Ga(\infty)$ and
\begin{equation*}
     \varepsilon_{\s_1,\s_2} = \begin{cases}
                              \;\;\,1 & \text{if $p < \s_2(\infty) < \s_1^{-1} p$}, \\
                             -1 & \text{if $p > \s_2 (\infty) > \s_1^{-1} p$},\\
                              \;\;\,0 & \text{otherwise}.
                             \end{cases}
\end{equation*}
Here $a_k$ and $b_k$ are edges of the fundamental domain $F$ for $\Ga$ in $\U$ (see
Section 2.2.1) with initial points $a_k(0)$ and $b_k(0)$, $\al_k$ and $\be_k$ are
corresponding generators of $\Ga$ and $\s_k =\al_k\be_k \al^{-1}_k\be^{-1}_k$. The
action functional does not depend on the choice of the fundamental domain $F$ for
$\Ga$, nor on the choice of $p\in\RR\setminus\Ga(\infty)$. Liouville action for
quasi-Fuchsian group $\Ga$ is defined by a similar construction where both components
of $\Omega$ are used (see Section 2.3.3).

Equation \eqref{F-QF} is global quasi-Fuchsian reciprocity. McMullen's
quasi-Fuchsian reciprocity, as well as the equation $\pa(P_F - P_{QF})=0$, immediately
follow from it. The classical action $S:\D(\Ga)\rightarrow\RR$ is symmetric with
respect to Riemann surfaces $X$ and $Y$,
\begin{equation} \label{Global}
S(X,Y)=S(\bar{Y},\bar{X}),
\end{equation}
where $\bar{X}$ is the mirror image of $X$, and this property manifests the global
quasi-Fuchsian reciprocity. Equation \eqref{McM} now follows from \eqref{F-QF} and
\eqref{I}. Its direct proof along the lines of~\cite{ZT87a,ZT87b} is given in Theorem
\ref{variation2}. As immediate corollary of \eqref{McM} and \eqref{F-QF}, we obtain
that function $-S$ is a K\"{a}hler potential of the Weil-Petersson metric on
$\D(\Ga)$.

Our second result is a precise relation between two and three-dimensional
constructions which proves the holography principle for the quasi-Fuchsian case. Let
$\up^3=\{Z=(x,y,t) \in\RR^3\,|\, t>0\}$ be the hyperbolic 3-space. The quasi-Fuchsian
group $\Ga$ acts discontinuously on $\up^3\cup\Omega$ and the quotient
$M\simeq\Ga\bk(\up^3\cup\Omega)$ is a hyperbolic 3-manifold with boundary
$\Ga\bk\Omega\simeq X\sqcup Y$. According to the holography principle (see, e.g.,
\cite{Manin} for mathematically oriented exposition), the regularized hyperbolic
volume of $M$ --- \emph{on-shell Einstein-Hilbert action with cosmological term}, is
related to the Liouville action functional $S[\phi]$.

In case when $\Ga$ is a classical Schottky group, i.e., when it has a fundamental
domain bounded by Euclidean circles, holography principle was established by
K.~Krasnov in~\cite{Krasnov}. Namely, let $M\simeq \Ga\bk(\up^3\cup\Omega)$ be the
corresponding hyperbolic 3-manifold (realized using the Ford fundamental region) with
boundary $X\simeq \Ga\bk\Omega$ --- a compact Riemann surface of genus $g>1$. For
every $ds^2 =e^{\phi}|dz|^2\in\mathcal{CM}(X)$ consider the family 
$\mathcal{H}_{\vep}$ of
surfaces given by the equation $f(Z)=t e^{\phi(z)/2} =\vep>0$ where $z=x+iy$, and let
$M_{\vep}=M\cap \mathcal{H}_\vep$. Denote by $V_{\vep}[\phi]$ the hyperbolic volume 
of $M_\vep$,
by $A_\vep[\phi]$ --- the area of the boundary of $M_{\vep}$ in the metric on 
$\mathcal{H}_\vep$
induced by the hyperbolic metric on $\up^3$, and by $A[\phi]$ --- the area of $X$ in
the metric $ds^2$. In~\cite{Krasnov} K.~Krasnov obtained the following formula
\begin{equation}\label{Krasnov}
 \lim_{\vep\rightarrow 0}\left(V_{\vep}[\phi] - \frac{1}{2} A_{\vep}[\phi] +
 (2g-2)\pi\log \vep\right)= - \frac{1}{4}\left(S[\phi] - A[\phi]\right).
\end{equation}
It relates three-dimensional data --- the regularized volume of $M$, to the
two-dimensional data --- the Liouville action functional $S[\phi]$, thus
establishing the holography principle. Note that the metric $ds^2$ on the boundary of
$M$ appears entirely through regularization by means of the surfaces
$\mathcal{H}_\vep$, which are not $\Ga$-invariant. As a result, arguments 
in~\cite{Krasnov} work only for classical Schottky groups.

We extend homological algebra methods in~\cite{AT} to the three-dimensional case when
$\Ga$ is a quasi-Fuchsian group. Namely, we construct $\Ga$-invariant cut-off function
$f$ using a partition of unity for $\Ga$, and prove in Theorem \ref{holography} that
on-shell regularized Einstein-Hilbert action functional
\begin{equation*}
\mathcal{E}[\phi] = -4\lim_{\vep \rightarrow 0} \left(V_{\ep}[\phi] -\frac{1}{2}
A_{\ep}[\phi] + 2\pi(2g-2)\log\vep \right),
\end{equation*}
is well-defined and satisfies the quasi-Fuchsian holography principle
\begin{equation*}
\mathcal{E}[\phi] = S[\phi] - \iint\limits_{\Ga\bk\Omega}\,e^{\phi}d^2 z - 8\pi
(2g-2)\log 2.
\end{equation*}
As immediate corollary we get another proof that the Liouville action functional
$S[\phi]$ does not depend on the choice of a fundamental domain $F$ of $\Ga$ in
$\Omega$, provided it is the boundary in $\Omega$ of a fundamental region of $\Gamma$
in $\up^3\cup\Omega$.

We also show that $\Ga$-invariant cut-off surfaces $\mathcal{H}_{\vep}$ can be 
chosen to be Epstein surfaces, which are naturally associated with the family of
metrics $ds^2_{\vep}=4\vep^{-2}e^\phi|dz|^2\in\mathcal{CM}(X)$ by the inverse
of the ``hyperbolic Gauss map'' \cite{Epstein1, Epstein2} (see also
~\cite{Anderson}). This construction also
gives a geometric interpretation of the density $|\phi_z|^2 + e^\phi$ in terms of 
Epstein surfaces. 

Schottky and quasi-Fuchsian groups considered above are basically the only examples of
geometrically finite, torsion-free, purely loxodromic Kleinian groups with finitely
many components. Indeed, according to the theorem of Maskit \cite{Maskit}, a 
geometrically finite, purely loxodromic Kleinian
group satisfying these properties has at most two components. The one-component case
corresponds to Schottky groups and the two-component case --- to Fuchsian or
quasi-Fuchsian groups and their $\ZZ_2$-extensions.

The third result of the paper is the generalization of main results for
quasi-Fuchsian groups --- Theorems \ref{variation1}, \ref{variation2} and
\ref{holography}, to Kleinian groups. Namely, we introduce a notion of a Kleinian
group of Class $A$ for which this generalization holds. By definition, a
non-elementary, purely loxodromic, geometrically finite Kleinian group is of Class $A$
if it has fundamental region $R$ in $\up^3\cup\Omega$ which is a finite
three-dimensional $CW$-complex with no vertices in $\up^3$. Schottky, Fuchsian,
quasi-Fuchsian groups, and their free combinations are of Class $A$, and Class $A$ is
stable under quasiconformal deformations. We extend three-dimensional homological
methods developed in Section 5 to the case of Kleinian group $\Ga$ of Class $A$ acting
on $\up^3\cup\Omega$. Namely, starting from the fundamental region $R$ for $\Ga$ in
$\up^3\cup\Omega$, we construct a chain of degree 3 in total homology complex
$\Tot\KKK$, whose boundary in $\Omega$ is the cycle $\Sigma$ of degree 2 for the
corresponding total homology complex of the region of discontinuity $\Omega$. In
Theorem \ref{K-holography} we establish holography principle for Kleinian groups: we
prove that the on-shell regularized Einstein-Hilbert action for the 3-manifold
$M\simeq\Ga\bk(\up^3\cup\Omega)$ is well-defined and is related to the Liouville
action functional for $\Ga$, defined by the evaluation map \eqref{Evaluation}. When
$\Ga$ is a Schottky group, we get the functional \eqref{Schottky} introduced
in~\cite{ZT87b}. As in the quasi-Fuchsian case, the Liouville action functional does
not depend on the choice of a fundamental domain $F$ for $\Ga$ in $\Omega$, as long as
it is the boundary in $\Omega$ of a fundamental region of $\Ga$ in $\up^3\cup\Omega$.
Denote by $\D(\Ga)$ the deformation space of the Kleinian group $\Ga$. To every point
$\Ga'\in\D(\Ga)$ with the region of discontinuity $\Omega'$ there corresponds a
disjoint union $X_1\sqcup\dots\sqcup X_n\simeq\Ga'\bk\Omega'$ of compact Riemann
surfaces simultaneously uniformized by the Kleinian group $\Ga'$. Conversely, by the
theorem of Maskit~\cite{Maskit}, for a given sequence of compact Riemann surfaces
$X_1,\dots,X_n$ there is a Kleinian group which simultaneously uniformizes them. Using
the same notation, we denote by $P_F$ projective connection on $X_1\sqcup\dots\sqcup
X_n$ given by the Fuchsian uniformization of these Riemann surfaces and by $P_K$ ---
projective connection given by their simultaneous uniformization  by a Kleinian group
($P_K=P_{QF}$ for the quasi-Fuchsian case). Let $S:\D(\Ga)\rightarrow \RR$ be the
classical Liouville action. Theorem \ref{K-global} states that
\begin{equation*}
P_{F}-P_K = \frac{1}{2}\,\pa S,
\end{equation*}
which is the ultimate generalization of \eqref{I}. Similarly, Theorem \ref{K-global-2}
is the statement
\begin{equation*}
\bar\pa(P_F - P_K) = - i\,\omega_{WP},
\end{equation*}
which implies that $-S$ is a \Ka potential of the Weil-Petersson metric on $\D(\Ga)$.
As another immediate corollary of Theorem \ref{K-global} we get McMullen's Kleinian
reciprocity
--- Theorem \ref{McM-2}.

Finally, we observe that our method and results, with appropriate modifications,
can be generalized to the case when quasi-Fuchsian and Class $A$ Kleinian groups 
have torsion
and contain parabolic elements. Our method also works for the Bers' universal \Te
space and the related infinite-dimensional K\"{a}hler manifold
$\Diff_{+}(S^1)/\Mob(S^1)$. We plan to discuss these generalizations elsewhere.

The content of the paper is the following. In Section 2 we give a construction of the
Liouville action functional following the method in \cite{AT}, which we review briefly
in 2.1. In Section 2.2 we define and establish the main properties of the Liouville
action functional in the model case when $\Ga$ is a Fuchsian group, and in Section 2.3
we consider technically more involved quasi-Fuchsian case. In Section 3 we recall all
necessary basic facts from the deformation theory. In Section 4 we prove our first
main result --- Theorems \ref{variation1} and \ref{variation2}. In Section 5 we prove
the second main result
--- Theorem \ref{holography} on quasi-Fuchsian
holography. Finally in Section 6 we generalize these results for Kleinian groups of
Class $A$: we define Liouville action functional and prove Theorems
\ref{K-holography}, \ref{K-global} and \ref{K-global-2}. \vspace{1mm}

\noindent \textbf{Acknowledgments.} We greatly appreciate stimulating discussions
with E.~Aldrovandi on homological methods, and with I.~Kra, 
M.~Lyubich and B. Maskit on
various aspects of the theory of Kleinian groups. We also are grateful to K. Krasnov
for valuable comments and to C. McMullen for the insightful suggestion to use 
Epstein surfaces for regularizing the Eistein-Hilbert action.
The work of the first author was partially supported by the NSF grant DMS-9802574.

\section{Liouville action functional}

Let $\Ga$ be a normalized, marked, purely loxodromic quasi-Fuchsian group of genus
$g>1$ with region of discontinuity $\Omega$, so that $\Ga\bk\Omega\simeq X\sqcup Y$,
where $X$ and $Y$ are compact Riemann surfaces of genus $g>1$ with opposite
orientations. Here we define Liouville action functional $S_\Ga$ for the group $\Ga$
as a functional on the space of smooth conformal metrics on $X\sqcup Y$ with the
property that its Euler-Lagrange equation is the Liouville equation on $X\sqcup Y$.
Its definition is based on the homological algebra methods developed in~\cite{AT}.

\subsection{Homology and cohomology set-up}
Let $\Ga$ be a group acting properly on a smooth manifold $M$. To this data one
canonically associates double homology and cohomoloy complexes (see, e.g., \cite{AT}
and references therein).

Let $\SSS_{\bu} \equiv \SSS_{\bu}(M)$ be the standard singular chain complex of $M$
with the differential $\pa'$. The group action on $M$ induces a left $\Ga$-action on
$\SSS_{\bu}$ by translating the chains and $\SSS_{\bu}$ becomes a complex of left
$\Ga$-modules. Since the action of $\Ga$ on $M$ is proper, $\SSS_{\bu}$ is a complex
of free left $\Z\Ga$-modules, where $\Z\Ga$ is the integral group ring of the group
$\Ga$. The complex $\SSS_{\bu}$ is endowed with a right $\Z\Ga$-module structure in
the standard fashion: $c\cdot\gamma= \gamma^{-1}(c)$.

Let $\BBB_{\bu} \equiv \BBB_{\bu}(\Z\Ga)$ be the canonical ``bar'' resolution complex
for $\Ga$ with differential $\pa''$. Each $\BBB_n (\Z\Ga)$ is a free left $\Ga$-module
on generators $[\s_1|\ldots|\s_n]$, with the differential $\pa'': \BBB_n
\longrightarrow \BBB_{n-1}$ given by
\begin{align*}
 \pa'' [\s_1|\ldots|\s_n] & = \s_1[\s_2|\ldots|\s_n] +
 \sum_{k=1}^{n-1}(-1)^k[\s_1|\ldots|\s_k\s_{k+1}|\ldots|\s_n]\\
    & \;\;\; + (-1)^n [\s_1|\ldots|\s_{n-1}]\,,\;\; n>1, \\
  \pa'' [\s]& = \s\,[\;] - [\;]\,,\;\; n=1,
\end{align*}
where $[\s_1|\ldots|\s_n]$ is zero if some $\s_i$ equals to the unit element $\id$ in
$\Ga$. Here $\BBB_0(\Z\Ga)$ is a $\Z\Ga$-module on one generator $[\;]$ and it can be
identified with $\Z\Ga$ under the isomorphism that sends $[\;]$ to $1$; by definition,
$\pa''[\;]=0$.

The double homology complex $\KKK_{\bu,\bu}$ is defined as $\SSS_{\bu}\otimes_{\Z\Ga}
\BBB_{\bu}$, where the tensor product over $\Z\Ga$ uses the right $\Ga$-module
structure on $\SSS_{\bu}$. The associated total complex $\Tot\KKK$ is equipped with
the total differential $\pa = \pa' + (-1)^p \pa''$ on $\KKK_{p,q}$, and the complex
$\SSS_{\bu}$ is identified with $\SSS_{\bu} \otimes_{\Z\Ga}\BBB_0$ by the isomorphism
$c\mapsto c\otimes [\;]$.

Corresponding double complex in cohomology is defined as follows. Denote by
$\AAA^{\bu} \equiv \AAA_{\C}^{\bu}(M)$ the complexified de Rham complex on $M$. Each
$\AAA^n$ is a left $\Ga$-module with the pull-back action of $\Ga$, i.e., $\s \cdot
\varpi = (\s^{-1})^* \varpi$ for $\varpi \in \AAA^{\bu}$ and $\s \in \Ga$. Define the
double complex $\CCC^{p,q}=\Hom_{\C}(\BBB_q, \AAA^p)$ with differentials $d$, the
usual de Rham differential, and $\de = (\pa'')^*$, the group coboundary. Specifically,
for $\varpi \in \CCC^{p,q}$,
\begin{align*}
(\de\varpi)_{\s_1,\cdots,\s_{q+1}} & = \s_1 \cdot \varpi_{\s_2,\cdots,\s_{q+1}} +
\sum_{k=1}^{q} (-1)^k \varpi_{\s_1,\cdots,\s_k\s_{k+1},\cdots,\s_{q+1}}\\
       &\;\;\; +(-1)^{q+1} \varpi_{\s_1,\cdots,\s_q}.
\end{align*}
We write the total differential on $\CCC^{p,q}$ as $D = d + (-1)^p \de$.

There is a natural pairing between $\CCC^{p,q}$ and $\KKK_{p,q}$ which assigns to the
pair $(\varpi, c\otimes [\s_1|\ldots|\s_q])$ the evaluation of the $p$-form
$\varpi_{\s_1,\cdots,\s_q}$ over the $p$-cycle $c$,
\begin{equation*}
\la \varpi, c\otimes [\s_1|\ldots|\s_q] \ra = \int_c \varpi_{\s_1,\cdots,\s_q}.
\end{equation*}
By definition,
\begin{equation*}
\la \de \varpi, c \ra = \la \varpi, \pa'' c \ra,
\end{equation*}
so that using Stokes' theorem we get
\begin{equation*}
\la D\varpi, c \ra = \la \varpi, \pa c \ra.
\end{equation*}
This pairing defines a non-degenerate pairing between corresponding cohomology and
homology groups $H^{\bu}(\Tot\CCC)$ and $H_{\bu}(\Tot\KKK)$, which we continue to
denote by $\la~,~\ra$. In particular, if $\Phi$ is a cocycle in $(\Tot\CCC)^n$ and $C$
is a cycle in $(\Tot\KKK)_n$, then the pairing $\la\Phi, C\ra$ depends only on
cohomology classes $[\Phi]$ and $[C]$ and not on their representatives.

It is this property that will allow us to define Liouville action functional by
constructing corresponding cocycle $\Psi$ and cycle $\Sigma$. Specifically, we
consider the following two cases.

\begin{itemize}
\item[\textbf{1.}] $\Ga$ is purely hyperbolic Fuchsian group of genus $g>1$ and
$M=\U$ --- the upper half-plane of the complex plane $\CC$. In this case, since $\U$
is acyclic, we have~\cite{AT}
\begin{equation*}
H_\bullet(X,\ZZ)\cong H_\bullet(\Gamma,\ZZ)\cong H_\bullet(\Tot\KKK)\,,
\end{equation*}
where the three homologies are: the singular homology of $X\simeq \Ga\bk\U$, a compact
Riemann surface of genus $g>1$, the group homology of $\Gamma$, and the homology of
the complex $\Tot\KKK$ with respect to the total differential $\del$. Similarly, for
$M=\lo$ --- the lower half-plane of the complex plane $\CC$, we have
\begin{equation*}
H_\bullet(\bar{X},\ZZ)\cong H_\bullet(\Gamma,\ZZ)\cong H_\bullet(\Tot\KKK)\,,
\end{equation*}
where $\bar{X}\simeq \Ga\bk\lo$ is the mirror image of $X$ --- a complex-conjugate of
the Riemann surface $X$.
\item[\textbf{2.}] $\Ga$ is purely loxodromic quasi-Fuchsian group of genus $g>1$
with region of discontinuity $\Omega$ consisting of two simply-connected components
$\Omega_1$ and $\Omega_2$ separated by a quasi-circle $\mathcal{C}$. The same
isomorphisms hold, where $X\simeq\Ga\bk\Omega_1$ and $\bar{X}$ is replaced by
$Y\simeq\Ga\bk\Omega_2$.
\end{itemize}
\subsection{The Fuchsian case}

Let $\Ga$ be a marked, normalized, purely hyperbolic Fuchsian group of genus $g>1$,
let $X\simeq\Ga\bk\U$ be corresponding marked compact Riemann surface of genus $g$,
and let $\bar{X}\simeq\Ga\bk\lo$ be its mirror image. In this case it is possible to
define Liouville action functionals on Riemann surfaces $X$ and $\bar{X}$ separately.
The definition will be based on the following specialization of the general
construction in Section 2.1.

\subsubsection{Homology computation}

Here is a representation of the fundamental class $[X]$ of the Riemann surface $X$ in
$H_2(X,\ZZ)$ as a cycle $\Sigma$ of total degree 2 in the homology complex $\Tot\KKK$
\cite{AT}.

Recall that the marking of $\Ga$ is given by a  system of $2g$ standard generators
$\alpha_1,\dots,\alpha_g, \beta_1,\dots , \beta_g$ satisfying the single relation
\begin{equation*}
\g_1\cdots\g_g=\id,
\end{equation*}
where $\g_k=[\al_k,\be_k]=\al_k\be_k\al_k^{-1}\be_k^{-1}$. The marked group $\Ga$ is
normalized, if the attracting and repelling fixed points of $\al_1$ are, respectively,
$0$ and $\infty$, and the attracting fixed point of $\be_1$ is $1$. Every marked
Fuchsian group $\Ga$ is conjugated in $\PSL(2,\RR)$ to a normalized marked Fuchsian
group. For a given marking there is a standard choice of the fundamental domain
$F\subset \U$ for $\Ga$ as a closed non-Euclidean polygon with $4g$ edges labeled by
$a_k,a_k',b_k',b_k$ satisfying $\al_k(a_k')=a_k,\,\be_k(b_k')=b_k, \,k=1,2,\dots,g$
(see Fig.~1). The orientation of the edges is chosen such that
\begin{equation*}
\pa' F=\sum_{k=1}^{g}(a_k+b_k'-a_k'-b_k).
\end{equation*}
Set $\pa' a_k=a_k(1)-a_k(0),\, \pa' b_k=b_k(1)-b_k(0)$, so that $a_k(0)=b_{k-1}(0)$.
The relations between the vertices of $F$ and the generators of $\Ga$ are the
following: $\al_k^{-1}(a_k(0))=b_k(1),\,\be_k^{-1}(b_k(0))=
a_k(1),\,\g_k(b_k(0))=b_{k-1}(0)$, where $b_0(0)=b_g(0)$.

According to the isomorphism $\SSS_{\bu}\simeq\KKK_{\bu,0}$, the fundamental domain
$F$ is identified with $F \otimes [\;] \in \KKK_{2,0}$. We have $\pa'' F = 0$ and, as
it follows from the previous formula,
\begin{equation*}
\pa' F = \sum_{k=1}^g \left(\be_k^{-1}(b_k) - b_k - \al_k^{-1}(a_k) + a_k\right) =
\pa'' L,
\end{equation*}
where $L \in \KKK_{1,1}$ is given by
\begin{equation} \label{L}
L = \sum_{k=1}^g \left(b_k\otimes[\be_k] - a_k\otimes[\al_k]\right).
\end{equation}
There exists $V\in\KKK_{0,2}$ such that $\pa'' V = \pa ' L$. A straightforward
computation gives the following explicit expression
\begin{align} \label{V}
V &= \sum_{k=1}^g \left(a_k(0)\otimes[\al_k|\be_k] - b_k(0)\otimes[\be_k|\al_k] +
b_k(0)\otimes\left[\g_k^{-1}|\al_k\be_k\right]\right)\\ & \;\;\; -\sum_{k=1}^{g-1}
b_g(0)\otimes\left[\g_g^{-1}\ldots\g_{k+1}^{-1}|\g_k^{-1}\right]. \nonumber
\end{align}

\begin{figure}
\centering \epsfxsize=.30\linewidth \epsffile{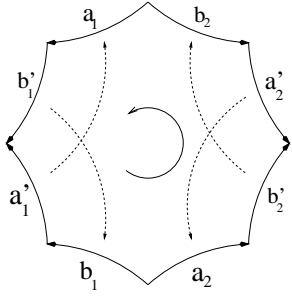}
\caption{Conventions for the fundamental domain $F$} \label{fig:conventions}
\end{figure}


Using $\pa''F=0$, $\pa' F = \pa'' L$,  $\pa'' V = \pa ' L$, and $\pa' V = 0$, we
obtain that the element $\Sigma = F + L - V$ of total degree 2 is a cycle in
$\Tot\KKK$, that is $\pa \Sigma = 0$. The cycle $\Sigma\in (\Tot\KKK)_2$ represents
the fundamental class $[X]$. It is proved in \cite{AT} that corresponding homology
class $[\Sigma]$ in $H_{\bu}(\Tot\KKK)$ does not depend on the choice of the
fundamental domain $F$ for the group $\Ga$.

\subsubsection{Cohomology computation}
Corresponding construction in cohomology is the following. Start with the space
$\mathcal{CM}(X)$ of all conformal metrics on $X\simeq \Ga\bk\U$. Every $ds^2\in
\mathcal{CM}(X)$ can be represented as $ds^2=e^\phi|dz|^2$, where $\phi\in
C^\infty(\U,\RR)$ satisfies
\begin{equation} \label{field}
\phi \circ \s + \log |\s'|^2 = \phi \;\;\; \text{for all} \; \s \in \Ga.
\end{equation}
In what follows we will always identify $\mathcal{CM}(X)$ with the affine subspace of
$C^\infty(\U,\RR)$ defined by \eqref{field}.

The ``bulk'' 2-form $\omega$ for the Liouville action is given by
\begin{equation} \label{omega}
\omega[\phi] = \left(|\phi_{z}|^2 +e^{\phi}\right)dz\wedge d\z,
\end{equation}
where $\phi\in\mathcal{CM}(X)$. Considering it as an element in $\CCC^{2,0}$ and
using~\eqref{field} we get
\begin{equation*}
\delta\omega[\phi]=d\theta[\phi],
\end{equation*}
where $\theta[\phi] \in \CCC^{1,1}$ is given explicitly by
\begin{equation} \label{theta}
\theta_{\s^{-1}}[\phi] = \left(\phi -  \frac{1}{2}\log|\s'|^2\right)
\left(\frac{\s''}{\s'} dz - \frac{\ov{\s''}}{\ov{\s'}}d\z\right).
\end{equation}

Next, set
\begin{equation*}
u=\de \theta[\phi] \in \CCC^{1,2}.
\end{equation*}
From the definition of $\theta$ and $\de^2=0$ it follows that the 1-form $u$ is
closed. An explicit calculation gives
\begin{align} \label{u}
u_{\s_1^{-1},\s_2^{-1}}= & -\frac{1}{2}\log|\s_1'|^2
\left(\frac{\s_2''}{\s_2'}\circ\s_1\, \s_1'\, dz -
\frac{\ov{\s_2''}}{\ov{\s_2'}}\circ\s_1\, \ov{\s_1'}\,d\z\right) \\
 & + \frac{1}{2}\log|\s_2'\circ\s_1|^2\left(\frac{\s_1''}{\s_1'} dz -
\frac{\ov{\s_1''}}{\ov{\s_1'}}d\z\right), \nonumber
\end{align}
and shows that $u$ does not depend on $\phi\in\mathcal{CM}(X)$.
\begin{remark}
The explicit formulas above are valid in the general case, when domain
$\Omega\subset\hat{\CC}$ is invariant under the action of a Kleinian group $\Ga$.
Namely, define the 2-form $\omega$ by formula \eqref{omega}, where $\phi$ satisfies
\eqref{field} in $\Omega$. Then solution $\theta$ to the equation
$\delta\omega[\phi]=d\theta[\phi]$ is given by the formula \eqref{theta} and $u=\de
\theta[\phi]$ --- by \eqref{u}.
\end{remark}

There exists a cochain $\Theta\in\CCC^{0,2}$ satisfying
\begin{equation*}
d\Theta =u~\text{and}~ \de\Theta =0.
\end{equation*}
Indeed, since the 1-form $u$ is closed and $\U$ is simply-connected, $\Theta$ can be
defined as a particular antiderivative of $u$ satisfying $\de\Theta=0$. This can be
done as follows. Consider the hyperbolic (Poincar\'{e}) metric on $\U$
\begin{equation*}
e^{\phi_{hyp}(z)}|dz|^2=\frac{|dz|^2}{y^2},\;\,z=x+iy\in\U.
\end{equation*}
This metric is $\PSL(2,\RR)$-invariant and its push-forward to $X$ is a hyperbolic
metric on $X$. Explicit computation yields
\begin{equation*}
\omega[\phi_{hyp}]=2 e^{\phi_{hyp}}\,dz\wedge d\z,
\end{equation*}
so that $\delta\omega[\phi_{hyp}]=0$. Thus the 1-form $\theta[\phi_{hyp}]$ on $\U$ is
closed and, therefore, is exact,
\begin{equation*}
\theta[\phi_{hyp}]=dl,
\end{equation*}
for some $l\in \CCC^{0,1}$. Set
\begin{equation} \label{Theta}
\Theta=\de l.
\end{equation}
It is now immediate that $\de \Theta = 0$ and $\de \theta[\phi]=u=d\Theta$ for all
$\phi\in\mathcal{CM}(X)$. Thus $\Psi[\phi]=\omega[\phi] -\theta[\phi] -\Theta$ is a
2-cocycle in the cohomology complex $\Tot\CCC$, that is, $D\Psi[\phi] =0$.
\begin{remark} \label{general} For every $\s\in\PSL(2,\RR)$ define the 1-form
$\theta_{\s}\left[\phi_{hyp}\right]$ by the same formula \eqref{theta},
\begin{equation} \label{theta-hyp}
\theta_{\s^{-1}}[\phi_{hyp}] = -\left(2\log y + \frac{1}{2}\log|\s'|^2\right)\left(
\frac{\s''}{\s'} dz - \frac{\ov{\s''}}{\ov{\s'}} d\z\right).
\end{equation}
Since for every $\s\in\PSL(2,\RR)$
\begin{equation*}
(\de\log y)_{\s^{-1}}=\log(y\circ\s) - \log y =\frac{1}{2}\log|\s'|^2,
\end{equation*}
the 1-form $u=\de\theta[\phi]$ is still given by \eqref{u} and is a
$\AAA^1(\up)$-valued group 2-cocycle for $\PSL(2,\RR)$, that is, $(\de
u)_{\s_1,\s_2,\s_3}=0$ for all $\s_1,\s_2,\s_3\in\PSL(2,\RR)$. Also 0-form $\Theta$
given by $\eqref{Theta}$ satisfies $d\Theta = u$ and is a $\AAA^0(\up)$-valued group
2-cocycle for $\PSL(2,\RR)$.
\end{remark}
\subsubsection{The action functional}
The evaluation map $\la\Psi[\phi],\Sigma\ra$ does not depend on the choice of the
fundamental domain $F$ for $\Ga$ \cite{AT}. It also does not depend on a particular
choice of antiderivative $l$, since by the Stokes' theorem
\begin{equation} \label{L2}
\la\Theta, V\ra = \la\de l,V\ra = \la l,\pa'' V\ra = \la l,\pa' L\ra =
\la\theta[\phi_{hyp}],L\ra.
\end{equation}
This justifies the following definition.
\begin{definition} The Liouville action functional $S[\,\cdot \, ;X]: \mathcal{CM}(X)
\rightarrow \RR$ is defined by the evaluation map
\begin{equation*}
S[\phi;X]=\frac{i}{2}\left\la\Psi[\phi],\Sigma\right\ra,~\phi\in\mathcal{CM}(X).
\end{equation*}
\end{definition}

For brevity, set $S[\phi]=S[\phi;X]$. The following lemma shows that the difference of
any two values of the functional $S$ is given by the bulk term only.
\begin{lemma} \label{1st-variation}
For all $\phi\in\mathcal{CM}(X)$ and $\sigma\in C^\infty(X,\RR)$,
\begin{equation*}
S[\phi+\sigma]-S[\phi]=\iint\limits_F\left(|\sigma_z|^2 + \left(e^\sigma + K\,\sigma -
1\right)e^\phi\right) d^2z,
\end{equation*}
where $d^2z=dx\wedge dy$ is the Lebesgue measure and $K=-2e^{-\phi}\phi_{z\z}$ is the
Gaussian curvature of the metric $e^{\phi}|dz|^2$.
\end{lemma}
\begin{proof}
We have
\begin{equation*}
\omega[\phi+\sigma]-\omega[\phi]= \omega[\phi;\sigma] + d\tilde\theta,
\end{equation*}
where
\begin{equation*}
\omega[\phi;\sigma]=\left(|\sigma_z|^2 + \left(e^\sigma + K\,\sigma - 1\right)
e^\phi\right) dz\wedge d\z,
\end{equation*}
and
\begin{equation*}
\tilde\theta = \sigma\left(\phi_{\z} d\z - \phi_z dz\right).
\end{equation*}
Since
\begin{equation*}
\de\tilde\theta_{\s^{-1}}= \sigma\left(\frac{\s''}{\s'} dz -
\frac{\ov{\s''}}{\ov{\s'}} d\z\right)=\theta[\phi + \sigma] - \theta[\phi],
\end{equation*}
the assertion of the lemma follows from the Stokes' theorem.
\end{proof}

\begin{corollary} \label{crit-value}
The Euler-Lagrange equation for the functional $S$ is the Liouville equation, the
critical point of $S$ --- the hyperbolic metric $\phi_{hyp}$, is non-degenerate, and
the classical action --- the critical value of $S$, is twice the hyperbolic area of
$X$, that is, $4\pi(2g-2)$.
\end{corollary}
\begin{proof}
As it follows from Lemma \ref{1st-variation},
\begin{equation*}
\left.\frac{dS[\phi + t\sigma] }{dt}\right|_{t=0}=\iint\limits_F\left(K+1\right)\sigma
e^\phi d^2z,
\end{equation*}
so that the Euler-Lagrange equation is the Liouville equation $K=-1$. Since
\begin{equation*}
\left.\frac{d^2 S[\phi_{hyp} + t\sigma]}{dt^2}\right|_{t=0} =\iint\limits_F
\left(2|\sigma_z|^2 + \sigma^2 e^{\phi_{hyp}}\right) d^2z >0\;\;\,\text{if}\;\;\,
\sigma\neq 0,
\end{equation*}
the critical point $\phi_{hyp}$ is non-degenerate. Using \eqref{L2} we get
\begin{equation*}
S[\phi_{hyp}]=\frac{i}{2}\la\Psi[\phi_{hyp}],\Sigma\ra=
\frac{i}{2}\la\omega[\phi_{hyp}],F\ra=2\iint\limits_F\frac{d^2z}{y^2}=4\pi(2g-2).
\end{equation*}
\end{proof}

\begin{remark} \label{Polyakov}
Let $\Delta[\phi]=-e^{-\phi}\pa_z \pa_{\z}$ be the Laplace operator of the metric
$ds^2=e^\phi |dz|^2$ acting on functions on $X$, and let $\det\Delta[\phi]$ be its
zeta-function regularized determinant (see, e.g., \cite{OPS} for details). Denote by
$A[\phi]$ the area of $X$ with respect to the metric $ds^2$ and set
\begin{equation*}
\mathcal{I}[\phi]=\log\frac{\det\Delta[\phi]}{A[\phi]}.
\end{equation*}
The Polyakov's ``conformal anomaly'' formula~\cite{Pol} reads
\begin{equation*}
\mathcal{I}[\phi + \sigma] - \mathcal{I}[\phi] = -\frac{1}{12\pi} \iint\limits_F
\left(|\sigma_z|^2 + K\sigma\,e^\phi\right)d^2z,
\end{equation*}
where $\sigma\in C^\infty(X,\RR)$ (see~\cite{OPS} for rigorous proof). Comparing it
with Lemma \ref{1st-variation} we get
\begin{equation*}
\mathcal{I}[\phi + \sigma] + \frac{1}{12\pi} \check{S}[\phi + \sigma]=
\mathcal{I}[\phi] +  \frac{1}{12\pi} \check{S}[\phi],
\end{equation*}
where $\check{S}[\phi]=S[\phi] - A[\phi]$.
\end{remark}

Lemma \ref{1st-variation}, Corollary \ref{crit-value} (without the assertion on
classical action) and Remark \ref{Polyakov} remain valid if $\Theta$ is replaced by
$\Theta + c$, where $c$ is an arbitrary group 2-cocycle with values in $\C$. The
choice \eqref{Theta}, or rather its analog for the quasi-Fuchsian case, will be
important in Section 4, where we consider classical action for families of Riemann
surfaces. For this purpose, we present an explicit formula for $\Theta$ as a
particular antiderivative of the 1-form $u$.

Let $p\in \ov{\U}$ be an arbitrary point on the closure of $\U$ in $\C$ (nothing will
depend on the choice of $p$). Set
\begin{equation} \label{l}
l_\s(z)=\int_p^z\theta_\s [\phi_{hyp}]~\text{for all}~\s\in\Ga,
\end{equation}
where the path of integration $P$ connects points $p$ and $z$ and, possibly except $p$,
lies entirely in $\U$. If $p\in\RR_{\infty}=\RR\cup\{\infty\}$, it is assumed that $P$
is smooth and is not tangent to $\RR_{\infty}$ at $p$. Such paths are called
admissible. A 1-form $\vartheta$ on $\U$ is called integrable along admissible path
$P$ with the endpoint $p\in\RR_\infty$, if the limit of $\int_{p'}^z\vartheta$, as
$p'\rightarrow p$ along $P$, exists. Similarly, a path $P$ is called $\Ga$-closed if
its endpoints are $p$ and $\s p$ for some $\s\in\Ga$, and $P\setminus\{p,\s p\}\subset
\U$. A $\Ga$-closed path $P$ with endpoints $p$ and $\s p$, $p\in\RR_\infty$, is
called admissible if it is not tangent to $\RR_{\infty}$ at $p$ and there exists
$p'\in P$ such that the translate by $\s$ of the part of $P$ between the points $p'$
and $p$ belongs to $P$. A 1-form $\vartheta$ is integrable along $\Ga$-closed
admissible path $P$, if the limit of $\int_{p'}^{\s p'}\vartheta$, as $p'\rightarrow
p$ along $P$, exists.

Let
\begin{align} \label{W}
W=& \sum_{k=1}^{g} \left( P_{k-1} \otimes [\al_k|\be_k] - P_k \otimes[\be_k|\al_k] +
P_k \otimes \left[\g_k^{-1}|\al_k\be_k\right]\right)\\ & -\sum_{k=1}^{g-1} P_g
\otimes\left[\g_g^{-1}\ldots\g_{k+1}^{-1}|\g_k^{-1}\right]\in\KKK_{1,2}, \nonumber
\end{align}
where $P_k$ is any admissible path from $p$ to $b_k(0)$, $k=1,\dots,g$, and $P_g=P_0$.
Since $P_k(1) = b_k(0) = a_{k+1}(0)$, we have
\begin{equation*}
\pa' W = V-U,
\end{equation*}
where
\begin{align} \label{U}
U & = \sum_{k=1}^g \left(p \otimes [\al_k|\be_k] - p \otimes[\be_k|\al_k] + p \otimes
\left[\g_k^{-1}|\al_k\be_k\right]\right)\\ &\;\;\;\; - \sum_{k=1}^{g-1} p
\otimes\left[\g_g^{-1}\ldots\g_{k+1}^{-1}|\g_k^{-1}\right]\in\KKK_{1,2}.\nonumber
\end{align}
We have the following statement.
\begin{lemma}\label{el}
Let $\vartheta \in \CCC^{1,1}$ be a closed 1-form on $\U$ and $p\in\ov{\U}$. In case
$p\in\RR_{\infty}$ suppose that $\de\vartheta$ is integrable along any admissible path
with endpoints in $\Ga\cdot p$ and $\vartheta$ is integrable along any $\Ga$-closed
admissible path with endpoints in $\Ga\cdot p$. Then
\begin{align*}
\la \vartheta, L \ra = &\;\;\, \la \de\vartheta, W \ra \\ & + \sum_{k=1}^g\left(
\int_p^{\al_k^{-1}p}\vartheta_{\be_k} - \int_p^{\be_k^{-1}p} \vartheta_{\al_k} +
\int_p^{\g_k p} \vartheta_{\al_k\be_k} - \int_p^{\g_{k+1}\ldots \g_g p}
\vartheta_{\g_k^{-1}}\right),
\end{align*}
where paths of integration are admissible if $p\in\RR_{\infty}$.
\end{lemma}
\begin{proof}
Since $\vartheta_{\s}$ is closed and $\U$ is simply-connected, we can define
function $l_{\s}$ on $\U$ by
\begin{displaymath}
l_{\s}(z) = \int_p^z \vartheta_{\s},
\end{displaymath}
where $p\in\U$. We have, using Stokes' theorem and $d(\de l)= \de (d l) =
\de\vartheta$,
\begin{align*}
\la \vartheta, L \ra &= \la dl, L \ra = \la l, \pa' L \ra
                   = \la l, \pa'' V \ra
                   = \la \de l, V \ra \\
                  &= \la \de l, \pa' W\ra + \la \de l, U \ra
                  = \la d(\de l), W \ra + \la \de l, U \ra\\
                  &= \la \de\vartheta, W \ra \;\;+ \;\; \la \de l, U\ra.
\end{align*}
Since
\begin{equation*}
(\de l)_{\s_1, \s_2}(p) = \int_p^{\s_1^{-1}p} \vartheta_{\s_2},
\end{equation*}
we get the statement of the lemma if $p\in\U$. In case $p\in\RR_\infty$, replace $p$
by $p'\in\up$. Conditions of the lemma guarantee the convergence of integrals as
$p'\rightarrow p$ along corresponding paths.
\end{proof}
\begin{remark} Expression $\la\de l,U\ra$, which appears in the statement of the lemma,
does not depend on the choice of a particular antiderivative of the closed 1-form
$\vartheta$. The same statement holds if we only assume that 1-form $\de\vartheta$ is
integrable along admissible paths with endpoints in $\Ga\cdot p$, and 
1-form $\vartheta$ has an
antiderivative $l$ (not necessarily vanishing at $p$) such that the limit of $(\de
l)_{\s_1, \s_2}(p')$, as $p'\rightarrow p$ along admissible paths, exists.
\end{remark}
\begin{lemma} \label{constants} We have
\begin{equation} \label{Theta-up}
\Theta_{\s_1,\s_2}(z)=\int_{p}^{z}u_{\s_1,\s_2} + \eta(p)_{\s_1,\s_2},
\end{equation}
where $p\in\RR\setminus\Ga(\infty)$ and integration goes along admissible paths. The
integration constants $\eta\in\CCC^{0,2}$ are given by
\begin{equation} \label{eta}
\eta(p)_{\s_1,\s_2}=4\pi i\varepsilon(p)_{\s_1,\s_2} (2\log 2+ \log|c(\s_2)|^2),
\end{equation}
and
\begin{equation*}
     \varepsilon(p)_{\s_1,\s_2} = \begin{cases}
                              \;\;\,1 & \text{if $p < \s_2(\infty) < \s_1^{-1} p$}, \\
                             -1 & \text{if $p > \s_2 (\infty) > \s_1^{-1} p$},\\
                              \;\;\,0 & \text{otherwise}.
                             \end{cases}
\end{equation*}
Here for $\gamma=\left(\begin{smallmatrix} a & b \\ c & d
\end{smallmatrix}\right)$ we set $c(\s)=c$.
\end{lemma}
\begin{proof} Since
\begin{equation*}
\Theta_{\s_1,\s_2}(z)=\int_{p}^{z}u_{\s_1,\s_2} +  \int_p^{\s_1^{-1}p}
\theta_{\s_2}[\phi_{hyp}],
\end{equation*}
it is sufficient to verify that
\begin{equation*}
     \frac{1}{2\pi i}\int_p^{\s_1 p} \theta_{\s_2^{-1}}[\phi_{hyp}] =
                            \begin{cases}
                              \;\;\,4\log 2 + 2\log|c(\s_2)|^2 & \text{if
                              $p < \s_2^{-1}(\infty) < \s_1\,p$},\\
                             -4\log 2 - 2\log|c(\s_2)|^2
                             & \text{if $p > \s_2^{-1}(\infty) > \s_1\,p$},\\
                              \;\;\,0 & \text{otherwise.}
                             \end{cases}
\end{equation*}
From \eqref{theta-hyp} it follows that $\theta_{\s^{-1}}[\phi_{hyp}]$ is a closed
1-form on $\U$, integrable along admissible paths with $p\in\RR\setminus
\{\s^{-1}(\infty)\}$. Denote by $\theta_{\s^{-1}}^{(\varepsilon)}$ its restriction on
the line $y=\varepsilon>0$, $z=x+iy$. When $x\neq \s_2^{-1}(\infty)$, we obviously
have
\begin{equation*}
\lim_{\varepsilon\rightarrow 0}\theta_{\s_2^{-1}}^{(\varepsilon)}=0,
\end{equation*}
uniformly in $x$ on compact subsets of $\RR\setminus\{\s_2^{-1}(\infty)\}$.

If $\s_2^{-1} (\infty)$ does not lie between points $p$ and $\s_1 p$ on $\RR$, we can
approximate the path of integration by the interval on the line $y=\vep$, which tends
to $0$ as $\varepsilon\rightarrow 0$. If $\s_2^{-1} (\infty)$ lies between points $p$
and $\s_1 p$, we have to go around the point $\s_2^{-1}(\infty)$ via a small
half-circle, so that
\begin{equation*}
\int_p^{\s_1 p} \theta_{\s_2^{-1}}[\phi_{hyp}] = \lim_{r\rightarrow 0}
\int_{C_r}\theta_{\s_2^{-1}}[\phi_{hyp}],
\end{equation*}
where $C_r$ is the upper-half of the circle of radius $r$ with center at
$\s_2^{-1}(\infty)$, oriented clockwise if $p<\s_2^{-1}(\infty)<\s_1 p$. 
Evaluating the limit using elementary formula
\begin{equation*}
\int_0^\pi\log\sin t\,d t= -\pi \log 2,
\end{equation*}
and Cauchy theorem, we get the formula.
\end{proof}
\begin{corollary}
The Liouville action functional has the following explicit representation
\begin{equation*}
S[\phi] = \frac{i}{2}\left(\la \omega[\phi], F \ra - \la \theta[\phi], L \ra + \la u,
W \ra + \la \eta, V\ra\right).
\end{equation*}
\end{corollary}
\begin{remark} \label{independence}
Since $\la \Theta, V \ra = \la u, W \ra + \la \eta, V\ra$, it immediately follows from
\eqref{L2} that the Liouville action functional does not depend on the choice of
point $p\in \RR\setminus \Ga(\infty)$ (actually it is sufficient to assume that $p\neq
\s_1(\infty), (\s_1\s_2)(\infty)$ for all $\s_1, \s_2\in \Ga$ such that
$V_{\s_1,\s_2}\neq 0$). This can also be proved by direct computation using Remark
\ref{general}. Namely, let $p'\in\RR_\infty$ be another choice, $p'=\sigma^{-1}
p\in\RR_{\infty}$ for some $\sigma\in\PSL(2,\RR)$. Setting $z=p$ in the equation
$(\de\Theta)_{\sigma,\s_1,\s_2}=0$ and using $(\de u)_{\sigma,\s_1,\s_2}=0$, where
$\s_1, \s_2\in \Ga$, we get
\begin{equation} \label{eta-integral-1}
\int_p^{\sigma^{-1}p}u_{\s_1,\s_2} = -(\de\eta(p))_{\sigma,\s_1,\s_2},
\end{equation}
where all paths of integration are admissible. Using
\begin{equation*}
\eta(p)_{\sigma\s_1,\s_2} = \eta(\sigma^{-1}p)_{\s_1,\s_2} + \eta(p)_{\sigma,\s_2},
\end{equation*}
we get from \eqref{eta-integral-1} that
\begin{equation*}
\int_p^z u_{\s_1,\s_2} + \eta(p)_{\s_1,\s_2} = \int_{p'}^z u_{\s_1,\s_2} +
\eta(p')_{\s_1,\s_2} + (\de \eta_\sigma )_{\s_1,\s_2},
\end{equation*}
where $(\eta_\sigma)_{\s} = \eta(p)_{\sigma,\s}$ is constant group 1-cochain. The
statement now follows from
\begin{equation*}
\la\de\eta_\sigma, V \ra =  \la\eta_\sigma,\pa'' V\ra = \la\eta_\sigma,\pa' L\ra = \la
d\eta_\sigma, L\ra =0.
\end{equation*}
\end{remark}

Another consequence of Lemmas \ref{el} and \ref{constants} is the following.
\begin{corollary} \label{varkappa-lemma}
Set
\begin{equation*}
\varkappa_{\s^{-1}}= \frac{\s''}{\s'}dz - \frac{\ov{\s''}}{\ov{\s'}}d\z \in
\CCC^{1,1}.
\end{equation*}
Then
\begin{equation*}
\la \varkappa,L\ra = 4\pi i\la\vep,V\ra= 4\pi i\,\chi(X),
\end{equation*}
where $\chi(X) = 2-2g$ is the Euler characteristic of Riemann surface $X\simeq
\Ga\bk\up$.
\end{corollary}
\begin{proof}
Since $\de\varkappa=0$, the first equation immediately follows from the proofs of
Lemmas \ref{el} and \ref{constants}. To prove the second equation, observe that
\begin{displaymath}
\varkappa = \de \varkappa_1, \quad \text{where}\quad \varkappa_1 = - \phi_z dz +
\phi_{\z} d\z \quad \text{and} \quad d \varkappa_1 = 2 \phi_{z \z}\,dz\wedge d\z.
\end{displaymath}
Therefore
\begin{equation*}
\la \varkappa, L \ra = \la \de \varkappa_1, L \ra = \la \varkappa_1, \pa'' L \ra = \la
\varkappa_1, \pa' F \ra = \la d \varkappa_1, F \ra.
\end{equation*}
The Gaussian curvature of the metric $ds^2 = e^{\phi} |dz|^2$ is $K = -2e^{-\phi}
\phi_{z\z}$, so by Gauss-Bonnet we get
\begin{equation*}
\la d \varkappa_1, F \ra = 2 \iint\limits_F\phi_{z \z}\,dz\wedge d\z = 2i
\iint\limits_{\Ga\bk\up} K e^\phi d^2z = 4\pi i \chi(X).
\end{equation*}
\end{proof}
Using this corollary, we can ``absorb'' the integration constants $\eta$ by shifting
$\theta[\phi]\in\CCC^{1,1}$ by a multiple of closed 1-form $\varkappa$. Indeed,
1-form $\theta[\phi]$ satisfies the equation $\de\omega[\phi] =d\theta[\phi]$ and is
defined up to addition of a closed 1-form. Set
\begin{equation} \label{theta-prime}
\ptheta_{\s}[\phi] = \theta_{\s}[\phi] - (2\log 2 +\log|c(\s)|^2) \varkappa_{\s},
\end{equation}
and define $\pu=\de\ptheta[\phi]$. Explicitly,
\begin{align} \label{u-prime}
\pu_{\s_1^{-1}, \s_2^{-1}} = u_{\s_1^{-1}, \s_2^{-1}} & - \log
\frac{|c(\s_2)|^2}{|c(\s_2\s_1)|^2} \left( \frac{\s_2''}{\s_2'}\circ \s_1\,\s'_1 dz -
\frac{\ov{\s_2''}}{\ov{\s_2'}}\circ \s_1\,\ov{\s'_1} d\z\right)  \\ & + \log
\frac{|c(\s_2\s_1)|^2}{|c(\s_1)|^2} \left(\frac{\s_1''}{\s_1'}\,dz -
\frac{\ov{\s_1''}}{\ov{\s_1'}}\,d\z\right), \nonumber
\end{align}
where $u$ is given by \eqref{u}. As it follows from Lemma \ref{el} and Corollary
\ref{varkappa-lemma},
\begin{equation} \label{shift-1}
S[\phi] = \frac{i}{2}\left(\la \omega[\phi], F \ra - \la \ptheta[\phi], L \ra + \la
\pu, W \ra \right).
\end{equation}

Liouville action functional for the mirror image $\bar{X}$ is defined similarly.
Namely, for every chain $c$ in the upper half-plane $\U$ denote by $\bar{c}$ its
mirror image in the lower half-plane $\lo$; chain $\bar{c}$ has an opposite
orientation to $c$. Set $\bar{\Sigma}=\bar{F}+\bar{L}-\bar{V}$, so that $\pa
\bar{\Sigma} = 0$. For $\phi\in\mathcal{CM}(\bar{X})$, considered as a smooth
real-valued function on $\lo$ satisfying \eqref{field}, define
$\omega[\phi]\in\CCC^{2,0},\, \theta[\phi]\in\CCC^{1,1}$ and $\Theta\in\CCC^{0,2}$ by
the same formulas \eqref{omega}, \eqref{theta} and \eqref{Theta}. Lemma
\ref{constants} has an obvious analog for the lower half-plane $\lo$, the analog of
formula \eqref{Theta-up} for $z\in\lo$ is
\begin{equation} \label{Theta-lo}
\Theta_{\s_1,\s_2}(z)=\int_{p}^{z}u_{\s_1,\s_2} - \eta(p)_{\s_1,\s_2},
\end{equation}
where the negative sign comes from the opposite orientation.
\begin{remark} \label{identity}
Similarly to \eqref{eta-integral-1} we get
\begin{equation} \label{eta-integral-2}
\int_p^{\sigma^{-1} p}u_{\s_1,\s_2} = (\de\eta(p))_{\sigma,\s_1,\s_2},
\end{equation}
where the path of integration, except the endpoints, lies in $\lo$. From
\eqref{eta-integral-1} and \eqref{eta-integral-2} we obtain
\begin{equation} \label{eta-integral}
\int_Cu_{\s_1,\s_2} = -2(\de\eta(p))_{\sigma,\s_1,\s_2},
\end{equation}
where the path of integration $C$ is a loop that starts at $p$, goes to $\sigma^{-1}
p$ inside $\U$, continues inside $\lo$ and ends at $p$. Note that formula
\eqref{eta-integral} can also be verified directly using Stokes' theorem. Indeed, the
$1$-form $u_{\s_1,\s_2}$ is closed and regular everywhere except points
$\s_1(\infty)~\text{and}~(\s_1\s_2)(\infty)$. Integrating over small circles
around these points if they lie inside $C$ and using \eqref{eta}, we get the result.
\end{remark}

Set $\Psi[\phi]=\omega[\phi]-\theta[\phi] - \Theta$, so that $D\Psi[\phi]=0$. The
Liouville action functional for $\bar{X}$ is defined by
\begin{equation*}
[\phi;\bar{X}]=-\frac{i}{2}\la \Psi[\phi],\bar{\Sigma}\ra.
\end{equation*}
Using an analog of Lemma \ref{el} in the lower half-plane $\lo$ and
\begin{equation*}
\la \eta,\bar{V}\ra=\la \eta,V\ra,
\end{equation*}
we obtain
\begin{equation*}
S[\phi;\bar{X}]= -\frac{i}{2}\left(\la \omega[\theta],\bar{F}\ra - \la
\theta[\phi],\bar{L}\ra + \la u,\bar{W}\ra - \la \eta,V\ra\right).
\end{equation*}

Finally, we have the following definition.

\begin{definition} The Liouville action functional $S_{\Ga}:
\mathcal{CM} (X\sqcup\bar{X})\rightarrow \RR$ for the Fuchsian group $\Ga$ acting on
$\up\cup\lo$ is defined by
\begin{align*}
S_{\Ga}[\phi]= & S[\phi;X] + S[\phi;\bar{X}]= \frac{i}{2}\la \Psi[\phi],\Sigma -
\bar{\Sigma}\ra \\ = &\frac{i}{2}\left(\la \omega[\phi],F-\bar{F}\ra - \la
\theta[\phi],L-\bar{L}\ra +\la u, W-\bar{W}\ra + 2\la \eta,V\ra\right),
\end{align*}
where $\phi\in \mathcal{CM}(X\sqcup\bar{X})$.
\end{definition}

The functional $S_{\Ga}$ satisfies an obvious analog of Lemma \ref{1st-variation}. Its
Euler-Lagrange equation is the Liouville equation, so that its single non-degenerate
critical point is the hyperbolic metric on $\up\cup\lo$. Corresponding classical
action is $8\pi(2g-2)$ --- twice the hyperbolic area of $X\sqcup\bar{X}$. Similarly to
\eqref{shift-1} we have
\begin{equation} \label{shift-2}
S_{\Ga}[\phi]= \frac{i}{2}\left(\la \omega[\phi],F-\bar{F}\ra - \la
\ptheta[\phi],L-\bar{L}\ra +\la \pu, W-\bar{W}\ra\right).
\end{equation}
\begin{remark}
In the definition of $S_{\Ga}$ it is not necessary to choose a fundamental domain for
$\Ga$ in $\lo$ to be the mirror image of the fundamental domain in $\up$ since the
corresponding homology class $[\Sigma-\bar{\Sigma}]$ does not depend on the choice of
the fundamental domain of $\Ga$ in $\up\cup\lo$.
\end{remark}

\subsection{The quasi-Fuchsian case}
Let $\Ga$ be a marked, normalized, purely loxodromic quasi-Fuchsian group of genus
$g>1$. Its region of discontinuity $\Omega$ has two invariant components $\Omega_1$
and $\Omega_2$ separated by a quasi-circle $\mathcal{C}$. By definition, there exists
a quasiconformal homeomorphism $J_1$ of $\hat\CC$ with the following properties.
\begin{itemize}
\item[\textbf{QF1}]
The mapping $J_1$ is holomorphic on $\U$ and $J_1(\up)=\Omega_1$, $J_1(\lo)=\Omega_2$,
and $J_1(\RR_{\infty})=\mathcal{C}$.
\item[\textbf{QF2}] The mapping $J_1$ fixes $0,1~\text{and}~\infty$.
\item[\textbf{QF3}] The group $\tilde\Ga=J_1^{-1}\circ \Ga\circ J_1$ is Fuchsian.
\end{itemize}
Due to the normalization, any two maps satisfying \textbf{QF1}-\textbf{QF3} agree on
$\U$, so that the group $\tilde\Ga$ is independent of the choice of the map $J_1$.
Setting $X\simeq \tilde\Ga\bk\U$, we get $\tilde\Ga\bk\up\cup\lo\simeq X\sqcup\bar{X}$
and $\Ga\bk\Omega\simeq X\sqcup Y$, where $X$ and $Y$ are marked compact Riemann
surfaces of genus $g>1$ with opposite orientations. Conversely, according to Bers'
simultaneous uniformization theorem~\cite{B3}, for any pair of marked compact Riemann
surfaces $X$ and $Y$ of genus $g>1$ with opposite orientations there exists a 
unique, up to a conjugation in $\PSL(2,\CC)$, quasi-Fuchsian group $\Ga$ such that
$\Ga\bk\Omega\simeq X\sqcup Y$.

\begin{remark}
It is customary (see, e.g., \cite{A2}) to define quasi-Fuchsian groups by
requiring that the map $J_1$ is holomorphic in the lower half-plane $\lo$. We will
see in Section 4 that the above definition is somewhat more convenient.
\end{remark}
Let $\mu$ be the Beltrami coefficient for the quasiconformal map $J_1$,
\begin{displaymath}
\mu = \frac{ (J_1)_{\z}}{ (J_1)_z},
\end{displaymath}
that is, $J_1=f^\mu$ --- the unique, normalized solution of the Beltrami equation on
$\hat{\CC}$ with Beltrami coefficient $\mu$. Obviously, $\mu=0$ on $\U$. Define
another Beltrami coefficient $\hat{\mu}$ by
\begin{equation*}
  \hat{\mu}(z) =\begin{cases}
                 \ov{\mu(\z)}& \text{if $z \in \U$},\\
                 \mu(z)& \text{if $z \in \lo$}.
              \end{cases}
\end{equation*}
Since $\hat{\mu}$ is symmetric, normalized solution $f^{\hat{\mu}}$ of the Beltrami
equation
\begin{equation*}
f^{\hat{\mu}}_{\z}(z) = \hat{\mu}(z)f^{\hat{\mu}}_z(z)
\end{equation*}
is a quasiconformal homeomorphism of $\hat\CC$ which preserves $\U$ and $\lo$. The
quasiconformal map $J_2= J_1\circ (f^{\hat{\mu}})^{-1}$ is then conformal on the lower
half-plane $\lo$ and has properties similar to \textbf{QF1}-\textbf{QF3}. In
particular, $J_2^{-1}\circ \Ga\circ J_2=\hat\Ga=f^{\hat{\mu}} \circ\tilde{\Ga}\circ
(f^{\hat{\mu}})^{-1}$ is a Fuchsian group and $\hat\Ga\bk\lo\simeq Y$. Thus for a
given $\Ga$ the restriction of the map $J_2$ to $\lo$ does not depend on the choice of
$J_2$ (and hence of $J_1$). These properties can be summarized by the following
commutative diagram
\begin{equation*}
\begin{CD}
\up\cup\RR_{\infty}\cup\lo  @> J_1=f^\mu >> \Omega_1\cup \mathcal{C}\cup\Omega_2 \\
@VV f^{\hat\mu}V     @AAJ_2A \\
 \up\cup\RR_{\infty}\cup\lo @>=>> \up\cup\RR_{\infty}\cup\lo
\end{CD}
\end{equation*}
where maps $J_1, J_2~\text{and}~f^{\hat\mu}$ intertwine corresponding pairs of groups
$\Ga,\tilde\Ga~\text{and}~\hat\Ga$.
\subsubsection{Homology construction}
The map $J_1$ induces a chain map between double complexes $\KKK_{\bu,\bu}=
\SSS_\bu\otimes_{\Z\Ga}\BBB_\bu$ for the pairs $\U\cup\lo, \tilde\Ga$ and $\Omega,
\Ga$, by pushing forward chains $S_{\bullet}(\up\cup\lo)\ni c\mapsto J_1(c)\in
S_\bullet(\Omega)$ and group elements $\tilde\Ga\ni\s\mapsto J_1\circ\s\circ
J_1^{-1}\in\Ga$. We will continue to denote this chain map by $J_1$. Obviously, the
chain map $J_1$ induces an isomorphism between homology groups of corresponding total
complexes $\Tot\KKK$.

Let $\Sigma = F + L - V$ be total cycle of degree 2 representing the fundamental
class of $X$ in the total homology complex for the pair $\U, \tilde\Ga$, constructed
in the previous section, and let $\Sigma^\prime = F^\prime + L^\prime - V^\prime$ be
the corresponding cycle for $\bar{X}$. The total cycle $\Sigma(\Ga)$ of degree 2
representing fundamental class of $X\sqcup Y$ in the total complex for the pair
$\Omega, \Ga$ can be realized as a push-forward of the total cycle $\Sigma(\tilde\Ga)=
\Sigma - \Sigma^\prime$ by $J_1$,
\begin{equation*}
\Sigma(\Ga) = J_1(\Sigma(\tilde\Ga))=J_1(\Sigma) - J_1(\Sigma^\prime).
\end{equation*}
We will denote push-forwards by $J_1$ of the chains $F, L, V$ in $\up$ by $F_1, L_1,
V_1$, and push-forwards of the corresponding chains $F^\prime, L^\prime, V^\prime$ in
$\lo$ --- by $F_2, L_2, V_2$, where indices $1$ and $2$ refer, respectively, to
domains $\Omega_1$ and $\Omega_2$.

The definition of chains $W_i$ is more subtle. Namely, the quasi-circle
$\mathcal{C}$ is not generally smooth or even rectifiable, so that an arbitrary path
from an interior point of $\Omega_i$ to $p\in\mathcal{C}$ inside $\Omega_i$ is not
rectifiable either. Thus if we define $W_1$ as a push-forward by $J_1$ of $W$
constructed using arbitrary admissible paths in $\U$, the paths in $W_1$ in general
will no longer be rectifiable. The same applies to the push-forward by $J_1$ of the
corresponding chain in $\lo$. However, the definition of $\la u,W_1\ra$ uses
integration of the 1-form $u_{\s_1,\s_2}$ along the paths in $W_1$, and these paths
should be rectifiable in order that $\la u,W_1\ra$ is well-defined. The invariant
construction of such paths in $\Omega_i$ is based on the following elegant observation
communicated to us by M.~Lyubich.

Since the quasi-Fuchsian group $\Ga$ is normalized, it follows from \textbf{QF2} that
the Fuchsian group $\tilde{\Ga} = J_1^{-1}\circ \Ga\circ J_1$ is also
normalized and $\tilde{\alpha}_1\in\tilde{\Ga}$ is a dilation $\tilde{\alpha}_1\,z =
\tilde{\lambda}z$ with the axis $i\RR_{\geq 0}$ and $0<\tilde{\lambda}<1$.
Corresponding loxodromic element $\alpha_1 = J_1\circ \tilde{\alpha_1}\circ
J_1^{-1}\in\Ga$ is also a  dilation $\alpha_1\,z = \lambda z$, where $0<|\lambda|<1$.
Choose $\tilde{z}_0 \in i\RR_{\geq 0}$ and denote by $\tilde{I}=[\tilde{z}_0,0]$ the
interval on $i\RR_{\geq 0}$ with endpoints $\tilde{z}_0$ and 0 --- the attracting
fixed point of $\tilde{\alpha}_1$. Set $z_0=J_1(\tilde{z}_0)$ and $I=J_1(\tilde{I})$.
The path $I$ connects points $z_0\in\Omega_1$ and $0=J_1(0)\in\mathcal{C}$ inside
$\Omega_1$, is smooth everywhere except the endpoint $0$, and is rectifiable. Indeed,
set $\tilde{I}_0=[\tilde{z}_0,\, \tilde{\lambda}\tilde{z}_0]\subset i\RR_{>0}$ and
cover the interval $\tilde{I}$ by subintervals $\tilde{I}_n$ defined by
$\tilde{I}_{n+1} = \tilde{\alpha}_1(\tilde{I}_n)$, $n=0,1, \dots,\infty$.
Corresponding paths $I_n=J_1(\tilde{I}_n)$ cover the path $I$, and due to the property
$I_{n+1} = \alpha_1(I_n)$, which follows from \textbf{QF3}, we have
\begin{equation*}
I=\bigcup_{n=0}^\infty \alpha_1^n(I_0).
\end{equation*}
Thus
\begin{displaymath}
l(I)=\sum_{n=0}^\infty |\lambda^n|l(I_0) = \frac{l(I_0)}{1-|\lambda|}<\infty,
\end{displaymath}
where $l(P)$ denotes the Euclidean length of a smooth path $P$.

The same construction works for every $p\in\mathcal{C}\setminus\{\infty\}$ which is a
fixed point of an element in $\Ga$, and we define $\Ga$-contracting paths in
$\Omega_1$ at $p$ as follows.
\begin{definition} \label{contracting1}
Path $P$ connecting points $z\in\Omega_1$ and $p\in\mathcal{C} \setminus\{\infty\}$
inside $\Omega_1$ is called $\Ga$-contracting in $\Omega_1$ at $p$, if the following
conditions are satisfied.
\begin{itemize}
\item[\textbf{C1}] Paths $P$ is smooth except at the point $p$.
\item[\textbf{C2}] The point $p$ is a fixed point for $\Ga$.
\item[\textbf{C3}] There exists $p'\in P$ and an arc $P_0$ on the path $P$
such that the iterates $\s^n(P_0),\,n\in\N$, where $\s\in\Ga$ has $p$ as the
attracting fixed point, entirely cover the part of $P$ from the point $p'$ to the
point $p$.
\end{itemize}
\end{definition}
As in Section 2.2, we define $\Ga$-closed paths and $\Ga$-closed contracting paths in
$\Omega_1$ at $p$. Definition of $\Ga$-contracting paths in $\Omega_2$ is analogous.
Finally, we define $\Ga$-contracting paths in $\Omega$ as follows.

\begin{definition} \label{contracting2}
Path $P$ is called $\Ga$-contracting in $\Omega$, if $P=P_1\cup P_2$, where $P_1\cap
P_2 =p\in \mathcal{C}$, and $P_1\setminus\{p\}\subset\Omega_1$ and
$P_2\setminus\{p\}\subset\Omega_2$ are $\Ga$-contracting paths at $p$ in the sense of
the previous definition.
\end{definition}

$\Ga$-contracting paths are rectifiable.

\begin{lemma} \label{Paths}
Let $\Ga$ and $\Ga'$ be two marked normalized quasi-Fuchsian groups with regions of
discontinuity $\Omega$ and $\Omega'$, and let $f$ be normalized quasiconformal
homeomorphism of $\CC$ which intertwines $\Ga$ and $\Ga'$ and is smooth in $\Omega$.
Then the push-forward by $f$ of a $\Ga$-contracting path in $\Omega$ is a
$\Ga'$-contracting path in $\Omega'$.
\end{lemma}
\begin{proof}
Obvious: if $p$ is the attracting fixed point for $\s\in\Ga$, then $p'=f(p)$ is the
attracting fixed point for $\s'=f\circ \s\circ f^{-1} \in\Ga'$.
\end{proof}
Now define a chain $W$ for the Fuchsian group $\tilde{\Ga}$ by first connecting points
$P_1(1),\dots, P_g(1)$ to some point $\tilde{z}_0\in i\RR_{>0}$ by smooth paths inside
$\U$ and then connecting this point to $0$ by $\tilde{I}$. The chain $W'$ in $\lo$ is
defined similarly. Setting $W_1=J_1(W)$ and $W_2=J_1(W')$, we see that the chain $W_1
- W_2$ in $\Omega$ consists of $\Ga$-contracting paths in $\Omega$ at $0$. Connecting
$P_1(1),\dots, P_g(1)$ to $0$ by arbitrary  
$\Ga$-contracting paths at $0$ results in 1-chains which are homotopic to 
the 1-chains $W_1$ and $W_2$ in components $\Omega_1$ and $\Omega_2$ respectively. 
Finally, we define chain
$U_1=U_2$ as push-forward by $J_1$ of the corresponding chain $U=U'$ with $p=0$.

\subsubsection{Cohomology construction}
Let $\mathcal{CM}(X\sqcup Y)$ be the space of all conformal metrics
$ds^2=e^\phi|dz|^2$ on $X\sqcup Y$, which we will always identify with the affine
space of smooth real-valued functions $\phi$ on $\Omega$ satisfying \eqref{field}. For
$\phi\in \mathcal{CM}(X\sqcup Y)$ we define cochains $\omega[\phi], \theta[\phi], u,
\eta~\text{and}~\Theta$ in the total cohomology complex $\Tot\CCC$ for the pair
$\Omega,\, \Ga$ by the same formulas \eqref{omega}, \eqref{theta}, \eqref{u},
\eqref{eta} and \eqref{Theta-up}, \eqref{Theta-lo} as in the Fuchsian case, where
$p=0\in \mathcal{C}$, integration goes over $\Ga$-contracting paths at 0, and
$\tilde\s\in \tilde\Ga$ are replaced by $\s=J_1\circ\tilde \s\circ J_1^{-1}\in\Ga$.
The ordering of points on $\mathcal{C}$ used in the definition \eqref{eta} of the
constants of integration $\eta_{\s_1, \s_2}$ is defined by the orientation of
$\mathcal{C}$.
\begin{remark} \label{paths-2}
Since 1-form $u$ is closed and regular in $\Omega_1\cup\Omega_2$, it follows from
Stokes' theorem that in the definition \eqref{Theta-up} and \eqref{Theta-lo} of the
cochain $\Theta\in\CCC^{0,2}$ we can use any rectifiable path from $z$ to $0$ inside
$\Omega_1$ and $\Omega_2$ respectively.
\end{remark}
As opposed to the Fuchsian case, we can no longer guarantee that the cochain
$\omega[\phi]-\theta[\phi] - \Theta$ is a 2-cocycle in the total cohomology complex
$\Tot\CCC$. Indeed, we have, using $\de u=0$,
\begin{equation} \label{Theta-qf}
(\de\Theta)_{\s_1,\s_2,\s_3}(z)=\begin{cases}\int_{P_1}u_{\s_2,\s_3} +
(\de\eta)_{\s_1,\s_2,\s_3}=(d_1)_{\s_1,\s_2,\s_3}& \text{if $z\in\Omega_1$}, \\
\int_{P_2}u_{\s_2,\s_3} - (\de\eta)_{\s_1,\s_2,\s_3}=(d_2)_{\s_1,\s_2,\s_3}& \text{if
$z\in\Omega_2$},
\end{cases}
\end{equation}
where paths of integration $P_1$ and $P_2$ are $\Ga$-closed contracting paths
connecting points $0$ and $\s_1^{-1}(0)$ inside $\Omega_1$ and $\Omega_2$
respectively. Since the analog of Lemma~\ref{constants} does not hold in the
quasi-Fuchsian case, we can not conclude that $d_1=d_2=0$. However, $d_1, d_2
\in\CCC^{0,3}$ are $z$-independent group 3-cocycles and
\begin{equation} \label{d-d}
(d_1-d_2)_{\s_1,\s_2,\s_3}=\int_{C} u_{\s_2,\s_3} + 2(\de\eta)_{\s_1,\s_2,\s_3},
\end{equation}
where $C=P_1 - P_2$ is a loop that starts at $0$, goes to $\s_1^{-1}(0)$ inside
$\Omega_1$, continues inside $\Omega_2$ and ends at $0$. In the Fuchsian case we have
the equation \eqref{eta-integral}, which can be derived using the Stokes' theorem (see
Remark \ref{identity}). The same derivation repeats verbatim for the quasi-Fuchsian
case,
 and we get
\begin{equation*}
\int_{C} u_{\s_2,\s_3} = -2(\de\eta)_{\s_1,\s_2,\s_3},
\end{equation*}
so that $d_1=d_2$. Since $H^3(\Ga,\CC)=0$, there exists a constant 2-cochain $\kappa$
such that $\de\kappa=-d_1=-d_2$. Then $\Theta + \kappa$ is a group 2-cocycle, that is,
$\de(\Theta + \kappa)=0$. As the result, we obtain that
\begin{equation*}
\Psi[\phi]=\omega[\phi]-\theta[\phi] - \Theta -\kappa \in(\Tot\CCC)^2
\end{equation*}
is a 2-cocycle in total cohomology complex $\Tot\CCC$ for the pair $\Omega,\, \Ga$,
that is, $D\Psi[\phi]=0$.
\begin{remark}
The map $J_1$ induces a cochain map between double cohomology complexes $\Tot\CCC$ for
the pairs $\up\cup\lo,\tilde\Ga$ and $\Omega,\,\Ga$, by pulling back cochains and
group elements,
\begin{equation*}
(J_1\cdot\varpi)_{\tilde{\s}_1,\dots,\tilde{\s}_q}=J_1^\ast \varpi_{\s_1,\dots,\s_q}
\in \CCC^{p,q}(\up\cup\lo),
\end{equation*}
where $\varpi\in\CCC^{p,q}(\Omega)$ and $\tilde{\s}=J_1^{-1}\circ\s\circ J_1$. This
cochain map induces an isomorphism of the cohomology groups of corresponding total
complexes $\Tot\CCC$. The map $J_1$ also induces a natural isomorphism between the
affine spaces $\mathcal{CM}(X\sqcup Y)$ and $\mathcal{CM}(X\sqcup\bar{X})$,
\begin{equation*}
J_1\cdot\phi=\phi\circ J_1 + \log|(J_1)_z|^2\in\mathcal{CM}(X\sqcup\bar{X}),
\end{equation*}
where $\phi\in \mathcal{CM} (X\sqcup Y)$. However,
\begin{displaymath}
\left|(J_1\cdot\phi)_z\right|^2 dz\wedge d\z \neq J_1^\ast\left(\left|\phi_z \right|^2
dz\wedge d\z\right),
\end{displaymath}
and cochains $\omega[\phi], \theta[\phi], u ~\text{and}~\Theta$ for the pair
$\Omega,\, \Ga$ are not pull-backs of cochains for the pair $\up\cup\lo$, 
$\tilde{\Ga}$ corresponding to $J_1\cdot\phi\in\mathcal{CM}(X\sqcup\bar{X})$.
\end{remark}
\subsubsection{The Liouville action functional}
Discussion in the previous section justifies the following definition.
\begin{definition} The Liouville action functional $S_{\Ga}:
\mathcal{CM} (X\sqcup Y)\rightarrow \RR$ for the quasi-Fuchsian group $\Ga$ is defined
by
\begin{align*}
S_{\Ga}[\phi]= & \frac{i}{2}\la \Psi[\phi],\Sigma(\Ga)\ra = \frac{i}{2}\la
\Psi[\phi],\Sigma_1 - \Sigma_2 \ra\\ = &\frac{i}{2}\left(\la \omega[\phi],F_1-F_2\ra -
\la \theta[\phi],L_1- L_2\ra + \la \Theta + \kappa, V_1-V_2 \ra \right),
\end{align*}
where $\phi\in \mathcal{CM} (X\sqcup Y)$.
\end{definition}
\begin{remark} \label{indep}
Since $\Psi[\phi]$ is a total 2-cocycle, the Liouville action functional $S_{\Ga}$
does not depend on the choice of fundamental domain for $\Ga$ in $\Omega$, i.e on
the choice of fundamental domains $F_1$ and $F_2$ for $\Ga$ in $\Omega_1$ and
$\Omega_2$. In particular, if $\Sigma_1$ and $\Sigma_2$ are push-forwards by the map
$J_1$ of the total cycle $\Sigma$ and its mirror image $\bar{\Sigma}$, then
$\la\kappa, V_1-V_2 \ra = 0$ and we have
\begin{equation} \label{qf-action}
S_{\Ga}[\phi]= \frac{i}{2}\left(\la \omega[\phi],F_1-F_2\ra - \la
\theta[\phi],L_1-L_2\ra +\la u, W_1-W_2\ra + 2\la \eta,V_1\ra\right).
\end{equation}
In general, the constant group 2-cocycle $\kappa$ drops out from the definition for
any choice of fundamental domains $F_1$ and $F_2$ which is associated with the same
marking of $\Ga$, i.e., when the same choice of standard generators
$\alpha_1,\dots,\alpha_g,\beta_1,\dots, \beta_g$ is used both in $\Omega_1$ and in
$\Omega_2$. Indeed, in this case $V_1$ and $V_2$ have the same
$\BBB_2(\ZZ\Ga)$-structure and $\la\kappa, V_1-V_2 \ra = 0$. Moreover, since 1-form
$u$ is closed and regular in $\Omega_1\cup\Omega_2$, we can use arbitrary rectifiable
paths with endpoint 0 inside $\Omega_1$ and $\Omega_2$ in the definition of chains
$W_1$ and $W_2$ respectively.
\end{remark}
\begin{remark}
We can also define chains $W_1$ and $W_2$ by using $\Ga$-contracting paths at any
$\Ga$-fixed point $p\in\mathcal{C}\setminus\{\infty\}$. As in Remark
\ref{independence} it is easy to show that
\begin{displaymath}
\la\Theta,V_1 - V_2 \ra = \la u,W_1-W_2\ra +2\la\eta,V_1\ra
\end{displaymath}
does not depend on the choice of a fixed point $p$.
\end{remark}

As in the Fuchsian case, the Euler-Lagrange equation for the functional $S_{\Ga}$ is
the Liouville equation and the hyperbolic metric $e^{\phi_{hyp}}|dz|^2$ on $\Omega$ is
its single non-degenerate critical point. It is explicitly given by
\begin{equation} \label{qf-hyp}
e^{\phi_{hyp}(z)} =
          \frac{|(J_i^{-1})^\prime(z)|^2}{(\im J_i^{-1}(z))^2}\;\;\,\text{if}\;\;\,
          z \in\Omega_i,\;\, i=1,2.
\end{equation}

\begin{remark}
Corresponding classical action $S_{\Ga}[\phi_{hyp}]$ is no longer twice the hyperbolic
area of $X\sqcup Y$, as it was in the Fuchsian case, but rather non-trivially depends
on $\Ga$. This is due to the fact that in the quasi-Fuchsian case the $(1,1)$-form
$\omega[\phi_{hyp}]$ on $\Omega$ is not a $(1,1)$-tensor for $\Ga$, as it was in the
Fuchsian case.
\end{remark}

Similarly to \eqref{shift-2} we have
\begin{equation} \label{shift-3}
S_{\Ga}[\phi]= \frac{i}{2}\left(\la \omega[\phi],F_1 - F_2\ra - \la \ptheta[\phi],L_1
- L_2\ra + \la \pu, W_1 - W_2\ra\right),
\end{equation}
where  $F_1$ and $F_2$ are fundamental domains for the marked group $\Ga$ in
$\Omega_1$ and $\Omega_2$ respectively.

\section{Deformation theory}
\subsection{The deformation space}
Here we collect the basic facts from deformation theory of Kleinian groups (see, e.g.,
\cite{A2,B1,B2,Kra2}). Let $\Ga$ be a non-elementary, finitely generated purely
loxodromic Kleinian group, let $\Omega$ be its region of discontinuity, and let
$\Lambda=\hat{\CC}\setminus \Omega$ be its limit set. The deformation space $\D(\Ga)$
is defined as follows. Let $\mathcal{A}^{-1,1}(\Ga)$ be the space of Beltrami
differentials for $\Ga$ --- the Banach space of $\mu\in L^{\infty}(\CC)$ satisfying
\begin{equation*}
    \mu(\s(z))\; \frac{\ov{\s'(z)}}{\s'(z)} = \mu (z)~\text{for all}~\s \in \Ga,
\end{equation*}
and
\begin{equation*}
\left. \mu \right|_{\Lambda}=0.
\end{equation*}
Denote by $\mathcal{B}^{-1,1}(\Ga)$ the open unit ball in $\mathcal{A}^{-1,1}(\Ga)$
with respect to $\parallel \cdot\parallel_{\infty}$ norm,
\begin{equation*}
\parallel \mu \parallel_{\infty} = \sup_{z\in \C} |\mu(z)| <1.
\end{equation*}
For each Beltrami coefficient $\mu\in\mathcal{B}^{-1,1}(\Ga)$ there exists a unique
homeomorphism $f^{\mu}: \hat{\C} \rightarrow \hat{\C}$ satisfying the Beltrami
equation
\begin{equation*}
f^{\mu}_{\z} = \mu f^{\mu}_z
\end{equation*}
and fixing the points $0,1~\text{and}~\infty$. Set $\Ga^{\mu} = f^{\mu}\circ\Ga
\circ(f^{\mu})^{-1}$ and define
\begin{equation*}
\D(\Ga)=\mathcal{B}^{-1,1}(\Ga)/ \thicksim\,,
\end{equation*}
where $\mu \thicksim \nu$ if and only if $f^\mu=f^\nu$ on $\Lambda$, which is
equivalent to the condition 
$f^\mu\circ\s\circ (f^\mu)^{-1}=f^\nu\circ\s\circ (f^\nu)^{-1}$ for all
$\s\in\Ga$.

Similarly, if $\Delta$ is a union of invariant components of $\Ga$, the deformation
space $\D(\Ga,\Delta)$ is defined using Beltrami coefficients supported on $\Delta$.

By Ahlfors finiteness theorem $\Omega$ has finitely many non-equivalent components
$\Omega_1,\ldots,\Omega_n$. Let $\Ga_i$ be the stabilizer subgroup of the component
$\Omega_i$, $\Ga_i=\{\s\in\Ga\,|\,\s(\Omega_i) = \Omega_i\}$ and let $X_i\simeq
\Ga_i\bk\Omega_i$ be the corresponding compact Riemann surface of genus $g_i>1$,
$i=1,\ldots,n$. The decomposition
\begin{equation*}
\Ga \bk \Omega = \Ga_1\bk\Omega_1\sqcup\cdots\sqcup\Ga_n\bk\Omega_n
\end{equation*}
establishes the isomorphism \cite{Kra2}
\begin{equation*} 
\D(\Ga)\simeq \D(\Ga_1,\Omega_1)\times\cdots\times \D(\Ga_n,\Omega_n).
\end{equation*}
\begin{remark}
When $\Ga$ is a purely hyperbolic Fuchsian group of genus $g>1$, $\D(\Ga,\up) =
\T(\Ga)$ --- the \Te space of $\Ga$. Every conformal bijection $\Ga \bk \up\rightarrow
X$ establishes isomorphism between $\T(\Ga)$ and $\T(X)$, the \Te space of marked
Riemann surface $X$. Similarly, $\D(\Ga,\lo)=\bar\T(\Ga)$, the mirror image of
$\T(\Ga)$ --- the complex manifold which is complex conjugate to $\T(\Ga)$.
Correspondingly, $\Ga \bk \lo\rightarrow \bar{X}$ establishes the isomorphism
$\bar\T(\Ga)\simeq \T(\bar{X})$, so that
\begin{equation*}
\D(\Ga)\simeq \T(X)\times \T(\bar{X}).
\end{equation*}
The deformation space $\D(\Ga)$ is ``twice larger'' than the \Te space $\T(\Ga)$
because its definition uses all Beltrami coefficients $\mu$ for $\Ga$, and not only
those satisfying the reflection property $\mu(\z)= \ov{\mu(z)}$, used in the
definition of $\T(\Ga)$.
\end{remark}

The deformation space $\D(\Ga)$ has a natural structure of a complex manifold,
explicitly described as follows (see, e.g., \cite{A2}). Let $\mathcal{H}^{-1,1}(\Ga)$
be the Hilbert space of Beltrami differentials for $\Ga$ with the following scalar
product
\begin{equation}\label{pairing}
(\mu_1, \mu_2) = \iint\limits_{\Ga \bk \Omega} \mu_1 \bar{\mu}_2 \rho =
\iint\limits_{\Ga \bk \Omega} \mu_1(z) \ov{\mu_2(z)} \rho(z)\,d^2z,
\end{equation}
where $\mu_1, \mu_2 \in \mathcal{H}^{-1,1}(\Ga)$ and $\rho=e^{\phi_{hyp}}$ is the
density of the hyperbolic metric on $\Ga \bk \Omega$. Denote by $\Omega^{-1,1}(\Ga)$
the finite-dimensional subspace of harmonic Beltrami differentials with respect to the
hyperbolic metric. It consists of $\mu\in\mathcal{H}^{-1,1}(\Ga)$ satisfying
\begin{equation*}
\pa_z(\rho \mu) = 0.
\end{equation*}
The complex vector space $\Omega^{-1,1}(\Ga)$ is identified with the holomorphic
tangent space to $\D(\Ga)$ at the origin. Choose a basis $\mu_1,\dots, \mu_d$ for
$\Omega^{-1,1}(\Ga)$, let $\mu = \vep_1 \mu_1 + \cdots + \vep_d \mu_d$, and let
$f^\mu$ be the normalized solution of the Beltrami equation. Then the correspondence
$(\vep_1, \dots, \vep_d) \mapsto \Ga^{\mu}= f^\mu \circ \Ga \circ (f^\mu)^{-1}$
defines complex coordinates in a neighborhood of the origin in $\D(\Ga)$, called Bers
coordinates. The holomorphic cotangent space to $\D(\Ga)$ at the origin can be
naturally identified with the vector space $\Omega^{2,0}(\Ga)$ of holomorphic
quadratic differentials --- holomorphic functions $q$ on $\Omega$ satisfying
\begin{equation*}
q(\s z)\s'(z)^2=q(z)~\text{for all}~\s\in\Ga.
\end{equation*}
The pairing between holomorphic cotangent and tangent spaces to $\D(\Ga)$ at the
origin is given by
\begin{equation*}
q(\mu) = \iint\limits_{\Ga\bk\Omega} q\mu = \iint\limits_{\Ga \bk \Omega}
q(z)\mu(z)\,d^2z.
\end{equation*}

There is a natural isomorphism $\Phi^\mu$ between the deformation spaces $\D(\Ga)$ and
$\D(\Ga^{\mu})$, which maps $\Ga^{\nu} \in \D(\Ga)$ to $(\Ga^{\mu})^{\lambda} \in
\D(\Ga^{\mu})$, where, in accordance with $f^\nu=f^\lambda \circ f^\mu$,
\[
\lambda = \left(\frac{\nu - \mu}{ 1- \nu \bar{\mu}}
\frac{f_z^{\mu}}{\bar{f}_{\z}^{\mu}} \right)\circ (f^{\mu})^{-1}.
\]
The isomorphism $\Phi^\mu$ allows us to identify the holomorphic tangent space to
$\D(\Ga)$ at $\Ga^{\mu}$ with the complex vector space $\Omega^{-1,1}(\Ga^{\mu})$, and
holomorphic cotangent space to $\D(\Ga)$ at $\Ga^{\mu}$ with the complex vector space
$\Omega^{2,0}(\Ga^{\mu})$. It also allows us to introduce the Bers coordinates in the
neighborhood of $\Ga^\mu$ in $\D(\Ga)$, and to show directly that these coordinates
transform complex-analytically. For the de Rham differential $d$ on $\D(\Ga)$ we
denote by $d=\pa + \bar\pa$ the decomposition into $(1,0)$ and $(0,1)$ components.

The differential of isomorphism $\Phi^\mu: \D(\Ga)\simeq \D(\Ga^{\mu})$ at
$\nu=\mu$ is given by the linear map $D^{\mu} : \Omega^{-1,1}(\Ga) \rightarrow
\Omega^{-1,1} (\Ga^{\mu})$,
\begin{equation*}
\nu \mapsto D^{\mu} \nu = P_{-1,1}^{\mu}\left[ \left(\frac{\nu}{ 1- |\mu|^2}
\frac{f_z^{\mu}}{\bar{f}_{\z}^{\mu}} \right)\circ (f^{\mu})^{-1} \right],
\end{equation*}
where $P_{-1,1}^{\mu}$ is orthogonal projection from
$\mathcal{H}^{-1,1}(\Ga^{\mu})$ to $\Omega^{-1,1}(\Ga^{\mu})$. The map $D^\mu$ allows
to extend a tangent vector $\nu$ at the origin of $\D(\Ga)$ to a local vector field
$\pa/\pa\vep_{\nu}$ on the coordinate neighborhood of the origin,
\begin{equation*}
\left.\frac{\pa}{\pa \vep_{\nu}}\right|_{\Ga^{\mu}}  = D^{\mu}
\nu\in\Omega^{-1,1}(\Ga^{\mu}).
\end{equation*}

The scalar product \eqref{pairing} in $\Omega^{-1,1}(\Ga^{\mu})$ defines a Hermitian
metric on the deformation space $\D(\Ga)$. This metric is called the Weil-Petersson
metric and it is \Ka. We denote its symplectic form by $\omega_{WP}$,
\begin{equation*}
\omega_{WP} \left.\left(\frac{\pa}{\pa \vep_\mu}, \frac{\pa}{\pa
\bar{\vep}_\nu}\right) \right|_{\Ga^{\lambda}} = \frac{i}{2} \left( D^{\lambda} \mu,
D^{\lambda} \nu \right),\quad \mu, \nu \in \Omega^{-1,1}(\Ga).
\end{equation*}
\subsection{Variational formulas}
Here we collect necessary variational formulas. Let $l$ and $m$ be integers. A tensor
of type $(l,m)$ for $\Ga$ is a $C^\infty$-function $\omega$ on $\Omega$ satisfying
\begin{equation*}
\omega(\s z)\s'(z)^l \ov{\s^\prime(z)}^m=\omega(z)~\text{for all}~\s\in\Ga.
\end{equation*}
Let $\omega^{\vep}$ be a smooth family of tensors of type $(l,m)$ for $\Ga^{\vep\mu}$,
where $\mu \in \Omega^{-1,1} (\Ga)$ and $\vep \in \C$ is sufficiently small. Set
\begin{equation*}
(f^{\vep\mu})^* (\omega^{\vep}) = \omega^{\vep} \circ f^{\vep\mu} (f_z^{\vep\mu})^l
(\bar{f}_{\z}^{\vep\mu})^m,
\end{equation*}
which is a tensor of type $(l,m)$ for $\Ga$ --- a pull-back of the tensor
$\omega^\vep$ by $f^{\vep\mu}$. The Lie derivatives of the family $\omega^{\vep}$
along the vector fields $\pa/\pa \vep_{\mu}$ and $\pa/\pa \bar{\vep}_{\mu}$ are
defined in the standard way,
\begin{equation*}
L_{\mu} \omega = \left.\frac{\pa}{\pa \vep} \right|_{\vep=0} (f^{\vep\mu})^*
(\omega^{\vep})~\text{and}~L_{\bar{\mu}} \omega = \left.\frac{\pa}{\pa \bar{\vep}}
\right|_{\vep=0}(f^{\vep\mu})^* (\omega^{\vep}).
\end{equation*}

When $\omega$ is a function on $\D(\Ga)$ --- a tensor of type $(0,0)$, Lie derivatives
reduce to directional derivatives
\begin{equation*}
L_\mu\omega = \pa\omega(\mu)~\text{and}~L_{\bar\mu}\omega=\bar\pa\omega(\bar\mu)
\end{equation*}
--- the evaluation of 1-forms $\pa\omega$ and $\bar\pa\omega$ on tangent
vectors $\mu$ and $\bar\mu$.

For the Lie derivatives of vector fields $\nu^{\ep \mu} = D^{\ep \mu} \nu$ we get
\cite{Wol} that $L_\mu\nu=0$ and $L_{\bar\mu}\nu$ is orthogonal to
$\Omega^{-1,1}(\Ga)$. In other words,
\begin{equation*}
\left[\frac{\pa}{\pa \vep_{\mu}},\frac{\pa}{\pa \vep_{\nu}}\right]=
\left[\frac{\pa}{\pa \vep_{\mu}}, \frac{\pa}{\pa \bar{\vep}_{\nu}}\right]=0
\end{equation*}
at the point $\Ga$ in $\D(\Ga)$.

For every $\Ga^\mu\in \D(\Ga)$, the density $\rho^\mu$ of the hyperbolic metric on
$\Omega^\mu$ is a $(1,1)$-tensor for $\Ga^\mu$. Lie derivatives of the smooth family
of $(1,1)$-tensors $\rho$ parameterized by $\D(\Ga)$ are given by the following lemma
of Ahlfors.
\begin{lemma} For every $\mu\in\Omega^{-1,1}(\Ga)$
\begin{equation*}
L_\mu\rho=L_{\bar\mu}\rho=0.
\end{equation*}
\end{lemma}
\begin{proof}
Let $\Omega_1,\ldots,\Omega_n$ be the maximal set of non-equivalent components of
$\Omega$ and let $\Ga_1,\ldots,\Ga_n$ be the corresponding stabilizer groups,
\begin{equation*}
\Ga \bk \Omega=\Ga_1\bk\Omega_1\sqcup\cdots \sqcup\Ga_n\bk\Omega_n\simeq
X_1\sqcup\cdots\sqcup X_n.
\end{equation*}
For every $\Omega_i$ denote by $J_i:\up\rightarrow\Omega_i$ the corresponding covering
map and by $\tilde{\Ga}_i$ --- the Fuchsian model of group $\Ga_i$, characterized by
the condition $\tilde{\Ga}_i\bk\up\simeq \Ga_i\bk\Omega_i\simeq X_i$ (see, e.g.,
\cite{Kra2}).

Let $\mu\in\Omega^{-1,1}(\Ga)$. For every component $\Omega_i$ the quasiconformal map
$f^{\vep\mu}$ gives rise to the following commutative diagram
\begin{equation} \label{cd}
\begin{CD}
\up @>F^{\vep\hat{\mu}_i}>> \up \\ @VVJ_iV     @VVJ_i^{\vep\mu}V \\ \Omega_i
@>f^{\vep\mu}>> \Omega_i^{\vep\mu}
\end{CD}
\end{equation}
where $F^{\vep\hat{\mu}_i}$ is the normalized quasiconformal homeomorphism of $\up$
with Beltrami differential $\hat\mu_i=J_i^*\mu$ for the Fuchsian group
$\tilde{\Ga}_i$. Let $\hat\rho$ be the density of the hyperbolic metric on $\up$; it
satisfies $\hat\rho=J_i^*\rho$, where $\rho$ is the density of the hyperbolic metric
on $\Omega_i$. Therefore, Beltrami differential $\hat{\mu}_i$ is harmonic with respect
to the hyperbolic metric on $\up$. It follows from the commutativity of the diagram
that
\begin{equation*}
(f^{\vep\mu})^* \rho^{\vep\mu}=((J^{\vep\mu}_i)^{-1}\circ f^{\vep\mu})^*
\hat\rho=(F^{\vep\hat{\mu}_i}\circ J^{-1}_i)^\ast\hat\rho=
 (J^{-1}_i)^\ast(F^{\vep\hat{\mu}_i})^\ast\hat\rho.
\end{equation*}
Now the assertion of the lemma reduces to
\begin{equation*}
\left.\frac{\pa}{\pa\vep}\right|_{\vep=0}(F^{\vep\hat{\mu}_i})^\ast\hat\rho=0,
\end{equation*}
which is the classical result of Ahlfors~\cite{Ahl}.
\end{proof}

Set
\begin{equation*}
\dot{f} = \left.\frac{\pa}{\pa \vep} \right|_{\vep=0} f^{\vep\mu},
\end{equation*}
then
\begin{equation} \label{qc-var}
\dot{f}(z)=-\frac{1}{\pi}\iint\limits_{\CC}\frac{z(z-1)\mu(w)}{(w-z)w(w-1)}\, d^2w.
\end{equation}
We have
\begin{equation*}
\dot{f}_{\z} = \mu
\end{equation*}
and also
\begin{equation*}
\left.\frac{\pa}{\pa \bar{\vep}} \right|_{\vep=0} f^{\vep\mu} = 0.
\end{equation*}

As it follows from Ahlfors lemma
\begin{equation*}
\left.\frac{\pa}{\pa \vep} \right|_{\vep=0} \left(\rho^{\vep\mu} \circ
f^{\vep\mu}\,|f^{\vep\mu}_z|^2 \right) = 0.
\end{equation*}
Using $\rho=e^{\phi_{hyp}}$ and the fact that $f^{\vep\mu}$ depends 
holomorphically on $\vep$, we get
\begin{equation} \label{Ahlfors1}
\left.\frac{\pa}{\pa \vep} \right|_{\vep=0} \left( \phi^{\vep\mu}_{hyp} \circ
f^{\vep\mu} \right) = - \dot{f}_z.
\end{equation}
Differentiation with respect to $z$ and $\z$ yields
\begin{equation} \label{Ahlfors2}
\left.\frac{\pa}{\pa \vep} \right|_{\vep=0} \left( \left(\phi^{\vep\mu}_{hyp}\right)_z
\circ f^{\vep\mu} f^{\vep\mu}_z \right) = - \dot{f}_{zz},
\end{equation}
and
\begin{equation} \label{Ahlfors3}
\left.\frac{\pa}{\pa \vep} \right|_{\vep=0} \left(
\left(\phi^{\vep\mu}_{hyp}\right)_{\z} \circ f^{\vep\mu} \bar{f}^{\vep\mu}_{\z}
\right) = - \left((\phi_{hyp})_z \dot{f}_{\z} + \dot{f}_{z\z}\right).
\end{equation}
For $\s\in\Ga$ set $\s^{\vep\mu} = f^{\vep\mu}\circ \s \circ (f^{\vep\mu})^{-1} \in
\Ga^{\vep\mu}$. We have
\begin{equation*}
(\s^{\vep\mu})^\prime \circ f^{\vep\mu} f^{\vep\mu}_z = f^{\vep\mu}_z \circ \s \; \s',
\end{equation*}
and
\begin{equation*}
\log|(\s^{\vep\mu})'\circ f^{\vep\mu}|^2 + \log|f^{\vep\mu}_z|^2 = \log|f^{\vep\mu}_z
\circ \s|^2 + \log|\s'|^2.
\end{equation*}
Therefore
\begin{equation} \label{Ahlfors4}
\left.\frac{\pa}{\pa \vep} \right|_{\vep=0} \left( \log|(\s^{\vep\mu})'\circ
f^{\vep\mu}|^2 \right) = \dot{f}_z \circ \s - \dot{f}_z,
\end{equation}
and, differentiating with respect to $z$,
\begin{equation} \label{Ahlfors5}
\left.\frac{\pa}{\pa \vep} \right|_{\vep=0}
\left(\frac{(\s^{\vep\mu})''}{(\s^{\vep\mu})'} \circ f^{\vep\mu} f^{\vep\mu}_z \right)
= \dot{f}_{zz} \circ \s  \; \s' - \dot{f}_{zz}.
\end{equation}

Denote by
\begin{equation*}
\mathcal{S}(h) = \left(\frac{h_{zz}}{h_z}\right)_z - \frac{1}{2}
\left(\frac{h_{zz}}{h_z} \right)^2=
\frac{h_{zzz}}{h_z}-\frac{3}{2}\left(\frac{h_{zz}}{h_z}\right)^2
\end{equation*}
the Schwarzian derivative of the function $h$.
\begin{lemma} \label{formula}
Set
\begin{equation*}
\dot{\s}=\left.\frac{\pa}{\pa\vep}\right|_{\vep=0}\s^{\vep\mu},~\s\in\Ga.
\end{equation*}
Then for all $\s\in\Ga$
\begin{equation} \tag{i}
\dot{f}_z\circ\s - \dot{f}_z = \frac{\dot{f}\circ\s \s^{\prime\prime}} {(\s^\prime)^2}
+ \left(\frac{\dot{\s}}{\s^\prime}\right)^\prime,
\end{equation}
and is well-defined on the limit set $\Lambda$. Also we have
\begin{align}
\dot{f}_{z\z} \circ \s \;\ov{\s'} - \dot{f}_{z\z} = &
\frac{\s^{\prime\prime}}{\s^\prime}\dot{f}_{\z}, \tag{ii} \\ \dot{f}_{zz} \circ \s
\;\s' - \dot{f}_{zz} = & \frac{1}{2} (\dot{f}_z \circ \s + \dot{f}_z) \frac{\s''}{\s'}
-\frac{2\dot{c}}{cz+d},\,~\text{for all}~\s\in\Ga. \tag{iii}
\end{align}
\end{lemma}
\begin{proof}
To prove formula (i), consider the equation
\begin{equation} \label{dot-gamma}
\dot{f}\circ\s=\dot{\s} + \s^\prime\dot{f},
\end{equation}
which follows from  $\s^{\vep\mu}\circ f^{\vep\mu}  = f^{\vep\mu}\circ \s$.
Differentiating with respect to $z$ gives (i). Since $\dot{f}$ is a homeomorphism of
$\hat{\CC}$ and $\dot{\s}/\s^\prime$ is a quadratic polynomial in $z$, formula (i)
shows that $\dot{f}_z\circ\s - \dot{f}_z$ is well-defined on $\Lambda$.

The formula (ii) immediately follows from $\dot{f}_{\z}=\mu$ and
\begin{equation*}
\mu\circ\s \frac{\ov{\s\prime}}{\s^\prime}=\mu,~\s\in\Ga.
\end{equation*}

To derive formula (iii), twice differentiating \eqref{dot-gamma} with respect to $z$
we obtain
\begin{equation*}
\dot{f}_z\circ \s \s^\prime = \dot{\s}^\prime + \s^{\prime\prime} \dot{f} +
\s^\prime\dot{f}_z,
\end{equation*}
and
\begin{equation*}
\s^\prime(\dot{f}_{zz}\circ \s \s^\prime -\dot{f}_{zz}) = \dot{\s}^{\prime\prime} +
\s^{\prime\prime\prime} \dot{f} + 2 \s^{\prime\prime}\dot{f}_z -\dot{f}_z\circ\s
\s^{\prime\prime}.
\end{equation*}
Since
\begin{equation*}
\s^{\prime\prime\prime}=\frac{3}{2}\frac{(\s^{\prime\prime})^2}{\s^\prime},
\end{equation*}
as it follows from $\mathcal{S}(\s)=0$, we can eliminate $\s^{\prime\prime} \dot{f}$
from the two formulas above and obtain
\begin{equation*}
\dot{f}_{zz} \circ \s \;\s' - \dot{f}_{zz} = \frac{1}{2} (\dot{f}_z \circ \s +
\dot{f}_z) \frac{\s''}{\s'} + \frac{\dot{\s}^{\prime\prime}}{\s^\prime}
-\frac{3}{2}\frac{\s^{\prime\prime}\dot{\s}^\prime} {(\s^\prime)^2}.
\end{equation*}
Using $2c = -\s^{\prime\prime}/(\s^\prime)^{3/2}$, we see that the last two terms in
this equations are equal to $-2\dot{c}/(cz+d)$, which proves the lemma. \end{proof}

Finally, we present the following formulas by Ahlfors~\cite{Ahl}. Let
$F^{\vep\hat\mu}$ be the quasiconformal homeomorphism of $\up$ with Beltrami
differential $\hat{\mu}$ for the Fuchsian group $\Ga$. If $\hat\mu$ is harmonic on
$\up$ with respect to the hyperbolic metric, then
\begin{align}
\left.\frac{\pa}{\pa \vep} \right|_{\vep=0} F_{zzz}^{\vep \hat{\mu}}(z) & = 0
\label{A1}, \\ \left.\frac{\pa}{\pa \bar{\vep}} \right|_{\vep=0} F_{zzz}^{\vep
\hat{\mu}}(z)& = -\frac{1}{2}\hat{\rho} \ov{\hat{\mu}(z)}. \label{A2}
\end{align}

\section{Variation of the classical action}
\subsection{Classical action}
Let $\Ga$ be a marked, normalized, purely loxodromic quasi-Fuchsian group of genus $g>1$
with region of discontinuity $\Omega=\Omega_1\cup\Omega_2$, let $X\sqcup
Y\simeq\Ga\bk\Omega$ be corresponding marked Riemann surfaces with opposite
orientations and let
\begin{equation*}
\D(\Ga)\simeq \D(\Ga,\Omega_1) \times \D(\Ga,\Omega_2)
\end{equation*}
be the deformation space of $\Ga$. Spaces $\D(\Ga,\Omega_1)$ and
$\D(\Ga,\Omega_2)$ are isomorphic to the \Te spaces $\T(X)$ and $\T(Y)$ --- they are
their quasi-Fuchsian models which use Bers' simultaneous uniformization of the varying
Riemann surface in $\T(X)$ and fixed $Y$ and, respectively, fixed $X$ and the varying
Riemann surface in $\T(Y)$. Therefore,
\begin{equation} \label{d-decomposition}
\D(\Ga)\simeq \T(X)\times \T(Y).
\end{equation}
Denote by $\Pro(\Ga)\rightarrow \D(\Ga)$ corresponding affine bundle of projective
connections, modeled over the holomorphic cotangent bundle of $\D(\Ga)$. We have
\begin{equation} \label{p-decomposition}
\Pro(\Ga)\simeq \Pro(X)\times \Pro(Y).
\end{equation}

For every $\Ga^\mu\in\D(\Ga)$ denote by $S_{\Ga^\mu}=S_{\Ga^\mu}[\phi_{hyp}]$ the
classical Liouville action. It follows from the results in Section 2.3.3 that
$S_{\Ga^\mu}$ gives rise to a well-defined real-valued function $S$ on $\D(\Ga)$.
Indeed, if $\mu\boldsymbol{\sim}\nu$, then corresponding total cycles 
$f^\mu(\Sigma(\Ga))$ and $f^\nu(\Sigma(\Ga))$ represent the same class in the 
total homology complex
$\Tot\KKK$ for the pair $\Omega^\mu,\,\Ga^\mu$, so that
\begin{equation*}
\la \Psi\left[\phi_{hyp}\right], f^\mu(\Sigma(\Ga)\ra = \la
\Psi\left[\phi_{hyp}\right], f^\nu(\Sigma(\Ga)\ra.
\end{equation*}
Moreover, real-analytic dependence of solutions of Beltrami equation on parameters
ensures that classical action $S$ is a real-analytic function on $\D(\Ga)$.

To every $\Ga'\in \D(\Ga)$ with the region of discontinuity $\Omega'$ there
corresponds a pair of marked Riemann surfaces $X'$ and $Y'$ simultaneously uniformized
by $\Ga'$, $X'\sqcup Y'\simeq \Ga'\bk\Omega'$. Set $S(X',Y')=S_{\Ga'}$ and denote by
$S_Y$ and $S_X$ restrictions of the function $S:\D(\Ga)\rightarrow\RR$ onto $\T(X)$
and $\T(Y)$ respectively. Let $\iota$ be the complex conjugation and let
$\bar\Ga=\iota(\Ga)$ be the quasi-Fuchsian group complex conjugated to $\Ga$. The
correspondence $\mu\mapsto \iota\circ\mu\circ\iota$ establishes complex-analytic
anti-isomorphism
\begin{displaymath}
\D(\Ga)\simeq \D(\bar\Ga)\simeq \T(\bar{Y})\times \T(\bar{X}).
\end{displaymath}
The classical Liouville action has the symmetry property
\begin{equation} \label{symmetry}
S(X',Y')=S(\bar{Y}',\bar{X}').
\end{equation}

For every $\phi\in\mathcal{CM}(\Ga\bk\Omega)$ set
\begin{equation*}
\vartheta[\phi]=2\phi_{zz} - \phi^2_z.
\end{equation*}
It follows from the Liouville equation that $\vartheta=\vartheta[\phi_{hyp}]
\in\Omega^{2,0}(\Ga)$, i.e., is a holomorphic quadratic differential for $\Ga$. It
follows from \eqref{qf-hyp} that
\begin{equation} \label{set}
\vartheta(z)=\begin{cases} 2\mathcal{S}\left({J^{-1}_1}\right)(z)~\text{if}~
z\in\Omega_1, \\ 2\mathcal{S}\left({J^{-1}_2}\right)(z)~\text{if}~ z\in\Omega_2.
\end{cases}
\end{equation}
Define a $(1,0)$-form $\boldsymbol{\vartheta}$ on the deformation space $\D(\Ga)$ by
assigning to every $\Ga'\in \D(\Ga)$ corresponding $\vartheta[\phi'_{hyp}]\in
\Omega^{2,0}(\Ga')$ --- a vector in the holomorphic cotangent space to $\D(\Ga)$ at
$\Ga'$.

For every $\Ga'\in\D(\Ga)$ let $P_{F}$ and $P_{QF}$ be Fuchsian and quasi-Fuchsian
projective connections on $X'\sqcup Y'\simeq \Ga'\bk\Omega'$, defined by the coverings
$\pi_F: \up\cup\lo\rightarrow X'\sqcup Y'$  and $\pi_{QF}:
\Omega_1\cup\Omega_2\rightarrow X'\sqcup Y'$ respectively.  We will continue to denote
corresponding sections of the affine bundle $\Pro(\Ga)\rightarrow\D(\Ga)$ by $P_F$ and
$P_{QF}$ respectively. The difference $P_F - P_{QF}$ is a $(1,0)$-form on $\D(\Ga)$.
\begin{lemma} On the deformation space $\D(\Ga)$,
\begin{equation*}
\boldsymbol{\vartheta} = 2(P_F - P_{QF}).
\end{equation*}
\end{lemma}
\begin{proof}
Consider the following commutative diagram
\begin{equation*}
\begin{CD}
\up\cup\lo @>J>> \Omega_1\cup\Omega_2 \\ @VV\pi_F V     @VV\pi_{QF} V \\ X\sqcup Y
@>=>> X\sqcup Y,
\end{CD}
\end{equation*}
where the covering map $J$ is equal to the map $J_1$ on component $\up$ and to the map
$J_2$ on component $\lo$. As explained in the Introduction,
$P_{F}=\mathcal{S}(\pi_{F}^{-1})$ and $P_{QF}= \mathcal{S}(\pi_{QF}^{-1})$, and it
follows from the property \textbf{SD1} and commutativity of the diagram that
\begin{equation*}
\left(\mathcal{S}\left(\pi_{F}^{-1}\right) - \mathcal{S}\left(\pi_{QF}^{-1}\right)
\right)\circ\pi_{QF} \,(\pi_{QF}^\prime)^2 = \mathcal{S}\left(J^{-1} \right).
\end{equation*}
\end{proof}
\subsection{First variation}
Here we compute the $(1,0)$-form $\pa S$ on $\D(\Ga)$.
\begin{theorem}\label{variation1} On the deformation space $\D(\Ga)$,
\begin{equation*}
\pa S = 2(P_F - P_{QF}).
\end{equation*}
\end{theorem}
\begin{proof}
It is sufficient to prove that for every $\mu\in\Omega^{-1,1}(\Ga)$
\begin{equation} \label{1st-var}
L_{\mu} S  = \boldsymbol\vartheta(\mu)=\iint\limits_{\Ga\bk\Omega}\vartheta\mu.
\end{equation}
Indeed, using the isomorphism $\Phi^\nu: \D(\Gamma)\rightarrow \D(\Gamma^{\nu})$, it
is easy to see that variation formula \eqref{1st-var} is valid at every point
$\Gamma^\nu\in \D(\Gamma)$ if it is valid at the origin. The actual computation of
$L_\mu S$ is quite similar to that in~\cite{ZT87b} for the case of Schottky groups,
with the clarifying role of homological algebra.

Let $\tilde{\Ga}$ be the Fuchsian group corresponding to $\Ga$ and let $\Sigma=F+L-V$
be the corresponding total cycle of degree 2 representing the fundamental class of $X$
in the total complex $\Tot\KKK$ for the pair $\up, \tilde{\Ga}$. As in Section 2.3.1,
set $\Sigma(\Ga)=J_1(\Sigma-\bar\Sigma)$. The corresponding total cycle for the pair
$\Omega^{\vep\mu}, \Ga^{\vep\mu}=f^{\vep\mu}\circ\Ga\circ(f^{\vep\mu})^{-1}$ can be
chosen as $\Sigma(\Ga^{\vep\mu}) = f^{\vep\mu}(\Sigma(\Ga))$. According to Remark
\ref{indep},
\begin{displaymath}
S_{\Ga^{\vep\mu}}
 =  \frac{i}{2}\left\langle \Psi\left[\phi_{hyp}^{\vep\mu}\right],
f^{\vep\mu}(\Sigma(\Ga)) \right\rangle.
\end{displaymath}
Moreover, as it follows from Lemma \ref{Paths}, we can choose
$\Ga^{\vep\mu}$-contracting at 0 paths of integration in the definition of
$\Theta^{\vep\mu}$ or, equivalently, paths in the definition of $W^{\vep\mu}_1 -
W^{\vep\mu}_2$, to be the push-forwards by $f^{\vep\mu}$ of the corresponding
$\Ga$-contracting at 0 paths. Denoting $\omega^{\vep\mu} =
\omega\left[\phi_{hyp}^{\vep\mu}\right],\,
\ptheta^{\vep\mu}=\ptheta\left[\phi_{hyp}^{\vep\mu}\right]$, and using \eqref{shift-3}
we have
\begin{displaymath}
S_{\Ga^{\vep\mu}}
 =  \frac{i}{2}\left(\la \omega^{\vep\mu},
F_1^{\vep\mu} - F_2^{\vep\mu} \ra - \la \check{\theta}^{\vep\mu}, L_1^{\vep\mu} -
L_2^{\vep\mu} \ra + \la \check{u}^{\vep\mu}, W_1^{\vep\mu} - W_2^{\vep\mu}\ra\right).
\end{displaymath}
Changing variables and formally differentiating under the integral sign in the term
$\la \check{u}^{\vep\mu}, W_1^{\vep\mu} - W_2^{\vep\mu}\ra$, we obtain
\begin{align*}
L_\mu S
 & = \left.\frac{\pa}{\pa\vep}\right|_{\vep=0}S_{\Ga^{\vep\mu}} \\
& = \frac{i}{2}\left(\langle L_\mu\omega, F_1 - F_2 \rangle - \langle L_\mu\ptheta,
L_1 - L_2 \rangle + \langle L_\mu \pu, W_1 - W_2  \rangle\right).
\end{align*}
We will justify this formula at the end of the proof. Here we observe that though
$\omega^{\vep\mu},\, \ptheta^{\vep\mu}$ and $\pu^{\vep\mu}$ are not tensors for
$\Ga^{\vep\mu}$, they are differential forms on $\Omega^{\vep\mu}$ so that their Lie
derivatives are given by the same formulas as in Section 3.2.

Using Ahlfors lemma and formulas \eqref{Ahlfors1}-\eqref{Ahlfors3}, we get
\begin{align*}
L_{\mu} \omega &= - \left(\left(\phi_{hyp}\right)_{\z} \dot{f}_{zz} +
\left(\phi_{hyp}\right)_{z}\left( \left(\phi_{hyp}\right)_z \dot{f}_{\z} +
\dot{f}_{z\z}\right)\right) dz\wedge d\z \\ &= \vartheta\mu\, dz \wedge d\z - d \xi,
\end{align*}
where
\begin{equation} \label{xi}
\xi = 2 \left(\phi_{hyp}\right)_z \dot{f}_{\z} d\z -\phi_{hyp}\, d\dot{f}_z.
\end{equation}
Since $\vartheta\mu$ is a $(1,1)$-tensor for $\Ga$, $\de (\vartheta\mu\,dz\wedge d\z)
= 0$, so that $\de L_{\mu} \omega = - \de d\xi$. We have
\begin{equation*}
\la d \xi, F_1-F_2 \ra = \la \xi, \pa' (F_1 - F_2) \ra = \la \xi, \pa'' (L_1 - L_2)
\ra = \la \de \xi, L_1-L_2 \ra.
\end{equation*}
Set $\chi=\de \xi + L_{\mu} \ptheta$. The 1-form $\chi$ on $\Omega$ is closed,
\begin{equation*}
d\chi = \de( d\xi) + L_{\mu} d \ptheta = \de(- L_{\mu} \omega) + L_{\mu} \de \omega =
0,
\end{equation*}
and satisfies
\begin{equation*}
\de \chi = \de (L_{\mu} \ptheta + \de \xi) = L_{\mu} \de \ptheta = L_\mu \pu.
\end{equation*}
Using \eqref{Ahlfors1}, \eqref{Ahlfors4}, \eqref{Ahlfors5} and part (ii) of Lemma
\ref{formula}, we get
\begin{align*}
L_{\mu} \ptheta_{\s^{-1}} =& -\dot{f}_z \frac{\s''}{\s'} \;dz +
\phi_{hyp}\left(\left(\dot{f}_{zz}\circ \s\; \s' - \dot{f}_{zz}\right)dz +
\frac{\s''}{\s'} \dot{f}_{\z} \;d\z\right) + \dot{f}_{z} \frac{\cl{\s''}}{\cl{\s'}}
\;d\z\\ &+\frac{1}{2} \left(-\left(\dot{f}_z \circ \s - \dot{f}_z\right)
\frac{\s''}{\s'}\; dz - \log |\s'|^2 \left(\dot{f}_{zz}\circ \s \;\s' -
\dot{f}_{zz}\right)dz \right. \\ & \left. - \log |\s'|^2 \frac{\s''}{\s'}
\dot{f}_{\z}\; d\z + \left(\dot{f}_{z}\circ \s - \dot{f}_z\right)
\frac{\cl{\s''}}{\cl{\s'}}\; d\z\right) \\ & - \left(\log|c(\g)|^2 + 2\log
2\right)d\left(\dot{f}_z \circ \g - \dot{f}_z\right) - \frac{\dot{c}(\s)}{c(\s)}
\left(\frac{\g''}{\g'} dz - \frac{\ov{\g''}}{\ov{\g'}} d\z\right) \\ =&  -\dot{f}_z
\frac{\s''}{\s'}dz + \dot{f}_z\circ\s \frac{\cl{\s''}}{\cl{\s'}}\; d\z -
d\left(\frac{1}{2}\log|\s^\prime|^2\left(\dot{f}_z \circ \s - \dot{f}_z\right) \right)
\\ & + \phi_{hyp}\,d\left(\dot{f}_z\circ\s -\dot{f}_z\right)
 - \left(\log|c(\g)|^2 + 2\log 2\right)d\left(\dot{f}_z \circ \g - \dot{f}_z\right) \\
& - \frac{\dot{c}(\s)}{c(\s)} \left(\frac{\g''}{\g'} dz - \frac{\ov{\g''}}{\ov{\g'}}
d\z\right).
\end{align*}
Using
\begin{equation*}
\de \xi_{\s^{-1}} = - 2 \frac{\s''}{\s'}\dot{f}_{\z}\; d\z
-\phi_{hyp}\,d\left(\dot{f}_z \circ \s - \dot{f}_z\right) +\log|\s^\prime|^2
d\left(\dot{f}_z\circ\s\right),
\end{equation*}
we get
\begin{align*}
\chi_{\s^{-1}} = & d\left(\frac{1}{2}\log|\s^\prime|^2 \left(\dot{f}_z \circ \s +
\dot{f}_z\right) -\left(\log|c(\s)|^2 + 2\log 2\right)\left(\dot{f}_z\circ \s -
\dot{f}_z\right)\right) \\ & - \left(\dot{f}_z \circ \s + \dot{f}_z\right)\frac{\s''}
{\s^\prime}dz - 2 \frac{\s''}{\s'}\dot{f}_{\z}\; d\z - \frac{\dot{c}(\s)}{c(\s)}
\left(\frac{\g''}{\g'} dz - \frac{\ov{\g''}}{\ov{\g'}} d\z\right).
\end{align*}
Using parts (ii) and (iii) of Lemma~\ref{formula} and
\begin{equation*}
-\frac{2\dot{c}}{cz+d}=\frac{\dot{c}}{c}\frac{\s''}{\s'}(z),
\end{equation*}
we finally obtain
\begin{align*}
\chi_{\s^{-1}} =\; & d\left( \frac{1}{2} \log|\s'|^2 \left(\dot{f}_z \circ \s +
\dot{f}_z + 2\frac{\dot{c}(\s)}{c(\s)}\right)\right. \\ & -\left.\left(\log|c(\s)|^2 +
2 + 2\log 2\right)\left(\dot{f}_z\circ \s - \dot{f}_z\right)\right) \\ =\; &
dl_{\s^{-1}}.
\end{align*}

We have
\begin{align*}
\la d\xi, F_1 - F_2\ra + \la L_\mu\ptheta, L_1 - L_2\ra & = \la \chi, L_1-L_2 \ra
 = \la d l, L_1 - L_2 \ra \\ & = \la l, \pa'(L_1 - L_2) \ra
 = \la l, \pa''(V_1 - V_2) \ra \\
& = \la \de l, V_1 - V_2 \ra.
\end{align*}
Using $L_\mu\pu = d\de l$ we get
\begin{equation*}
\la L_\mu\pu, W_1 - W_2 \ra  = \la \de l, \pa'(W_1 - W_2) \ra = \la \de l, V_1 - V_2
\ra
\end{equation*}
so that
\begin{equation*}
L_{\mu} S = \frac{i}{2} \la \vartheta\mu\, dz \wedge d\z, F_1 - F_2 \ra,
\end{equation*}
as asserted.

Finally, we justify the differentiation under the integral sign. Set
\begin{equation*}
l_{\s} = l^{(0)}_{\s} + l^{(1)}_{\s},
\end{equation*}
where
\begin{align*}
l^{(0)}_{\s^{-1}} & = \frac{\dot{c}(\s)}{c(\s)}\log|\s'|^2 - \left(\log|c(\s)|^2 +
2 + 2\log 2\right)(\dot{f}_z\circ \s - \dot{f}_z),\\ l^{(1)}_{\s^{-1}} & = \frac{1}{2}
\log|\s'|^2 \left(\dot{f}_z \circ \s + \dot{f}_z\right).
\end{align*}
Next, we use part (i) of Lemma \ref{formula}. According to it, the function
$l^{(0)}_{\s}$ is continuous on $\mathcal{C}\setminus\{\s(\infty)\}$. Since
\begin{equation*}
(\de l^{(1)})_{\s_1^{-1},\s_2^{-1}} = \frac{1}{2}\left(\log |\s_2'\circ
\s_1|^2(\dot{f}_z\circ\s_1-\dot{f}_z)- \log |\s_1'|^2(\dot{f}_z\circ \s_2\s_1 -
\dot{f}_z \circ \s_1)\right),
\end{equation*}
we also conclude that $(\de l^{(1)})_{\s_1,\s_2}$, and hence the function $(\de
l)_{\s_1,\s_2}$, are continuous on $\mathcal{C}
\setminus\{\s_1(\infty),\,(\s_1\s_2)(\infty)\}$. Now let $W_1^{(n)}\Subset W_1$ and
$W_2^{(n)}\Subset W_2$ be a sequence of 1-chains in $\Omega_1$ and $\Omega_2$ obtained
from $W_1$ and $W_2$ by ``cutting'' $\Ga$-contracting at 0 paths at points
$p'_n\in\Omega_1$ and $p''_n\in\Omega_2$, where $p_n',\,p''_n\rightarrow 0$ as
$n\rightarrow \infty$. Clearly,
\begin{equation*}
S= \lim_{n\rightarrow\infty}S_n,
\end{equation*}
where
\begin{equation*}
S_n = \frac{i}{2}\left(\la \omega, F_1 -F_2\ra -\la \ptheta, L_1 - L_2\ra + \la \pu,
W_1^{(n)} - W_2^{(n)}   \ra \right).
\end{equation*}
Our previous arguments show that
\begin{equation*}
L_\mu S_n = \frac{i}{2} \la \vartheta\mu\, dz \wedge d\z, F_1 - F_2 \ra - \la (\de
l)(p'_n), U_1 \ra + \la (\de l)(p''_n), U_2 \ra.
\end{equation*}
Since function $\de l$ is continuous at $p=0$ and $U_1=U_2$, we get
\begin{equation*}
\lim_{n\rightarrow\infty}L_\mu S_n = \frac{i}{2}\,\la \vartheta\mu, F_1 -F_2\ra.
\end{equation*}
Moreover, the convergence is uniform in some neighborhood of $\Ga$ in $\D(\Ga)$, since
$f^{\vep\mu}$ is holomorphic at $\vep=0$. Thus
\begin{equation*}
L_\mu S = \lim_{n\rightarrow\infty}L_\mu S_n,
\end{equation*}
which completes the proof.
\end{proof}
For fixed Riemann surface $Y$ denote by $P_F$ and $P_{QF}$ sections of
$\Pro(X)\rightarrow\T(X)$ corresponding to the Fuchsian uniformization of $X'\in\T(X)$
and to the simultaneous uniformization of $X'\in\T(X)$ and $Y$ respectively.
\begin{corollary} \label{f-qf} On the \Te space $\T(X)$,
\begin{equation*}
P_{F} - P_{QF} = \frac{1}{2}\,\pa S_Y.
\end{equation*}
\end{corollary}
\begin{remark} Conversely, Theorem \ref{variation1} follows from the Corollary
\ref{f-qf} and the symmetry property~\eqref{symmetry}.
\end{remark}
\begin{remark} In the Fuchsian case the maps $J_1$ and $J_2$ are identities and
similar computation shows that $\boldsymbol\vartheta=0$, in accordance with
$S=8\pi(2g-2)$ being a constant function on $\T(X)\times \T(\bar{X})$.
\end{remark}
\subsection{Second variation}
Here we compute $d\boldsymbol\vartheta = \bar\pa\boldsymbol\vartheta$. First, we have
the following statement.
\begin{lemma} \label{proj}
The quasi-Fuchsian projective connection $P_{QF}$ is a holomorphic section of the
affine bundle $\Pro(\Ga)\rightarrow \D(\Ga)$.
\end{lemma} \label{qf-hol}
\begin{proof}
Consider the following commutative diagram
\begin{equation*}
\begin{CD}
\Omega @>f^{\vep\mu}>> \Omega^{\vep\mu} \\ @VV\pi_{QF} V @VV\pi_{QF}^{\vep\mu} V
\\ X\sqcup Y @>F^{\vep\mu}>> X^{\vep\mu}\sqcup Y^{\vep\mu}
\end{CD}
\end{equation*}
where $\mu\in\Omega^{-1,1}(\Ga)$. We have
\begin{displaymath}
\mathcal{S}\left(\left(\pi_{QF}^{\vep\mu} \right)^{-1}\right)\circ F^{\vep\mu}
\left(F_z^{\vep\mu}\right)^2 + \mathcal{S}\left(F^{\vep\mu}\right) =
\mathcal{S}\left(f^{\vep\mu}\right)\circ \pi_{QF}^{-1}\left(\pi_{QF}^{-1} \right)_z^2
+ \mathcal{S}\left (\pi_{QF}^{-1} \right).
\end{displaymath}
Since $f^{\vep\mu}$ and, obviously $F^{\vep\mu}$, are holomorphic at $\vep=0$, we get
\begin{displaymath}
\left.\frac{\pa}{\pa\bar{\vep}}\right|_{\vep=0}
\mathcal{S}\left(\left(\pi_{QF}^{\vep\mu} \right)^{-1}\right) = 0.
\end{displaymath}
\end{proof}

Using Corollary \ref{f-qf}, Lemma \ref{qf-hol}  and the result~\cite{ZT87b}
\begin{equation*}
\bar\pa P_{F} = -i\,\omega_{WP},
\end{equation*}
which follows from \eqref{II} since $P_S$ is a holomorphic section of $\Pro_g
\rightarrow \mathfrak{S}_g$, we immediately get
\begin{corollary} \label{f-qf2} For fixed $Y$
\begin{equation*}
\pa \bar\pa S_Y = -2\bar\pa(P_F - P_{QF}) = -2 d(P_F - P_{QF})  = 2\,i\, \omega_{WP},
\end{equation*}
so that $-S_Y$ is a K\"{a}hler potential for the Weil-Petersson metric on $\T(X)$.
\end{corollary}

\begin{remark}
The equation $d(P_{F} - P_{QF}) = -i\,\omega_{WP}$ was first proved in~\cite{McM} and
was used for the proof that moduli spaces are K\"{a}hler hyperbolic (note that
symplectic form $\omega_{WP}$ used there is twice the one we are using here, and there
is a missing factor $1/2$ in the computation in~\cite{McM}). Specifically, the
Kraus-Nehari inequality
asserts that $P_F - P_{QF}$ is a bounded antiderivative of
$-\,i\,\omega_{WP}$ with respect to \Te and Weil-Petersson metrics~\cite{McM}. In this
regard, it is interesting to estimate the K\"{a}hler potential $S_Y$ on $\T(X)$. From
the basic inequality of the distortion theorem (see, e.g., \cite{Dur})
\begin{equation*}
\left|\frac{h^{\prime\prime}(z)}{h^\prime(z)} - \frac{2\z}{(1-|z|^2)}\right| \leq
\frac{4}{(1-|z|^2)},
\end{equation*}
where $h$ is a univalent function in the unit disk, we immediately get
\begin{displaymath}
|(\phi_{hyp})_z|^2\leq 4e^{\phi_{hyp}},
\end{displaymath}
so that the bulk term in $S_Y$ is bounded on $\T(X)$ by $20\pi(2g-2)$. It can also be
shown that other terms in $S_Y$ have at most ``linear growth'' on $\T(X)$, in
accordance with the boundness of $\pa S_Y$.
\end{remark}
The following result follows from the Corollary \ref{f-qf2} and the symmetry
property~\eqref{symmetry}. For completeness, we give its proof in the form that
is generalized verbatim to Kleinian groups.
\begin{theorem}\label{variation2} The following formula holds on $\D(\Ga)$,
\begin{equation*}
d \boldsymbol\vartheta = \bar\pa \pa S= -2i\,\omega_{WP},
\end{equation*}
so that $-S$ is a K\"{a}hler potential of the Weil-Petersson metric on $\D(\Ga)$.
\end{theorem}
\begin{proof}
Let $\mu, \nu \in \Omega^{-1,1} (\Gamma)$. First, using Cartan formula, we get
\begin{align*}
d \boldsymbol\vartheta \left(\frac{\pa}{\pa \vep_{\mu}}, \frac{\pa}{\pa
\vep_{\nu}}\right)  = & L_{\mu}(\boldsymbol\vartheta (\nu)) - L_{\nu}
(\boldsymbol\vartheta(\mu)) - \boldsymbol\vartheta\left(\left[\frac{\pa}{\pa
\vep_{\mu}}, \frac{\pa}{\pa \vep_{\nu}}\right]\right) \\
  = & L_{\mu}(L_{\nu} S) - L_{\nu}(L_{\mu} S) = 0,
\end{align*}
which just manifests that $\pa^2 = 0$. On the other hand,
\begin{align*}
d \boldsymbol\vartheta \left(\frac{\pa}{\pa \vep_{\mu}}, \frac{\pa}{\pa
\bar{\vep}_{\nu}}\right) = & L_{\mu}(\boldsymbol\vartheta (\bar{\nu})) - L_{\bar{\nu}}
(\boldsymbol\vartheta(\mu)) - \boldsymbol\vartheta\left(\left[\frac{\pa}{\pa
\vep_{\mu}}, \frac{\pa}{\pa \bar\vep_\nu}\right] \right) \\ = & - L_{\bar{\nu}}
\iint\limits_{\Ga\bk\Omega} \vartheta \mu \\ = & -
\iint\limits_{\Ga\bk\Omega} (L_{\bar{\nu}} \boldsymbol\vartheta) \mu,
\end{align*}
since $\boldsymbol\vartheta$ is a $(1,0)$-form.

The computation of $L_{\bar{\nu}} \boldsymbol\vartheta$ repeats verbatim the one given
in \cite{ZT87b}. Namely, consider the commutative diagram \eqref{cd} with $i=1,2$,
and, for brevity, omit the index $i$. Since $ (J^{\vep \nu})^{-1} \circ f^{\vep \nu} =
F^{\vep \hat{\nu}} \circ J^{-1}$, the property \textbf{SD1} of the Schwarzian
derivative (applicable when at least one of the functions is holomorphic) yields
\begin{align}\label{schwarzian}
\mathcal{S}(J^{\vep\nu})^{-1} \circ f^{\vep \nu} \; (f_z^{\vep \nu})^2 +
\mathcal{S}(f^{\vep \mu}) = \mathcal{S}(F^{\vep\hat{\nu}}) \circ J^{-1} (J_z^{-1})^2 +
\mathcal{S}(J^{-1}).
\end{align}
We obtain
\begin{align*}
\left.\frac{\pa}{\pa \bar{\vep}_{\nu}} \right|_{\vep=0} \mathcal{S}(J^{\vep
\nu})^{-1}
\circ f^{\vep \nu} \; (f_z^{\vep \nu})^2 = & \left.\frac{\pa}{\pa
\bar{\vep}_{\nu}}\right|_{\vep=0} \mathcal{S}(F^{\vep {\hat{\nu}}}) \circ J^{-1}
(J_z^{-1})^2 \\ = & \left.\frac{\pa}{\pa \bar{\vep}_{\nu}} \right|_{\vep=0}
F_{zzz}^{\vep \hat{\nu}} \circ J^{-1} (J_z^{-1})^2 \\ = & -\frac{1}{2}\rho
\ov{\nu(z)},
\end{align*}
where in the last line we have used Ahlfors formula \eqref{A2}. Finally,
\begin{align*}
d \boldsymbol\vartheta \left(\frac{\pa}{\pa \vep_{\mu}}, \frac{\pa}{\pa
\bar{\vep}_{\nu}}\right) = \iint\limits_{\Ga\bk\Omega} \mu \bar{\nu} \rho
 =  -2i\,\omega_{WP} \left(\frac{\pa}{\pa \vep_{\mu}}, \frac{\pa}{\pa \bar{\vep}_{\nu}}\right).
\end{align*}
\end{proof}
\subsection{Quasi-Fuchsian reciprocity}
The existence of the function $S$ on the deformation space $\D(\Ga)$ satisfying the
statement of Theorem \ref{variation1} is a global form of quasi-Fuchsian reciprocity.
Quasi--Fuchsian reciprocity of McMullen \cite{McM} follows from it as immediate
corollary.

Let $\mu, \nu \in \Omega^{-1,1}(\Gamma)$ be such that $\mu$ vanishes outside
$\Omega_1$ and $\nu$ --- outside $\Omega_2$, so that Lie derivatives $L_{\mu}$ and
$L_\nu$ stand for the variation of $X$ for fixed $Y$ and variation of $Y$ for fixed
$X$ respectively.

\begin{theorem} (McMullen's quasi-Fuchsian reciprocity)
\begin{equation*}
\iint_{X} (L_{\nu} \mathcal{S}(J_1^{-1}))\;\mu = \iint_{Y} (L_{\mu}
\mathcal{S}(J_2^{-1}))\; \nu.
\end{equation*}
\end{theorem}
\begin{proof} Immediately follows from Theorem \ref{variation1}, since
\begin{align*}
L_{\nu}L_{\mu} S &=  2\iint_{X} (L_{\nu} \mathcal{S}(J_1^{-1}))\; \mu, \\
L_{\mu}L_{\nu} S &=  2\iint_{Y} (L_{\mu} \mathcal{S}(J_2^{-1}))\; \nu, \\
\end{align*}
and $[L_\mu, L_\nu]=0$.
\end{proof}

In \cite{McM}, quasi-Fuchsian reciprocity was used to prove that $d(P_F - P_{QF}) =
-i\,\omega_{WP}$. For completeness, we give here another proof of this result using
earlier approach in~\cite{ZT87a}, which admits generalization to other deformation
spaces.

\begin{proposition} On the deformation space $\D(\Ga)$,
\begin{equation*}
\pa\boldsymbol\vartheta=0.
\end{equation*}
\end{proposition}
\begin{proof}
Using the same identity ~\eqref{schwarzian} which follows from the commutative diagram
\eqref{cd}, we have
\begin{align*}
\left.\frac{\pa}{\pa \vep_{\nu}} \right|_{\vep=0} \mathcal{S}(J^{\vep\nu})^{-1} \circ
f^{\vep \nu} \; (f_z^{\vep\nu})^2 = & \left. \frac{\pa}{\pa \vep_{\nu}}
\right|_{\vep=0} \mathcal{S}(F^{\vep \hat{\nu}}) \circ J^{-1} (J_z^{-1})^2 -
\left.\frac{\pa}{\pa \vep_{\nu}} \right|_{\vep=0} \mathcal{S}(f^{\vep \nu})\\ = &
\left.\frac{\pa}{\pa \vep_{\nu}} \right|_{\vep=0} F_{zzz}^{\vep \hat{\nu}} \circ
J^{-1} (J_z^{-1})^2 - \left.\frac{\pa}{\pa \vep_{\nu}} \right|_{\vep=0} f_{zzz}^{\vep
\nu},\\
\end{align*}
where we replaced $\mu$ by $\nu$ and omit index $i=1,2$. Differentiating
\eqref{qc-var} three times with respect to $z$ we get
\begin{equation} \label{zzz}
\left.\frac{\pa}{\pa \vep_{\nu}} \right|_{\vep=0} f_{zzz}^{\vep \nu} (z) =
-\frac{6}{\pi} \iint\limits_{\CC} \frac{\nu(w)}{(z-w)^4} \;\; d^2w = -\frac{6}{\pi}
\iint\limits_{\Ga\bk\Omega} K(z,w) \nu(w) d^2w,
\end{equation}
where
\begin{equation*}
K(z,w) = \sum_{\s \in \Gamma} \frac{\s'(w)^2}{(z-\s w)^4}.
\end{equation*}
It is well-known that for harmonic $\nu$ the integral in \eqref{zzz} is understood in
the principal value sense (as $\lim_{\delta\rightarrow 0}$ of integral over
$\CC\setminus\{|w-z|\leq \delta\}$). Therefore, using Ahlfors formula \eqref{A1} we
obtain
\begin{align*}
(L_{\nu} \boldsymbol\vartheta)(z) = \frac{12}{\pi} \iint\limits_{\Ga\bk\Omega} K(z,w)
\nu(w) d^2w,
\end{align*}
and
\begin{align*}
\pa \boldsymbol\vartheta (\mu, \nu) &= L_{\mu} \boldsymbol\vartheta(\nu) - L_{\nu}
\boldsymbol\vartheta(\mu)\\
 &=\iint\limits_{\Ga\bk\Omega}(L_{\mu} \boldsymbol\vartheta)(z) \, \nu(z) d^2z -
 \iint\limits_{\Ga\bk\Omega} (L_{\nu} \boldsymbol\vartheta)(w) \, \mu(w) d^2w
= 0,
\end{align*}
since kernel $K(z,w)$ is obviously symmetric in $z$ and $w$, $K(z,w)=K(w,z)$.
\end{proof}

\section{Holography}
Let $\Gamma$ be marked, normalized, purely loxodromic quasi-Fuchsian group of genus
$g>1$. The group $\Gamma\subset \PSL(2,\CC)$ acts on the closure $\ov{\up}^3=\up^3\cup
\hat{\CC}$ of the hyperbolic 3-space $\up^3 = \{Z = (x,y,t) \in \RR^3\,|\, t>0 \}$.
The action is discontinuous on $\up^3 \cup \Omega$ and $M = \Gamma \bk (\up^3 \cup
\Omega)$ is a hyperbolic 3-manifold, compact in the relative topology of $\ov{\up}^3$,
with the boundary $X \sqcup Y \simeq\Gamma\bk\Omega$. According to the holography
principle, the on-shell gravity theory on $M$, given by the Einstein-Hilbert action
functional with the cosmological term, is equivalent to the ``off-shell'' gravity
theory on its boundary $X \sqcup Y$, given by the Liouville action functional. Here we
give a precise mathematical formulation of this principle.
\subsection{Homology and cohomology set-up} We start by generalizing
homological algebra methods in Section 2 to the three-dimensional case.
\subsubsection{Homology computation} Denote by $\SSS_{\bu} \equiv \SSS_{\bu}(\up^3
\cup\Omega)$ the standard singular chain complex of $\up^3 \cup\Omega$, and let $R$ be
a fundamental region of $\Gamma$ in $\up^3 \cup\Omega$ such that $R \cap \Omega$ is
the fundamental domain $F=F_1-F_2$ for the group $\Ga$ in $\Omega$ (see Section 2). To
have a better picture, consider first the case when $\Gamma$ is a Fuchsian group. Then
$R$ is a region in $\ov{\up}^3$ bounded by the hemispheres which intersect $\hat{\C}$
along the circles that are orthogonal to $\RR$ and bound the fundamental domain $F$
(see Section 2.2.1). The fundamental region $R$ is a three-dimensional $CW$-complex
with a single 3-cell given by the interior of $R$. The 2-cells --- the  faces $D_k,
D_k', E_k$ and $E_k'$, $k=1, \dots, g$, are given by the parts of the boundary of $R$
bounded by the intersections of the hemispheres and the arcs $a_k - \bar{a}_k,
a_k'-\bar{a}_k', b_k-\bar{b}_k$ and $b_k'-\bar{b}_k'$ respectively (see Fig.~1). The
1-cells --- the edges, are given by the 1-cells of $F_1-F_2$ and by $e_k^0, e^1_k,
f_k^0, f^1_k$ and $d_k$, $k=1, \dots, g$, defined as follows. The edges $e^0_k$ are
intersections of the faces $E_{k-1}$ and $D_k$ joining the vertices $\bar{a}_k(0)$ to
$a_k(0)$, the edges $e_k^1$ are intersections of the faces $D_k$ and $E_k'$ joining
the vertices $\bar{a}_k(1)$ to $a_k(1)$; $f^0_k=e^0_{k+1}$ are intersections of $E_k$
and $D_{k+1}$ joining $\bar{b}_k(0)$ to $b_k(0)$, $f_k^1$ are intersections of $D_k'$
and $E_k$ joining $\bar{b}_k(1)$ to $b_k(1)$, and $d_k$ are intersections of $E'_k$
and $D'_k$ joining $\bar{a}'_k(1)$ to $a'_k(1)$. Finally, the 0-cells --- 
the vertices, are
given by the vertices of $F$. This property means that the edges of $R$ do not
intersect in $\up^3$. When $\Ga$ is a quasi-Fuchsian group, the fundamental region $R$
is a topological polyhedron homeomorphic to the geodesic polyhedron for the
corresponding Fuchsian group $\tilde\Ga$.

As in the two-dimensional case, we construct the 3-chain representing $M$ in the total
complex $\Tot\KKK$ of the double homology complex $\KKK_{\bu,\bu} =
\SSS_\bu\otimes_{\Z\Ga}\BBB_\bu$ as follows. First, identify $R$ with $R\otimes [\,]
\in \KKK_{3,0}$. We have $\pa'' R=0$ and
\begin{align*}
\pa' R & = - F + \sum_{k=1}^g \left(D_k - D_k' - E_k + E_k'\right)\\
       & = - F + \pa''S,
\end{align*}
where $S \in \KKK_{2,1}$ is given by
\begin{equation*}
S = \sum_{k=1}^g \left(E_k\otimes[\be_k] - D_k\otimes[\al_k]\right).
\end{equation*}
Secondly,
\begin{align*}
\pa' S  & =\sum_{k=1}^g \left((b_k - \bar{b}_k)\otimes[\be_k] - (a_k - \bar{a}_k)
\otimes[\al_k]\right) \\
       &- \sum_{k=1}^g \left(\left(f_k^1 - f_k^0\right) \otimes[\be_k] - \left(e_k^1 -
       e_k^0\right) \otimes[\al_k]\right)\\
      &= L - \pa''E,
\end{align*}
where $L= L_1 - L_2$ and $E\in \KKK_{1,2}$ is given by
\begin{align*}
E &= \sum_{k=1}^g \left(e_k^0\otimes[\al_k|\be_k] - f_k^0\otimes[\be_k|\al_k] +
f_k^0\otimes[\g_k^{-1}|\al_k\be_k]\right)\\ & \;\;\; -\sum_{k=1}^{g-1}
f_g^0\otimes[\g_g^{-1}\ldots\g_{k+1}^{-1}|\g_k^{-1}].
\end{align*}
Therefore $\pa' E = V = V_1 - V_2$ and the 3-chain $R-S+E \in (\Tot\KKK)_3$ satisfies
\begin{equation} \label{R-Sigma}
\pa (R - S + E) = -F - L + V = -\Sigma,
\end{equation}
as asserted.

\subsubsection{Cohomology computation}
The $\PSL(2,\CC)$-action on $\up^3$ is the following. Represent $Z=(z,t)\in\up^3$ by a
quaternion
\begin{equation*}
\boldsymbol{Z}= x\cdot\boldsymbol{1} + y\cdot\boldsymbol{i} + t\cdot\boldsymbol{j} =
\begin{pmatrix} z & -t \\ t & \z \end{pmatrix},
\end{equation*}
and for every $c\in\CC$ set
\begin{equation*}
\boldsymbol{c}=\re c\cdot\boldsymbol{1} + \im c \cdot\boldsymbol{i} =
\begin{pmatrix} c & 0 \\ 0 & \bar{c} \end{pmatrix}.
\end{equation*}
Then for $\s = \ma{a}{b}{c}{d} \in\PSL(2, \C)$ the action $Z \mapsto\s Z$ is given by
\begin{equation*}
\boldsymbol{Z}  \mapsto (\boldsymbol{a}\boldsymbol{Z} + \boldsymbol{b})
(\boldsymbol{c}\boldsymbol{Z} + \boldsymbol{d})^{-1}.
\end{equation*}

Explicitly, for $Z=(z,t)\in\up^3$ setting $z(Z)=z$ and $t(Z)=t$ gives
\begin{align}
z(\s Z) & = \left((az+b)(\ov{cz+d}) + a\bar{c}\,t^2 \right) J_\s(Z), \label{z}
\\ t(\s Z) & = t\,J_\s(Z), \label{t}
\end{align}
where
\begin{equation*}
J_\s(Z) = \frac{1}{|cz+d|^2 + |ct|^2}.
\end{equation*}
Note that $J_\s^{3/2}(Z)$ is the Jacobian of the map $ Z \mapsto \s Z$, hence it
satisfies the transformation property
\begin{equation} \label{jacobian}
J_{\s_1 \circ \s_2}(Z) = J_{\s_1}(\s_2 Z) J_{\s_2}(Z).
\end{equation}
From \eqref{z} and \eqref{t} we get the following formulas for the derivatives
\begin{align}
\frac{\pa z(\s Z)}{\pa z} & = (\ov{cz+d})^2  J_{\s}^2(Z), \label{d-z} \\ \frac{\pa
z(\s Z)}{\pa\z} &= - (\bar{c}\,t)^2 J_{\s}^2(Z), \label{d-barz} \\ \frac{\pa z(\s
Z)}{\pa t} &= 2t\bar{c}(\ov{cz+d}) J_{\s}^2(Z) \label{d-t}.
\end{align}
In particular,
\begin{equation} \label{as t to 0}
\frac{\pa z(Z)}{\pa z} = \s'(z) + O(t^2),\quad \frac{\pa z(\s Z)}{\pa \z}
=O(t^2),\quad \frac{\pa z(Z)}{\pa t} = O(t),
\end{equation}
as $t\rightarrow 0$ and $z\in\hat{\CC}\setminus\{\s^{-1}(\infty)\}$, where for
$z\in\CC$ we continue to use the two-dimensional notations
\begin{equation*}
\s(z) = \frac{az+b}{cz+d}\quad \text{and} \quad \s'(z) = \frac{1}{(cz+d)^2},\quad
\frac{\s''}{\s'}(z) = \frac{-2c}{cz+d}.
\end{equation*}

The hyperbolic metric on $\up^3$ is given by
\begin{equation*}
ds^2=\frac{|dz|^2 + dt^2}{t^2},
\end{equation*}
and is $\PSL(2,\CC)$-invariant. Denote by
\begin{equation*}
w_3 = \frac{1}{t^3}\,dx\wedge dy\wedge dt = \frac{i}{2t^3}\,dz\wedge d\z \wedge dt
\end{equation*}
the corresponding volume form on $\up^3$. The form $w_3$ is exact on $\up^3$,
\begin{equation} \label{w-2}
w_3=dw_2,\quad \text{where} \quad w_2= -\frac{i}{4t^2}\,dz\wedge d\z.
\end{equation}
The 2-form $w_2 \in \CCC^{2,0}$ is no longer  $\PSL(2,\CC)$-invariant. A
straightforward computation using \eqref{d-z}-\eqref{d-t} gives for $\s =
\ma{a}{b}{c}{d}\in\PSL(2,\CC)$,
\begin{align*}
(\de w_2)_{\s^{-1}} = & \s^* w_2 - w_2 \\ = & \frac{i}{2}\,J_{\s}(Z) \left(|c|^2
dz\wedge d\z - \frac{c(\ov{cz+d})}{t}\,dz\wedge dt +
                    \frac{\bar{c}(cz+d)}{t}\,d\z \wedge dt \right).
\end{align*}
Since $d\de w_2 = \de dw_2 = \de w_3 = 0$ and $\up^3$ is simply connected, this
implies that there exists $w_1 \in \CCC^{1,1}$ such that $dw_1 = \de w_2$. Explicitly,
\begin{equation} \label{w-1}
(w_1)_{\s^{-1}} = -\frac{i}{8}\,\log \left(|ct|^2 J_{\s}(Z)\right)
\left(\frac{\s''}{\s'} dz  - \frac{\ov{\s''}}{\ov{\s'}}d\z  \right).
\end{equation}
Using \eqref{jacobian} and \eqref{as t to 0} we get for $\de w_1 \in C^{1,2}$
\begin{align} \label{delta w-1}
(\de w_1)_{\s_1^{-1}, \s_2^{-1}} =& - \frac{i}{8}\left(\log J_{\s_1}(Z) + \log
\frac{|c(\s_2)|^2}{|c(\s_2\s_1)|^2}\right) \left( \frac{\s_2''}{\s_2'}\circ \s_1
\s_1'\,dz - \frac{\ov{\s_2''}}{\ov{\s_2'}}\circ \s_1 \ov{\s'_1}\,d\z\right)\nonumber
\\ & - \frac{i}{8}\left(\log J_{\s_2}(\s_1 Z) + \log
\frac{|c(\s_2\s_1)|^2}{|c(\s_1)|^2}\right) \left(\frac{\ov{\s_1''}}{\ov{\s_1'}}\,d\z -
\frac{\s_1''}{\s_1'}\,dz\right) \nonumber \\ & + B_{\s_1^{-1}, \s_2^{-1}}(Z).
\end{align}
Here $B_{\s_1^{-1}, \s_2^{-1}}(Z) = O(t\log t)$ as $t \rightarrow 0$, uniformly on
compact subsets of $\CC\setminus\{\s_1^{-1}(\infty),(\s_2\s_1)^{-1}(\infty)\}$.

Clearly the 1-form $\de w_1$ is closed,
\begin{equation*}
d (\de w_1) = \de (dw_1) = \de (\de w_2) = 0.
\end{equation*}
Since $\up^3$ is simply connected, there exists $w_0\in \CCC^{0,2}$ such that $w_1 =
dw_0$. Moreover, using $H^3(\Gamma, \C)=0$ we can always choose the antiderivative
$w_0$ such that $\de w_0 = 0$. Finally, set $\Phi = w_2 - w_1 - w_0 \in (\Tot\CCC)^2$,
so that
\begin{equation} \label{D-Phi}
D\Phi=w_3.
\end{equation}
\subsection{Regularized Einstein-Hilbert action}
In two dimensions, the critical value of the Liouville action for a Riemann surface
$X\simeq \Ga\bk\up$ is proportional to the hyperbolic area of the surface (see Section
2). It is expected that in three dimensions the critical value of the Einstein-Hilbert
action functional with cosmological term is proportional to the hyperbolic volume of
the 3-manifold $M\simeq \Ga\bk(\up^3\cup\Omega)$ (plus a term proportional to the
induced area of the boundary). However, the hyperbolic metric diverges at the boundary
of $\ov{\up}^3$ and for quasi-Fuchsian group $\Ga$ (as well as for general Kleinian
group \footnote{Note that we are using definition of  Kleinian groups as in
\cite{Maskit}. In the theory of hyperbolic 3-manifolds these groups are called
Kleinian groups of the second kind.}) the hyperbolic volume of
$\Ga\bk(\up^3\cup\Omega)$ is infinite. In \cite{Witten}, Witten proposed a
regularization of the action functional by truncating the 3-manifold $M$ by surface
$f=\vep$, where the cut-off function $f \in C^{\infty} (\up^3,\RR_{>0})$ vanishes to
the first order on the boundary of $\ov{\up}^3$. Every choice of the function $f$
defines a metric  on $\up^3$
\begin{equation*}
ds^2 = \frac{f^2}{t^2}(|dz|^2+dt^2),
\end{equation*}
belonging to the conformal class of the hyperbolic metric. On the boundary of
$\ov{\up}^3$ it induces the metric
\begin{displaymath}
\lim_{t\rightarrow 0} \frac{f^2(z,t)}{t^2}|dz|^2.
\end{displaymath}

Clearly for the case of quasi-Fuchsian group $\Gamma$ (or for the general Kleinian
case considered in the next section), the cut-off function $f$ should be
$\Ga$-automorphic. Existence of such function is guaranteed by the following result,
which we formulate for the general Kleinian case.

\begin{lemma} \label{partition}
Let $\Gamma$ be non-elementary purely loxodromic, geometrically finite
 Kleinian group with region of discontinuity $\Omega$, normalized so that
$\infty\notin\Omega$. For every $\phi \in \mathcal{CM}(\Gamma\bk\Omega)$ there
exists $\Gamma$-automorphic function $f \in C^{\infty}(\up^3\cup\Omega)$ which is
positive on $\up^3$ and satisfies
\begin{equation*}
f(Z) = t e^{\phi(z)/2} + O(t^3),\quad \text{as}~t\rightarrow 0,
\end{equation*}
uniformly on compact subsets of $\Omega$.
\end{lemma}

\begin{proof}
Note that $\Ga\bk\Omega$ is isomorphic to a finite disjoint union of
compact Riemann surfaces. Let $R$ be a fundamental region of $\Ga$ in
$\up^3\cup\Omega$ which is compact in the relative topology of $\ov{\up}^3$.
I.~Kra has proved in \cite{Kra} (the construction in \cite{Kra} suggested
by M.~Kuga generalizes verbatim to our case) that there exist a bounded 
open set $V$ in $\ov{\up}^3$ such that $R\subset V$ and a function 
$\eta \in C^{\infty}(\up^3 \cup \Omega)$ ---
partition of unity for $\Gamma$ on $\up^3 \cup\Delta$, satisfying 
the following properties.
\begin{itemize}
\item[(i)]
$0 \leq \eta \leq 1$ and $\text{supp}\,\eta\subset \ov{V}$.
\item[(ii)]
For each $Z\in \up^3 \cup \Omega$ there is a neighborhood $U$ of $Z$ and a finite
subset $J$ of $\Gamma$ such that $\eta\vert_{\s(U)} = 0$ for each $\s \in
\Gamma\setminus J$.
\item[(iii)]
$\sum_{\s \in \Gamma} \eta(\s Z) = 1$ for all $Z\in \up^3 \cup \Omega$.
\end{itemize}

Let $B = \ov{V} \cap \{(z,t) \;\vert\; z \in \Lambda\}$. Since $R\cap\Omega$ is 
compact, then (shrinking $V$ if necessary) there exists a $t_0 > 0$ 
such that $B$ does not intersect the region $ \{(z,t)\in\ov{V}\;\vert \; t 
\leq t_0\}$. Define the function $\hat{f} : \ov{V} \rightarrow \RR$ by
\begin{equation*}
\hat{f}(z,t) = \begin{cases} 
               t e^{\phi(z)/2} & \text{if}\;\; (z,t) \in \ov{V} \;\; \text{and} 
\;t\leq t_0/2, \\
1 & \text{if} \;\; (z,t) \in \ov{V} \;\; \text{and} \;\;t\geq t_0,
\end{cases}
\end{equation*}
and extend it to smooth function $\hat{f}$ on $\ov{V}$, positive on $\ov{V}\cap\up^3$.
Set
\begin{equation*}
f(Z) = \sum_{\g \in \Gamma} \eta (\g Z) \hat{f}(\g Z).
\end{equation*}
By the property (ii), for every $Z\in \up^3 \cup \Omega$ this sum contains only
finitely many non-zero terms, so that the function $f$ is well-defined. By properties
(i) and (iii) it is positive on $\up^3$. To prove the asymptotic behavior, we use
elementary formulas
\begin{align*}
z(\g Z) &= \frac{az+b}{cz+d} + O(t^2) = \g(z) + O(t^2),\\ t(\g Z) &=
\frac{t}{|cz+d|^2} + O(t^3)\quad \text{as}~t \rightarrow 0,
\end{align*}
where $z \neq \g^{-1}(\infty)$. Since $\phi$ is smooth on $\Omega$ and
\begin{displaymath}
e^{\phi(\g z)/2} = e^{\phi(z)/2} |cz+d|^2,
\end{displaymath}
we get for $z\in \Omega$ such that $\s Z\in\ov{V}$ and 
$t$ is small enough,
\begin{align*}
\hat{f}(\g Z) &= \left(\frac{t}{|cz+d|^2} + O(t^3)\right)\left(e^{\phi(\g z)/2} +
O(t^2) \right)\\ &= t e^{\phi(z)/2} + O(t^3),
\end{align*}
where the O-term depends on $\g$. Using properties (ii) and (iii) we finally obtain
\begin{align*}
f(Z) &= \sum_{\g \in \Gamma} \eta(\g Z)\left(t e^{\phi(z)/2} + O(t^3)\right)\\
           &= t e^{\phi(z)/2} + O(t^3),
\end{align*}
uniformly on compact subsets of $\Omega$.
\end{proof}

Returning to the case when $\Gamma$ is a normalized purely loxodromic quasi-Fuchsian
group, for every $\phi \in \mathcal{CM}(\Gamma \bk \Omega)$
let $f$ be a function given by the lemma. For $\vep >0$ let $R_{\vep} = R \cap \{f
\geq \vep\}$ be the truncated fundamental region. For every chain $c$ in $\up^3$ let
$c_\vep = c\cap\{f\geq \vep\}$ be the corresponding truncated chain. Also let
$F_{\vep} = \pa' R_{\ep} \cap \{f = \vep\}$ be the boundary of $R_{\vep}$ on the
surface $f=\vep$ and define chains $L_{\vep}$ and  $V_{\vep}$ on $f=\vep$ by the same
equations $\pa' F_\vep =\pa''L_\vep$ and $\pa'L_\vep = \pa'' V_\vep$ as chains $L$ and
$V$ (see Sections 2.2.1 and 2.3.1). Since the truncation is $\Gamma$-invariant, for
every chain $c \in \SSS_{\bu}(\up^3)$ and $\g \in \Gamma$ we have
\begin{displaymath}
(\g c)_{\vep} = \g c_{\vep}.
\end{displaymath}
In particular, relations between the chains, derived in Section 5.1, hold for
truncated chains as well.

Let $M_{\vep}$ be the truncated 3-manifold with the  boundary $\pa' M_{\vep}$. For
$\vep$ sufficiently small $\pa' M_{\vep}=X_{\vep}\sqcup Y_{\vep}$ is diffeomorphic to
$X\sqcup Y$. Denote by $V_{\vep}[\phi]$ the hyperbolic volume of $M_{\vep}$. The
hyperbolic metric induces a metric on $\pa' M_{\vep}$, and $A_{\vep}[\phi]$ denotes the
area of $\pa' M_{\vep}$ in the induced metric.
\begin{definition}
The regularized on-shell Einstein-Hilbert action functional is defined by
\begin{align*}
\mathcal{E}_{\Gamma}[\phi] = -4\lim_{\vep \rightarrow 0} \left(V_{\vep}[\phi]
-\frac{1}{2} A_{\vep}[\phi] - 2\pi\chi(X)\log\vep\right),
\end{align*}
where $\chi(X)=\chi(Y)=2-2g$ is the Euler characteristic of $X$.
\end{definition}
The main result of this section is the following.
\begin{theorem} (Quasi-Fuchsian holography) \label{holography} For every
$\phi\in\mathcal{CM}(\Ga\bk\Omega)$ the regularized Einstein-Hilbert action is
well-defined and
\begin{displaymath}
\mathcal{E}_{\Gamma}[\phi] = \check{S}_{\Gamma}[\phi],
\end{displaymath}
where $\check{S}_{\Gamma}[\phi]$ is the modified Liouville action functional without
the area term,
\begin{equation*}
\check{S}_{\Gamma}[\phi] = S_{\Gamma}[\phi] - \iint\limits_{\Ga\bk\Omega}\,e^{\phi}d^2
z - 8\pi\,(2g-2)\log 2.
\end{equation*}
\end{theorem}
\begin{proof}
It is sufficient to verify the formula,
\begin{equation}\label{regqcvol}
 V_{\vep}[\phi] - \frac{1}{2} A_{\vep}[\phi] = 2\pi\chi(X)\log\vep - \frac{1}{4}
 \check{S}_{\Gamma}[\phi] + o(1)\quad \text{as}~\vep \rightarrow 0,
\end{equation}
which is a counter-part of the formula \eqref{Krasnov} for quasi-Fuchsian groups.

The area form induced by the hyperbolic metric on the surface $f(Z)=\vep$ is given by
\begin{displaymath}
\sqrt{1+\left(\frac{f_x}{f_t}\right)^2 + \left(\frac{f_y}{f_t}\right)^2} \frac{dx
\wedge dy} {t^2}.
\end{displaymath}
Using
\begin{displaymath}
\frac{f_x}{f_t}(Z) = \frac{t}{2}\,\phi_x(z) + O(t^3)\quad \text{and} \quad
\frac{f_y}{f_t}(Z) = \frac{t}{2}\,\phi_y(z) + O(t^3),
\end{displaymath}
we have as $\vep\rightarrow 0$
\begin{align*}
A_{\vep}[\phi] & = \iint\limits_{F_{\vep}}  \sqrt{1+ \frac{t^2}{4} (\phi_x^2 +
\phi_y^2)(z) + O(t^4)}\, \frac{dx\wedge dy}{t^2} \\ & = \iint\limits_{F_{\vep}}
\frac{dx\wedge dy}{t^2} + \frac{1}{2}\iint\limits_{F} \phi_z \phi_{\z}\,dx \wedge dy +
o(1)\\ & = \iint\limits_{F_\vep} \frac{dx\wedge dy}{t^2} + \frac{i}{4} \la
\check\omega[\phi], F \ra + o(1).
\end{align*}
Here we have introduced
\begin{equation} \label{omega-prime}
\check\omega[\phi] = \omega[\phi] - e^{\phi} dz \wedge d\z =|\phi_z|^2 dz \wedge d\z,
\end{equation}
and have used that for $Z \in F_{\vep}$
\begin{equation}\label{boundary}
t = \vep e^{-\phi(z)/2} + O(\vep^3),
\end{equation}
uniformly for $Z=(z,t)$ where $z\in F$.

Next, using \eqref{R-Sigma} and \eqref{D-Phi} we have,
\begin{align*}
V_{\vep}[\phi] & = \;\;\, \left\la w_3, R_{\vep} \right\ra \\ & = \;\;\, \left\la w_3,
R_{\vep} - S_{\vep} + E_{\vep} \right\ra \\
        & = \;\;\, \left\la D(w_2 - w_1 - w_0), R_{\vep} - S_{\vep} + E_{\vep} \right\ra\\
        & = \;\;\,  \left\la w_2 - w_1 - w_0, \pa \left(R_{\vep} - S_{\vep} + E_{\vep}\right)
        \right\ra \\
        & = - \la w_2, F_{\vep} \ra + \la w_1, L_{\vep} \ra - \la w_0, V_{\vep} \ra.
\end{align*}
The terms in this formula simplify as $\vep\rightarrow 0$. First of all, it follows
from \eqref{w-2} that
\begin{align*}
- \la w_2, F_{\vep} \ra = \frac{1}{2}\iint\limits_{F_{\vep}} \frac{dx \wedge dy}{t^2}.
\end{align*}
Secondly, using \eqref{boundary} and $J_{\s}(Z)=|\s'(z)| + O(t^2)$ as $t\rightarrow
0$, we have on $L_{\vep}$
\begin{align*}
(w_1)_{\s^{-1}}& = - \frac{i}{8}\log \left(\left|c\vep\right|^2
e^{-\phi}\,|\s'(z)|\right) \left(\frac{\s''}{\s'} dz  - \frac{\ov{\s''}}{\ov{\s'}} d\z
\right) + o(1) \\ & =  -\frac{i}{8} \left(2\log\vep - \phi + \frac{1}{2} \log |\s'|^2
+ \log |c(\s)|^2\right) \left(\frac{\s''}{\s'}dz - \frac{\ov{\s''}}{\ov{\s'}} d\z
\right) + o(1).
\end{align*}
Therefore, as $\vep\rightarrow 0$,
\begin{align*}
\la w_1, L_{\vep} \ra = -\frac{i}{4}\la \varkappa, L \ra (\log \vep - \log 2)
 + \frac{i}{8} \la \ptheta[\phi], L\ra + o(1),
\end{align*}
where 1-forms $\varkappa_\s$ and  $\ptheta_\s[\phi]$ were introduced in Corollary
\ref{varkappa-lemma} and formula \eqref{theta-prime} respectively. Finally,
\begin{displaymath}
\la w_0, V_{\vep}\ra = \la w_0, \pa' E_{\vep} \ra = \la dw_0, E_{\vep}\ra = \la \de
w_1, E_{\vep}\ra = \la \de w_1, E\ra + o(1),
\end{displaymath}
where we used that 1-form $\de w_1$ is smooth on $\up^3$ and continuous on
$\CC\setminus \Ga(\infty)$. Since it is closed, we can replace the 1-chain $E$ by the
1-chain $W=W_1-W_2$ consisting of $\Ga$-contracting paths at 0 (see Section 2.3). It
follows from \eqref{delta w-1} that $\de w_1 =\frac{i}{8} \pu + o(1)$  as
$t\rightarrow 0$, where the 1-form $\pu_{\s_1,\s_2}$ was introduced in
\eqref{u-prime}, so that
\begin{displaymath}
- \la w_0, V_{\vep} \ra = - \frac{i}{8} \la \pu, W \ra + o(1).
\end{displaymath}

Putting everything together, we have as $\vep\rightarrow 0$
\begin{align*}
V_{\vep}[\phi] - \frac{1}{2}A_{\vep}[\phi] = & -\frac{i}{4}\la \varkappa, L \ra ( \log
\vep - \log 2) - \frac{i}{8}\left( \la \check\omega[\phi], F \ra - \la \ptheta[\phi],
L \ra \right.\\ & \left.+ \la \pu, W \ra \right) + o(1).
\end{align*}

Using Corollary \ref{varkappa-lemma}, trivially modified for the quasi-Fuchsian case,
and \eqref{shift-3} concludes the proof.
\end{proof}

A fundamental domain $F$ for $\Gamma$ in $\Omega$ is called admissible, if it is the
boundary in $\C$ of a fundamental region $R$ for $\Ga$ in $\up^3\cup\Omega$. As an
immediate consequence of the theorem we get the following.
\begin{corollary}
The Liouville action functional $S_{\Gamma}[\phi]$ is independent of the choice of 
admissible fundamental domain.
\end{corollary}
\begin{proof}
Since $V_{\vep}[\phi]$, $A_{\vep}[\phi]$ are intrinsically associated with the
quotient manifolds $M\simeq \Ga\bk(\up^3\cup\Omega)$ and $X\sqcup Y\simeq
\Ga\bk\Omega$, the statement follows from the definition of the Einstein-Hilbert
action and the theorem.
\end{proof}

Although we proved the same result in Section 2 using methods of homological algebra,
the above argument easily generalizes to other Kleinian groups.

\begin{remark}
The truncation of the 3-manifold $M$ by the function $f$ does depend on the choice of
the realization of the fundamental group of $M$ as a normalized discrete subgroup
$\Ga$ of $\PSL(2, \C)$. Different realizations of $\pi_1(M)$ result in different
choices of the function $f$, since $f$ has to satisfy the asymptotic behavior in Lemma
\ref{partition}, where the leading term $te^{\phi(z)/2}$ is not a well-defined
function on $M$.
\end{remark}
\begin{remark}
The cochain $w_0 \in \CCC^{0,2}$ was defined as a solution of the equation $dw_0=w_1$
satisfying $\de w_0=0$. However, in the computation in Theorem \ref{holography} this
condition is not needed --- any choice of an antiderivative for $w_1$ will suffice.
This is due to the fact that the chain in $(\Tot\KKK)_3$ that starts with 
$R\in\KKK_{3,0}$
does not contain a term in $\KKK_{0,3}$, hence $\pa' E = V$. Thus we can trivially add
the term $\la\de w_0, R_{\vep} - S_{\vep} + E_{\vep}\ra=0$ to $V_{\vep}[\phi]$, which
through the equation $D\Phi=w_3-\de w_0$ still gives $\la w_0, V \ra = \la dw_0, E
\ra$. Thus the absence of $\KKK_{0,3}$-components in the chain in $(\Tot\KKK)_3$
implies that each term in $E$ produces two boundary terms in $V$ which cancel out the
integration constants in definition of $w_0$. As a result, $S_{\Ga}[\phi]$ does not
depend on the choice of $w_0$. In the next section we generalize the Liouville action
functional to Kleinian groups having the same property.
\end{remark}
\subsection{Epstein map}
Construction of the regularized Einstein-Hilbert action in the previous section
works for a larger class of cut-off surfaces than those given by equation
$f=\vep$. Namely, it follows from the proof of Theorem \ref{holography}, 
that the statement holds for any family $S_{\vep}$ of cut-off surfaces such that for
$Z=(z,t)\in S_{\vep}$
\begin{equation} \label{boundary2}
t =  \vep e^{-\phi(z)/2} + O(\vep^3)\quad \text{as}~\vep\rightarrow 0,
\end{equation}
uniformly for $z\in F$. 

Given a conformal metric $ds^2=e^{\phi(z)}|dz|^2$ on $\Omega\subset
\hat{\CC}$ there is a natural surface in $\up^3$ associated 
to it through the inverse of the hyperbolic Gauss map. Corresponding construction 
is due to C.~Epstein~\cite{Epstein1,Epstein2} (see also~\cite{Anderson}) and is 
the following.
For every $Z \in\up^3$ and $z\in\hat{\CC}$ there is a unique horosphere $H$ based at
point $z$ and passing through the point $Z$: $H$ is a Euclidean sphere in $\up^3$ 
tangent to $z\in\CC$ and passing through $Z$, or is an Euclidean plane parallel to 
the complex plane for $z=\infty$. Denote by $[Z,z]$ an affine parameter of the 
horosphere $H$ --- the hyperbolic distance 
between the point $(0,1)\in\up^3$ and the horosphere $H$ 
considered as positive if the point $(0,1)$ is outside $H$
and negative otherwise. Denote the corresponding horosphere by $H(z,[Z,z])$. 
The Epstein map
$\mathcal{G}: \Omega \rightarrow \up^3$ is defined by
\begin{align*}
e^{\phi(z)/2} |dz| = e^{[\mathcal{G}(z), z]} \frac{2|dz|}{1+|z|^2} 
\end{align*} 
and it is $\Ga$-invariant
\begin{equation*}
\mathcal{G}\circ\s =\s\circ \mathcal{G}~\text{for all}~\s\in\Ga.
\end{equation*}
if $\phi\in\mathcal{CM}(\Ga\bk\Omega)$.
\begin{remark}
Note that our definition of the Epstein map corresponds to the case $f=id$ in 
the Definition 3.9 in \cite{Anderson}.
\end{remark}

Geometrically, the image of the Epstein map is the Epstein surface $\mathcal{H}
=\mathcal{G}(\Omega)$, which is the envelope of the family of horospheres
$H(z, \varrho(z))$ with
\begin{displaymath}
\varrho(z) = \log \left(\frac{1}{2}(1+ |z|^2)\right) + \frac{\phi(z)}{2},
\end{displaymath}
parametrized by $z\in\Omega$, where 
$\mathcal{G}(z)$ is the point of tangency 
of the horosphere $H(z, \varrho(z))$ with the surface $\mathcal{H}$. 
Explicit computation gives
\begin{equation}\label{InverseGauss}
\mathcal{G}(w) = \left( w + \frac{2 \phi_{\bar{w}}(w)}{e^{\phi(w)}+  
\vert \phi_w(w) \vert^2},\; \frac{2 e^{\phi(w)/2}}{e^{\phi(w)}+ 
\vert \phi_w(w) \vert^2} \right), \;\; w\in \Omega.
\end{equation}

\begin{remark}
The square of Euclidean distance between points $w$ and $\mathcal{G}(w)$
in $\ov{\up}^3$ is $4/(e^{\phi(w)} + |\phi_w(w)|^2)$. This gives a geometric
interpretation of the density $|\phi_z|^2 + e^\phi$ of the (1,1)-form 
$\omega$ \eqref{omega}. 
\end{remark}

Now to a given $\phi\in\mathcal{CM}(\Gamma\bk \Omega)$ we associate the family
$\phi_{\vep}=\phi + 2\log 2 -2\log\vep\in \mathcal{CM}(\Gamma\bk\Omega)$ with
$\vep>0$, which corresponds to the family of
conformal metrics $ds^2_{\vep}=4\vep^{-2} e^{\phi(w)}|dw|^2$, 
and consider the corresponding
$\Ga$-invariant family $\mathcal{H}_{\vep}$ of Epstein surfaces. It follows 
from the parametric representation
\begin{align}\label{parametric}
z & = w + \frac{2 \vep^2 \phi_{\bar{w}}(w)}
{4 e^{\phi(w)}+ \vep^2 \vert \phi_w (w)\vert^2},\\
t & = \frac{4\vep e^{\phi(w)/2}}{4e^{\phi(w)} + \vep^2 
\vert \phi_w (w)\vert^2},
\end{align}
that for $\vep$ small the surfaces $\mathcal{H}_{\vep}$ embed smoothly in 
$\up^3$ and as $\vep \rightarrow 0$, 
\begin{align*}
z &= w + O(\vep^2),\\
t &= \vep e^{-\phi(w)/2} + O(\vep^3),
\end{align*}
uniformly for $w$ in compact subsets of $\Omega$. These formulas 
immediately give the desired asymptotic behavior \eqref{boundary2}. 
The choice of Epstein surfaces 
$\mathcal{H}_{\vep}$ as cut-off surfaces for definition of the regularized 
Einstein-Hilbert action seems to be the most natural. It is quite remarkable 
that independently equations \eqref{parametric} appear in \cite{Krasnov2} in
relation with a general solution of ``asymptotically AdS three-dimensional 
gravity''.  

\section{Generalization to Kleinian groups}
\subsection{Kleinian groups of Class $A$}
Let $\Ga$ be a finitely generated Kleinian  group with the region of discontinuity
$\Omega$, a maximal set of non-equivalent components $\Omega_1,\dots,\Omega_n$ of
$\Omega$, and the limit set $\Lambda=\hat\C\setminus\Omega$. As in the quasi-Fuchsian
case, a path $P$ is called $\Ga$-contracting in $\Omega$, if $P=P_1\cup P_2$, where
$p\in\Lambda \setminus\{\infty\}$ is a fixed point for $\Ga$, paths
$P_1\setminus\{p\}$ and $P_2\setminus\{p\}$ lie entirely in distinct components of
$\Omega$ and are $\Ga$-contracting at $p$ in the sense of Definition
\ref{contracting1}. It follows from arguments in Section 2.3.1 that $\Ga$-contracting
paths in $\Omega$ are rectifiable.
\begin{definition}
A Kleinian group $\Ga$ is of Class $A$ if it satisfies the following conditions.
\begin{enumerate}
\item[\textbf{A1}]
$\Gamma$ is non-elementary and purely loxodromic.
\item[\textbf{A2}] $\Gamma$ is geometrically finite.
\item[\textbf{A3}] $\Ga$ has a fundamental region $R$ in $\up^3\cup\Omega$ which is a
finite three-dimensional $CW$-complex with no 0-dimensional cells in $\up^3$ and such
that $R\cap\Omega\subset\Omega_1\cup\dots\cup\Omega_n$.
\end{enumerate}
\end{definition}

In particular, property \textbf{A1} implies that $\Ga$ is torsion-free and does not
contain parabolic elements, and property \textbf{A2} asserts that $\Ga$ has a
fundamental region $R$ in $\up^3\cup\Omega$ which is a finite topological polyhedron.
Property \textbf{A3} means that the region $R$ can be chosen such that the vertices of
$R$ --- endpoints of edges of $R$, lie on $\Omega\in\hat{\C}$ and the boundary of $R$
in $\hat{\C}$, which is a fundamental domain for $\Ga$ in $\Omega$, is not too
``exotic''.

The class $A$ is rather large: it clearly contains all purely loxodromic Schottky
groups (for which the property \textbf{A3} is vacuous), Fuchsian groups,
quasi-Fuchsian groups, and free combinations of these groups.

As in the previous section, we say that Kleinian group $\Ga$ is normalized if
$\infty\in\Lambda$.

\subsection{Einstein-Hilbert and Liouville functionals}
For a finitely generated Kleinian group $\Gamma$ let $M\simeq \Gamma \bk(\up^3\cup
\Omega)$ be corresponding hyperbolic 3-manifold, and let $\Gamma_1, \dots, \Gamma_n$
be the stabilizer groups of the maximal set $\Omega_1, \dots, \Omega_n$  of
non-equivalent components of $\Omega$. We have
\begin{align*}
\Gamma \bk \Omega = \Gamma_1\bk \Omega_1 \sqcup \dots \sqcup \Gamma_n \bk \Omega_n
\simeq X_1 \sqcup \dots \sqcup X_n,
\end{align*}
so that Riemann surfaces $X_1, \dots, X_n$ are simultaneously uniformized by $\Gamma$.
Manifold $M$ is compact in the relative topology of $\ov{\up}^3$ with the disjoint
union $X_1\sqcup\dots\sqcup X_n$ as the boundary.

\subsubsection{Homology and cohomology set-up} Let $\SSS_{\bu}\equiv
\SSS_{\bu}(\up^3 \cup \Omega)$, $\BBB_{\bu} \equiv \BBB_{\bu}(\Z \Gamma)$ be standard
singular chain and bar-resolution homology complexes and $\KKK_{\bu, \bu} \equiv
\SSS_{\bu} \otimes_{\Z \Gamma}\BBB_{\bu}$ --- the corresponding double complex. 
When $\Ga$
is a Kleinian group of Class $A$, we can generalize homology construction from the
previous section and define corresponding chains $R, S, E, F, L, V$ in total complex
$\Tot \KKK$ as follows. Let $R$ be a fundamental region for $\Ga$ in $\up^3 \sqcup
\Omega$
--- a closed topological polyhedron in $\ov{\up}^3$
satisfying property \textbf{A3}. The group $\Ga$ is generated by side pairing
transformations of $R\cap\up^3$ and we define the chain $S\in\KKK_{2,1}$ as the sum of
terms $-s \otimes \g^{-1}$ for each pair of sides $s, s'$ of $R\cap\up^3$ identified
by a transformation $\g$, i.e., $s'=-\g s$. The sides are oriented as components of
the boundary and the negative sign stands for the opposite orientation. We have
\begin{equation} \label{R}
\pa' R = -F + \pa'' S,
\end{equation}
where $F = \pa'R\cap\Omega\in \KKK_{2,0}$. Note that it is immaterial whether we
choose $-s \otimes \g^{-1}$ or $-s' \otimes \g$ in the definition of $S$, since these
terms differ by a $\pa''$-coboundary. Next, relations between generators of $\Ga$
determine the $\Ga$-action on the edges of $R$, which, in turn, determines the chain
$E\in \KKK_{1,2}$ through the equation
\begin{equation} \label{S}
\pa' S = L - \pa'' E.
\end{equation}
Here $L = \pa' S\cap\Omega\in \KKK_{1,1}$. Finally, property \textbf{A3} implies that
\begin{equation} \label{E}
\pa' E = V,
\end{equation}
where the chain $V \in \KKK_{0,2}$ lies in $\Omega$.

Next, let the 1-chain $W \in \KKK_{1,2}$ be a ``proper projection'' of the 1-chain $E$
onto $\Omega$, i.e., $W$ is defined by connecting every two vertices belonging to the
same edge of $R$ either by a smooth path lying entirely in one component of $\Omega$,
or by a $\Ga$-contracting path, so that $\pa' W = V$. The existence of such 1-chain
$W$ is guaranteed by the property \textbf{A3} and the following lemma, which is of
independent interest.
\begin{lemma} \label{chains} Let $\Ga$ be a normalized, geometrically finite, purely
loxodromic Kleinian group, and let $R$ be the fundamental region of $\Ga$ in
$\up^3\cup\Omega$ such that $R\cap\Omega\subset\Omega_1\cup\dots\cup\Omega_n$
--- a union of the maximal set of non-equivalent components of $\Omega$. If an
edge $e$ of $R\cap\up^3$ has endpoints $v_0$ and $v_1$ belonging to two distinct
components $\Omega_i$ and $\Omega_j$, then there exists a $\Ga$-contracting path in
$\Omega$ joining vertices $v_0$ and $v_1$. In particular, $\Omega_i$ and $\Omega_j$
has at least one common boundary point, which is a fixed point for $\Ga$.
\end{lemma}
\begin{proof} There exist sides $s_1$ and $s_2$ of $R$ such that $e\subset s_1\cap
s_2$. For each of these sides there exists a group element identifying it with another
side of $R$. Let $\s\in\Ga$ be such element for $s_1$. Since $\Ga$ is torsion-free and
$v_0, v_1\in\Omega$, element $\s$ identifies the edge $e$ with the edge
$e'$ of $R$ with endpoints $\s(v_0)\neq v_0$ and $\s(v_1)\neq v_1$. Since, by
assumption, $R\cap\Omega\subset \Omega_1\cup\dots\cup\Omega_n$,  we have that
$\s(v_0)\in\Omega_i$ and $\s$ fixes $\Omega_i$. Similarly, $\s(v_1)\in\Omega_j$ and
$\s$ fixes $\Omega_j$. Now assume that attracting fixed point $p$ of $\s$ is not
$\infty$ (otherwise we replace $\s$ by $\s^{-1}$). Join $v_0$ and $\s(v_0)$ by a
smooth path $P_1^0$ inside $\Omega_i$, and let $P_1^n=\s^n(P_1^0)$ be its $n$-th
$\s$-iterate. Since $\s$ fixes $\Omega_i$, the path $P_1^n$ lies entirely inside
$\Omega_i$. Since $\lim_{n\rightarrow\infty}\s^n(v_0)= p$, the path
$P_1=\cup_{n=0}^\infty P_1^n$ joins $v_0$ and $p$, and except for the endpoint $p$ lies
entirely in $\Omega_i$. Clearly path $P_1^0$ can be chosen so that the path $P_1$ is
smooth everywhere except at $p$. The path $P_2$ joining points $v_1$ and $p$ inside
$\Omega_j$ is defined similarly, and the path $P=P_1\cup P_2$ is $\Ga$-contracting in
$\Omega$.
\end{proof}
Setting $\Sigma=F+L-V$ we get from \eqref{R}-\eqref{E} that
\begin{displaymath}
\pa(R-S+E) = -\Sigma.
\end{displaymath}

\begin{remark} Since $\up^3$ is acyclic, it follows from general arguments in
\cite{AT} that for any geometrically finite purely loxodromic Kleinian group $\Ga$
with fundamental region $R$ given by a closed topological polyhedron, there exist
chains $S\in\KKK_{2,1}, E\in\KKK_{1,2}, T\in\KKK_{0,3}$ and chains $F\in\KKK_{2,0},
L\in\KKK_{1,1}, V\in\KKK_{0,2}$ on $\Omega$, satisfying
\begin{align*}
\pa' R &= -F + \pa'' S\\ \pa' S &= L - \pa'' E\\ \pa' E &= V + \pa''T.
\end{align*}
Property \textbf{A3} asserts that $T=0$, and we get equations \eqref{R}-\eqref{E}.
\end{remark}

Correspondingly, let $\AAA^{\bu} \equiv \AAA_{\C}^{\bu}(\up^3 \cup \Omega)$ and
$\CCC^{\bu, \bu} \equiv \Hom(\BBB_{\bu}, \AAA^{\bu})$ be the de Rham complex on $\up^3
\cup \Omega$ and the bar-de Rham complex respectively. The cochains $w_3, w_2, w_1,
\de w_1, w_0$ are defined by the same formulas as in Section 5.1. For $\phi \in
\mathcal{CM}(\Gamma \bk \Omega)$ define the cochains $\omega[\phi], \theta[\phi], u$
by the same formulas \eqref{omega}, \eqref{theta}, \eqref{u}, with the group elements
belonging to $\Gamma$. Finally, define the cochains $\ptheta[\phi], \pu$ by
\eqref{theta-prime} and \eqref{u-prime}.

\subsubsection{Action functionals}
Let $\Ga$ be a normalized Class $A$ Kleinian group. For each $\phi \in
\mathcal{CM}(\Gamma \bk \Omega)$ let $f$ be the function constructed in Lemma
\ref{partition}. As in Section 5.2, we truncate the manifold $M$ by the cut-off
function $f$ and define $V_{\vep}[\phi]$, $A_{\vep}[\phi]$.
\begin{definition}
The regularized on-shell Einstein-Hilbert action functional for a normalized Class $A$
Kleinian group $\Ga$ is defined by
\begin{align*}
\mathcal{E}_{\Gamma}[\phi] = -4\lim_{\vep \rightarrow 0}\left(V_{\ep}[\phi]
-\frac{1}{2} A_{\ep}[\phi] - \pi(\chi(X_1)+\dots+\chi(X_n))\log \vep\right).
\end{align*}
\end{definition}

As in the quasi-Fuchsian case, a fundamental domain $F$ for a Kleinian group $\Gamma$
in $\Omega$ is called admissible, if it is the boundary in $\C$ of a fundamental
region $R$ for $\Ga$ in $\up^3$ satisfying property \textbf{A3}.

\begin{definition}
The Liouville action functional $S_{\Gamma}: \mathcal{CM} (\Gamma \bk \Omega)
\rightarrow \RR$ for a normalized Class $A$ Kleinian group $\Gamma$ is defined by
\begin{equation}
S_{\Gamma} [\phi] = \frac{i}{2} \left(\la \omega[\phi], F \ra - \la \ptheta[\phi], L
\ra + \la \pu, W\ra\right),
\end{equation}
where $F$ is an admissible fundamental domain for $\Ga$ in $\Omega$.
\end{definition}

\begin{remark}When $\Gamma$ is a purely loxodromic Schottky group (not necessarily
classical Schottky group), the Liouville action functional defined above is, up to the
constant term $4\pi(2g-2)\log 2$, the functional \eqref{Schottky}, introduced by P.
Zograf and the first author \cite{ZT87b}.
\end{remark}

Using these definitions and repeating verbatim arguments in Section 5 we have the
following result.
\begin{theorem} \label{K-holography} (Kleinian holography) For every
$\phi\in\mathcal{CM}(\Ga\bk\Omega)$ the regularized Einstein-Hilbert action is
well-defined and
\begin{displaymath}
\mathcal{E}_{\Gamma}[\phi] = \check{S}_{\Gamma}[\phi] = S_{\Gamma}[\phi] -
\iint\limits_{\Ga\bk\Omega} e^\phi d^2 z + 4\pi\left(\chi(X_1)+\dots+\chi(X_n)
\right)\log 2.
\end{displaymath}
\end{theorem}
\begin{corollary} \label{well-defined}
The definition of a Liouville action functional does not depend on the choice of
admissible fundamental domain $F$ for $\Gamma$.
\end{corollary}

As in the Fuchsian and quasi-Fuchsian cases, the Euler-Lagrange equation for the
functional $S_{\Gamma}$ is the Liouville equation, and its single critical point given
by the hyperbolic metric $e^{\phi_{hyp}} |dz|^2$ on $\Ga\bk\Omega$ is non-degenerate.
For every component $\Omega_i$ denote by $J_i : \U \rightarrow \Omega_i$ the
corresponding covering map (unique up to a $\PSL(2, \RR)$-action on $\up$). Then the
density $e^{\phi_{hyp}}$ of the hyperbolic metric is given by
\begin{equation} \label{hyp-kleinian}
e^{\phi_{hyp}(z)} = \frac{|(J_i^{-1})'(z)|^2}{(\im \,J_i^{-1}(z))^2} \quad \text{if}
\;\; z \in \Omega_i \;, \;\; i=1,\dots, n.
\end{equation}
\begin{remark} \label{Polyakov-2}
As in Remark \ref{Polyakov}, let $\Delta[\phi]=-e^{-\phi}\pa_z \pa_{\z}$ be the
Laplace operator of the metric $ds^2=e^\phi |dz|^2$ acting on functions on
$X_1\sqcup\dots\sqcup X_n$, let $\det\Delta[\phi]$ be its zeta-function regularized
determinant, and let
\begin{equation*}
\mathcal{I}[\phi]=\log\frac{\det\Delta[\phi]}{A[\phi]}.
\end{equation*}
Polyakov's ``conformal anomaly'' formula and Theorem \ref{K-holography} give the
following relation between Einstein-Hilbert action $\mathcal{E}[\phi]$ for $M\simeq
\Ga\bk(\up^3\cup\Omega)$ and ``analytic torsion'' $\mathcal{I}[\phi]$ on its boundary
$X_1\sqcup\dots\sqcup X_n\simeq \Ga\bk\Omega$,
\begin{equation*}
\mathcal{I}[\phi + \sigma] + \frac{1}{12\pi} \mathcal{E}[\phi + \sigma]=
\mathcal{I}[\phi] +  \frac{1}{12\pi} \mathcal{E}[\phi],\;\; \sigma\in
C^\infty(X_1\sqcup\dots\sqcup X_n,\RR).
\end{equation*}
\end{remark}
\subsection{Variation of the classical action}
Here we generalize theorems in Section 4 for quasi-Fuchsian groups to Kleinian groups.
\subsubsection{Classical action}
Let $\Gamma$ be a normalized Class $A$ Kleinian group and let $\D(\Gamma)$ be its
deformation space. For every Beltrami coefficient $\mu \in \mathcal{B}^{-1,1}(\Gamma)$
the normalized quasiconformal map $f^{\mu} : \hat{\C} \rightarrow \hat{\C}$ descends
to an orientation preserving homeomorphism of the quotient Riemann surfaces $\Gamma
\bk \Omega$ and $\Gamma^{\mu} \bk \Omega^{\mu}$. This homeomorphism extends to
homeomorphism of corresponding $3$-manifolds $\Gamma \bk (\up^3 \cup \Omega)$ and
$\Gamma^{\mu} \bk (\up^3 \cup \Omega^{\mu})$, which can be lifted to orientation
preserving homeomorphism of $\up^3$. In particular, a fundamental region of $\Gamma$
is mapped to a fundamental region of $\Gamma^{\mu}$. Hence property \textbf{A3} is
stable and every group in $\D(\Gamma)$ is of Class $A$. Moreover, since $\infty$
is a fixed point of $f^{\mu}$, every group in $\D(\Gamma)$ is normalized.

For every $\Gamma' \in\D(\Gamma)$ let $S_{\Gamma'} = S_{\Gamma'}[\phi'_{hyp}]$ be the
classical Liouville action for $\Gamma'$. Since the property of the fundamental domain
$F$ being admissible is stable, Corollary \ref{well-defined} asserts that the
classical action gives rise to a well-defined real-analytic function $S:
\D(\Gamma)\rightarrow\RR$.

As in Section 4, let $\vartheta \in \Omega^{2,0}(\Gamma)$ be the holomorphic quadratic
differential $\Ga\bk\Omega$, defined by
\begin{equation*}
\vartheta = 2 (\phi_{hyp})_{zz} - (\phi_{hyp})_z^2.
\end{equation*}
It follows from \eqref{hyp-kleinian} that
\begin{align*}
\vartheta(z) &= 2 \mathcal{S}(J_i^{-1})(z)\quad\text{if}\;\;z\in\Omega_i,\;\;
i=1,\dots,n.
\end{align*}
Define a $(1,0)$-form $\boldsymbol{\vartheta}$ on $\D(\Gamma)$ by assigning to every
$\Gamma' \in \D(\Gamma)$ corresponding $\vartheta \in\Omega^{2,0}(\Gamma')$.

For every $\Ga'\in\D(\Ga)$ let $P_{F}$ and $P_K$ be Fuchsian and Kleinian projective
connections on $X'_1\sqcup\dots\sqcup X'_n \simeq \Ga'\bk\Omega'$, defined by the
Fuchsian uniformizations of Riemann surfaces $X'_1,\dots, X'_n$ and by their
simultaneous uniformization by Kleinian group $\Ga'$. We will continue to denote
corresponding sections of the affine bundle $\Pro(\Ga)\rightarrow\D(\Ga)$ by $P_F$ and
$P_K$ respectively. The difference $P_F - P_K$ is a $(1,0)$-form on $\D(\Ga)$. As in
the Section 4.1,
\begin{equation*}
\boldsymbol{\vartheta} = 2(P_F - P_K).
\end{equation*}
Correspondingly, the isomorphism
\begin{displaymath}
\D(\Gamma) \simeq \D(\Gamma_1, \Omega_1) \times \dots \times \D(\Gamma_n, \Omega_n)
\end{displaymath}
defines embeddings
\begin{equation*}
\D(\Gamma_i, \Omega_i) \hookrightarrow \D(\Gamma)
\end{equation*}
and pull-backs $S_i$ and $(P_F - P_K)_i$ of the function $S$ and the $(1,0)$-form $P_F
- P_K$. The deformation space $\D(\Gamma_i, \Omega_i)$ describes simultaneous Kleinian
uniformization of Riemann surfaces $X_1, \cdots, X_n$ by varying the complex structure
on $X_i$ and keeping the complex structures on other Riemann surfaces fixed, and the
$(1,0)$-form $(P_F - P_K)_i$ is the difference of corresponding projective
connections.
\subsubsection{First variation} Here we compute the $(1,0)$-form $\pa S$ on $\D(\Ga)$.
\begin{theorem}   \label{K-global}
On the deformation space $\D(\Gamma)$,
\begin{displaymath}
\pa S = 2(P_F - P_K).
\end{displaymath}
\end{theorem}
\begin{proof}
Since $F^{\vep \mu}= f^{\vep \mu}(F)$ is an admissible fundamental domain for
$\Gamma^{\vep\mu}$, and, according to Lemma \ref{Paths}, the 1-chain $W^{\vep\mu}=
f^{\vep \mu}(W)$ consists of $\Ga^{\vep\mu}$-contracting paths in $\Omega^{\vep\mu}$,
the proof repeats verbatim the proof of Theorem \ref{variation1}. Namely, after the
change of variables we get
\begin{equation*}
L_{\mu}S = \frac{i}{2} \left(\la L_{\mu} \omega, F \ra - \la L_{\mu} \ptheta, L \ra +
\la L_\mu,\pu,W \ra\right),
\end{equation*}
where
\begin{align*}
L_{\mu} \omega &= \vartheta \mu dz\wedge d\z - d \xi
\end{align*}
and 1-form $\xi$ is given by \eqref{xi}. As in the proof of Theorem \ref{variation1},
setting $\chi=\de \xi + L_{\mu} \ptheta$ we get that the 1-form $\chi$ on $\Omega$ is
closed,
\begin{equation*}
d\chi = \de( d\xi) + L_{\mu} d \ptheta = \de(- L_{\mu} \omega) + L_{\mu} \de \omega =
0,
\end{equation*}
and satisfies
\begin{equation*}
\de \chi = \de (L_{\mu} \ptheta + \de \xi) = L_{\mu} \de \ptheta = L_\mu \pu =d\de l.
\end{equation*}
Since the $1$-chain $W$ consists either of smooth paths or of $\Ga$-contracting paths 
in $\Omega$, and function
$\de l$ is continuous on $W$, the same arguments as in the proof of Theorem
\ref{variation1} allow to conclude that
\begin{align*}
L_{\mu}S &= \frac{i}{2} \la\vartheta \mu dz \wedge d\z, F \ra.
\end{align*}
\end{proof}
\begin{corollary} Let $X_1,\dots, X_n$ be Riemann surfaces simultaneously
uniformized by a Kleinian group $\Ga$ of Class $A$. Then on $\D(\Omega_i,\Ga_i)$
\begin{equation*}
(P_F - P_K)_i = \frac{1}{2}\,\pa S_i.
\end{equation*}
\end{corollary}
\subsubsection{Second variation}
\begin{theorem} \label{K-global-2}
On the deformation space $\D(\Gamma)$,
\begin{align*}
d \boldsymbol{\vartheta} = \bar\pa \pa S = -2i\,\omega_{WP},
\end{align*}
so that $-S$ is a \Ka potential of the Weil-Petersson metric on $\D(\Gamma)$.
\end{theorem}
The proof is the same as the proof of Theorem \ref{variation2}.

\subsection{Kleinian Reciprocity}

Let $\mu \in \Omega^{-1,1}(\Gamma)$ be a harmonic Beltrami differential, $f^{\vep\mu}$
be corresponding normalized solution of the Beltrami equation, and let
$\upsilon=\dot{f}$ be corresponding vector field on $\hat{\C}$,
\begin{align*}
\upsilon (z) = - \frac{1}{\pi} \iint\limits_{\C} \frac{\mu(w) z(z-1)}{(w-z)w(w-1)}
d^2w
\end{align*}
(see Section 3.2). Then
\begin{displaymath}
\varphi_{\mu}(z)= \frac{\pa^3}{\pa z^3} \upsilon (z) = - \frac{6}{\pi}
\iint\limits_{\C} \frac{\mu(w) }{(w-z)^4}d^2w
\end{displaymath}
is a quadratic differential on $\Gamma \bk \Omega$, holomorphic outside the support of
$\mu$.

In \cite{McM} McMullen proposed the following generalization of quasi-Fuchsian
reciprocity.
\begin{theorem}(McMullen's Kleinian Reciprocity) \label{McM-2}
Let $\Gamma$ be a finitely generated Kleinian group. Then for every $\mu, \nu \in
\Omega^{-1,1}(\Gamma)$
\begin{align*}
\iint\limits_{\Gamma\bk \Omega} \varphi_{\mu} \nu = \iint \limits_{\Gamma\bk
\Omega}\varphi_{\nu} \mu.
\end{align*}
\end{theorem}
The proof in \cite{McM} is based on the symmetry of the kernel $K(z,w)$, defined in
Section 4.2. Here we note that Theorem \ref{K-global} provides a global form of
Kleinian reciprocity for Class $A$ groups from which Theorem \ref{McM-2} follows
immediately.

Indeed, when $\Gamma$ is a normalized Class $A$ Kleinian group, Kleinian reciprocity
is the statement
\begin{displaymath}
L_{\mu} L_{\nu} S =  L_{\nu} L_{\mu} S,
\end{displaymath}
since, according to \eqref{zzz},
\begin{displaymath}
-\frac{1}{2}L_{\mu} \boldsymbol{\vartheta} (z)  = - \frac{6}{\pi} \iint\limits_{\C}
\frac{\mu(w) }{(w-z)^4}d^2w = \varphi_{\mu}(z)
\end{displaymath}
and
\begin{displaymath}
\iint\limits_{\Ga\bk\Omega}\varphi_\mu\nu = -\frac{1}{2}
\iint\limits_{\Ga\bk\Omega}L_\mu L_\nu S.
\end{displaymath}

\bibliographystyle{amsalpha}
\bibliography{action}

\end{document}